\documentclass[a4paper,11pt,oneside]{article}
\usepackage[utf8]{inputenc}
\usepackage[T1]{fontenc}
\usepackage{amsmath,amsthm,amsfonts,amssymb,amscd}
\usepackage{a4wide}
\usepackage[algoruled]{algorithm2e}
\usepackage[capitalize]{cleveref}
\usepackage{mathtools}
\usepackage{booktabs}
\usepackage{pgfplots}
\pgfplotsset{compat=1.14,grid style={gray!15}}
\usepackage{siunitx}[=v2]
\usepackage{ao-math-std}
\usepackage{ao-math-fields}
\usepackage{xcolor}
\usepackage{authblk}

\newcommand\qfor{\quad\forall}
\newcommand\tfor{\ \forall}
\let\bs\boldsymbol
\newcommand\boldu{\bs u}
\newcommand\boldv{\bs v}
\newcommand\boldvh{\bs v_h}
\newcommand\boldz{\bs z}
\newcommand\bolduu{\bs{uu}}
\newcommand\bolduz{\bs{uz}}
\newcommand\boldzu{\bs{zu}}
\newcommand\boldPhi{\bs\Phi}
\newcommand\bolddeltau{\bs{\delta u}}
\newcommand\bolddeltaz{\bs{\delta z}}
\newcommand\Phiu[1][]{\Phi_{u\optindex{#1}}}
\newcommand\Phiuy{\Phi_{u:y}}
\newcommand\Phiphi[1][]{\Phi_{\varphi\optindex{#1}}}

\newcommand\Phiphih[1][]{\Phi_{\varphi,h\optindex{#1}}}
\newcommand\zu[1][]{z_{u\optindex{#1}}}

\newcommand\zphi[1][]{z_{\varphi\optindex{#1}}}
\newcommand\J{\mathcal J}
\newcommand\Lagr{\mathcal L}
\newcommand\dphi{\partial_t \varphi}
\newcommand\subphi[2]{_{\set{\varphi(#1) > \varphi(#2)}}}
\newcommand\subphih[2]{_{\set{\varphi_h(#1) > \varphi_h(#2)}}}
\newcommand\subphim[2]{_{\set{\varphi^-_{#1} > \varphi^-_{#2}}}}
\newcommand\subdphi[1]{_{\set{\dphi(#1) > 0}}}
\newcommand\subdphiI[1]{_{\set{\dphi > 0, #1}}}
\newcommand\CIQ{W}
\newcommand\dtm[1][m]{\Delta t_#1}
\newcommand\notch{\mathcal N}
\newcommand\diff{\,\mathrm{d}}
\newcommand\dt{\diff t}
\newcommand\dx{\diff x}
\newcommand\free{_{\text{free}}}

\theoremstyle{plain} 
\newtheorem{theorem}{Theorem}[section]

\newtheorem{Remark}[theorem]{Remark}

\definecolor{ao(english)}{rgb}{0.0, 0.5, 0.0}

\usepackage[top=2cm, bottom=2cm, left=2cm, right=2cm]{geometry}

\begin{document}

\title
{Space-time formulation, discretization, and computational performance
studies for phase-field fracture optimal control problems}

\author[1]{D. Khimin}
\author[1]{M. C. Steinbach}
\author[1,2]{T. Wick}

\affil[1]{Leibniz Universit\"at Hannover, Institut f\"ur Angewandte
  Mathematik, Welfengarten 1, 30167 Hannover, Germany}

\affil[2]{Universit\'e Paris-Saclay, ENS Paris-Saclay, LMT -- Laboratoire de M\'ecanique et Technologie, 91190 Gif-sur-Yvette, France}

\date{\today}
\maketitle

\begin{abstract}
  The purpose of this work is the development of space-time discretization
  schemes for phase-field optimal control problems.
  Specifically in the optimal control
  minimization problem, a tracking-type cost functional is minimized
  to steer the crack via the phase-field variable into a desired pattern.
  To achieve such optimal solutions, Neumann type boundary conditions 
  need to be determined.
  First, a time discretization of the forward problem is derived
  using a discontinuous Galerkin formulation.
  Here, a challenge is to include regularization terms
  and the crack irreversibility constraint.
  The optimal control setting is formulated by means of the Lagrangian approach
  from which the primal part, adjoint, tangent and adjoint Hessian are derived.
  Herein the overall Newton algorithm is based
  on a reduced approach by eliminating the state constraint, 
  namely the displacement and phase-field unknowns, but 
  keeping the control variable as the only unknown.
  From the low-order discontinuous Galerkin discretization,
  adjoint time-stepping schemes are finally obtained.
  Both our formulation and algorithmic developments 
  are substantiated and illustrated
  with six numerical experiments.
  \textbf{AMS}:\\
  74R10,
  65N30,
  49M15,
  49K20,
  35Q74
\end{abstract}

\section{Introduction}
Fracture propagation using variational approaches and phase-field methods
is currently an important topic in applied mathematics and engineering.
The approach for the mathematical and mechanical literature
was established in \cite{FraMar98,BourFraMar00,KuMue10} and in physics
by \cite{KaKeLe01,HaKa09,SpBrKa11}
and overview articles and monographs include
\cite{BourFraMar08,BouFra19,WuNgNgSuBoSi19,Wi20_book,Fra21,DiLiWiTy22}
with numerous further references cited therein.

It is well-known that the efficient and robust numerical solution of
the nonlinear and linear subproblems in phase-field fracture is challenging.
This is mainly due to the nonlinear structure of the coupled problem
and the interaction of model, discretization and material
parameters. In spite of the development of
robust preconditioning and parallel, scalable, iterative algorithms,
the forward solution remains costly in general for
both two-dimensional and three-dimensional settings.

While the major amount of work concentrates on forward modeling
of phase-field fracture, more recently some work started
on parameter identification employing Bayesian inversion
\cite{KhoNoPaAbWiHei20,WuRoLoMa21,NoKhoUlAlWiFrWr21,NoKhoWi21},
topology optimization \cite{DESAI2022111048},
stochastic phase-field modeling \cite{GERASIMOV2020113353},
and optimal control \cite{NeiWiWo17,NeiWiWo19,MoWo20}.
Solving phase-field fracture problems using methods from shape optimization
was proposed in \cite{Allaire20115010}.

The main objective of this work is to design a mathematical
framework including computational performance studies
for phase-field fracture optimal control problems. In optimal
control some cost functional shall be minimized where the forward problem
(here the phase-field fracture weak formulation) acts as constraint and
the control (often involving boundary conditions or right hand side forces) is designed
in such a way that the minimization goal is achieved as well as possible.
Specifically, we consider a tracking type cost functional in which a desired
phase-field crack pattern shall be realized by controlling Neumann type boundary forces. From
an engineering viewpoint, such cost functionals and controls are reasonable since often
either a desired crack path shall be achieved (for instance in hydraulic fracturing), or
in the case of preventing fracture/damage either no fractures should develop or at least once
they start developing they should be steered with appropriate forces into directions that cause minimal
damage. Since volume forces such as gravitational forces play minor roles in such settings,
we mainly control optimal fracture patterns by boundary forces and concentrate on Neumann type
conditions.
In prior work \cite{NeiWiWo17,NeiWiWo19} 
the emphasis was on mathematical analysis and
a brief illustration in terms of a numerical simulation 
for a fixed fracture by
a tracking type functional for a desired displacement field.
However, computational details have not yet been discussed therein,
but are necessary in order to apply and investigate the methodology
for more practical applications such as propagating fractures.

Due to the irreversibility constraint on the fracture growth,
optimization problems subject to such an evolution become
mathematical programs with complementarity
constraints (MPCC) \cite{Barbu1984,MIGNOT1976,MIGNOT1984}
so that standard constraint qualifications like 
\cite{Robinson1976,Zowe1979} cannot hold.
Our computational approach requires stronger regularity
and hence we replace the complementarity constraint
with a suitable penalty term.

Designing a computational framework for phase-field fracture optimal
control is novel and challenging
because appropriate
function spaces and weak formulations need to be determined, and
robust forward and optimization solvers are required.

Specifically, we are interested in a rigorous mathematical framework,
which is the reason why we concentrate on one type of cost functional (here tracking type)
and one type of controls (here Neumann boundary controls) in this work. However,
from an engineering perspective other controls such as Dirichlet controls or right hand
side controls would be possible, too. Moreover, other cost functionals controling for instance
the bulk or crack energies could be employed. Technically, such implementations
can be realized in our software as shown for other numerical experiments in \cite{dope,DOpElib}.

For the forward solver, as intensively discussed in the literature, the linear
and nonlinear solutions are demanding because of the non-convexity of the
governing energy functional of the forward phase-field fracture model
and the relationship of discretization and regularization parameters.
For the nonlinear solution various methods were proposed
such as alternating minimization (staggered solution) \cite{Bour07,BuOrSue10},
quasi-monolithic solutions \cite{HeWheWi15,JoLaWi20,Wi20_book},
and fully monolithic schemes
\cite{GeLo16,Wi17_SISC,Wi17_CMAME,KoKr20,wambacq2020interiorpoint}.
Nonetheless, monolithic solutions as adopted here remain difficult
and we add an additional viscous regularization term
as originally proposed in \cite{KneRoZa13}
and used in our governing model from \cite{NeiWiWo19}.

The optimal control problem is formulated in terms of the reduced approach
by eliminating the state variable with a control-to-state operator;
see for instance \cite{Troe09,HiPiUlUl09}.
In this work, the state variable consists 
of the vector-valued displacement field and the 
phase-field variable. The control 
variable is a function defined on the boundary of the domain.
By eliminating the state variable, we obtain a so-called reduced 
cost functional 
defined in terms of the control variable only, which results 
in an unconstrained optimization problem.  The numerical solution is obtained
via the first-order necessary optimality condition.

Applying Newton-type methods requires the second derivative
of the reduced cost functional, and needs in practice 
the evaluation of the
\textit{adjoint, tangent}, and \textit{adjoint Hessian equations}.
The latter requires the evaluation of second-order derivatives;
see, e.g., \cite{BeMeVe07} and \cite[Chapter 4]{Me08} 
for parabolic optimization problems.

The paper \cite{BeMeVe07} serves as point of departure for our approach
in the current work.
Specifically, we employ Galerkin formulations in time
and discuss in detail how the crack irreversibility constraint is formulated
using a penalization \cite{MiWheWi13a,NeiWiWo17}
and an additional viscous regularization \cite{NeiWiWo19,KneRoZa13}.
Based on these settings, concrete time-stepping schemes are derived.
As usual, the primal and tangent problem run forward in time
whereas the adjoint and adjoint Hessian equations run backward in time.
The main emphasis is to establish robust
numerical solvers in terms of the nonlinear forward solver
and the nonlinear optimization loop. 

We then perform extensive tests by
means of six numerical experiments with different complexities.
First, we notice that propagating fractures
for such optimization problems were not addressed
in the prior work \cite{NeiWiWo17,NeiWiWo19}.
In the current work, considering now propagating fractures,
the overall goals are computational investigations of the performance
of the reduced Newton algorithm (NLP), the linear conjugate gradient (CG) method,
and convergence of the residuals, cost functionals,
tracking parts, Tikhonov parts and optimal controls.
We recall upfront that such investigations
are even challenging for forward phase-field fracture problems due to the
interaction of model, discretization and
material parameters (see \cite{BourFraMar00,KOLDITZ2022100047}
and closely related work on image segmentation \cite{Bou99},
and the prior seminal work on Gamma convergence \cite{AmTo90,AmTo92,Braides1998}),
and possibly also penalization parameters
for treating the crack irreversibility constraint \cite{Wi20_book}.
All of them have an impact on mathematical well-posedness
\cite{BourFraMar08} (and references cited therein),
and on numerical approximations
and nonlinear and linear solution algorithms
 \cite{BourFraMar08,BouFra19,WuNgNgSuBoSi19,Wi20_book,Fra21,DiLiWiTy22}.
The parameters include:
phase-field regularization $\varepsilon$ and bulk regularization $\kappa$,
crack irreversibility penalization $\gamma$ and viscous regularization $\eta$,
mesh size $h$ and loading step size $\Delta t$,
critical energy release rate $G_c$ and Lam\'e parameters $\lambda$ and $\mu$.

The extension to optimization adds further levels of complexity:
the forward problem, with all its own challenges,
must be solved numerous times,
more parameters enter such as the Tikhonov
regularization $\alpha$, and
in order to guarantee a well-posed optimization setting,
the adjustment of $\alpha$ is delicate for weighting the
physical tracking functional against the Tikhonov regularization term.
Our experiments below encompass
propagating fractures, non-constant controls
on one or more boundary sections,
multiple (propagating) fractures,
an adaptation of Winkler's \cite{winkler2001}
L-shaped panel test,
and using controls to prevent crack growth.
These tests provide novel insight for both
the capabilities of the phase-field method for fracture
from a numerical viewpoint as well as for applications.
On the other hand, limitations and opportunities for
future work also become visible, such as the need to further improve
the linear solver's cost complexity
(e.g., by parallel multigrid methods \cite{JoLaWi20_parallel}
and model order reduction \cite{BeCoOhlWill15})
as for fine meshes the forward solver becomes prohibitivly expensive.
Some further preliminary results
(yet with a stationary, non-propagating fracture) are published in
the book chapter \cite{KhiSteiWi21}.

The outline of this paper is as follows:
In \cref{sec_PFF}, the phase-field fracture forward model is introduced.
Furthermore, a Galerkin time discretization is derived.
Next, in \cref{sec_opt}, the optimization problem is stated,
including the reduced space approach.
In the key \cref{sec_Lag_aux}, the Lagrangian and three auxiliary equations
are carefully derived in great detail, and with 
our final complete algorithm. Then,
in \cref{SEC:Numtests}, extensive studies with six numerical experiments are discussed
to substantiate our algorithmic developments.
Our work is summarized in \cref{sec_conclusions}.

\section{Space-time phase-field fracture forward model}
\label{sec_PFF}
To formulate the space-time forward problem, 
we first introduce some basic notation
and then proceed with the construction of function spaces and 
a space-time weak formulation. Afterwards, a space-time Galerkin 
discretization is derived with discontinuous (dG) functions in time and 
a classical continuous Galerkin (cG) method in space.

\subsection{Notation}

We consider a bounded domain $\Omega \subset \R^2$.
The boundary is partitioned as
$\partial \Omega = \Gamma_N \stackrel.\cup \Gamma_D$
where both the Dirichlet boundary $\Gamma_D$ and
the Neumann boundary $\Gamma_N$ have nonzero Hausdorff measure.
Next we define two function spaces,
$V \coloneqq H_D^1(\Omega;\R^2) \times H^1(\Omega)$
for the displacement field $u$ and the phase-field $\varphi$,
and $Q \coloneqq L^2(\Gamma_N)$ for the control $q$, where
\begin{align*}
   H^1(\Omega;\R^2)
  &\coloneqq \defset{v \in L^2(\Omega;\R^2)}
    {D^{\alpha}v \in L^2(\Omega;\R^2) \tfor
    \alpha \in \mathbb{N}_0^2, \ \abs\alpha \le 1}, \\
  H^1_D(\Omega; \R^2)
  &\coloneqq \defset{v \in H^1(\Omega; \R^2)}{v|_{\Gamma_D} = 0}.
\end{align*}
Moreover we consider a bounded time interval $I = [0, T]$
and introduce the spaces
\begin{equation*}
  X \coloneqq
  \defset{\boldu = (u, \varphi)}
  {\boldu \in L^2(I, V), \,
    \dphi \in L^2(I, H^{1}(\Omega)^*)}, \qquad
  \CIQ \coloneqq C(I, Q).
\end{equation*}
Here, $H^{1}(\Omega)^*$ denotes the dual space
to $H^1(\Omega)$, which can be identified via the well-known Hilbert
space isomorphism with $H^1(\Omega)$ such that
$H^{1}(\Omega)^* \simeq H^1(\Omega)$.

On $V$ respectively $X$ we define the scalar products
\begin{align*}
  (\boldu, \boldv)
  &\coloneqq \int_\Omega \boldu \cdot \boldv \dx
  \qfor \boldu, \boldv \in V, \\
  (\boldu, \boldv)_I
  &\coloneqq \int_I \int_\Omega \boldu \cdot \boldv \dx \dt
    = \int_I (\boldu(t), \boldv(t)) \dt
    \qfor \boldu, \boldv \in X,
\end{align*}
with induced norms $\norm{}$ and $\norm{}_I$,
and furthermore the restricted inner products
\begin{align*}
  (\boldu(t), \boldv(t))\subdphi{t}
  &\coloneqq \begin{cases}
    (\boldu(t), \boldv(t)), & \dphi(t) > 0, \\
    0, & \text{else},
  \end{cases} \\
  (\boldu, \boldv)\subdphiI{I}
  &\coloneqq \int_I (\boldu(t), \boldv(t))\subdphi{t} \dt
    \qfor \boldu, \boldv \in X.
\end{align*}

Later we also work with
$(\fcdot, \fcdot)\subphi{t_i}{t_j}$,
defined like $(\fcdot, \fcdot)\subdphi{t}$,
and with a semi-linear
form $a(\fcdot)(\fcdot)$ in which the first argument is nonlinear
and the second argument is linear.

\subsection{Weak formulation}
We deal with the following weak formulation:
Given $\boldu_0 \in V$ and $q \in \CIQ$,
find $\boldu \in X$ such that
\begin{equation}
  \label{model}
  \begin{aligned}
    (g(\varphi) \C e(u), e(\Phiu))_I - (q, \Phi_{u:\perp})_{\Gamma_N,I} &= 0, \\
    G_c \varepsilon (\nabla \varphi, \nabla \Phiphi)_I -
    \smash[b]{\frac{G_c}{\varepsilon}} (1 - \varphi, \Phiphi)_I +
    (1 - \kappa) (\varphi \C e(u) : e(u), \Phiphi)_I \\
    {} + \gamma (\dphi, \Phiphi)\subdphiI{I} + \eta (\dphi, \Phiphi)_I &= 0,
  \end{aligned}
\end{equation}
for every test function $\boldPhi = (\Phiu,\Phiphi) \in X$.
Herein $ \Phi_{u:\perp}$ denotes the component of $\Phiu$
that is orthogonal to $\Gamma_N$.
The critcial energy release rate is denoted by $G_c>0$.
The so-called degradation function
$g(\varphi)\coloneqq (1-\kappa)\varphi^2 + \kappa$
helps to extend the displacements to the entire domain $\Omega$.
The bulk regularization parameter is $\kappa >0$,
the phase-field regularization parameter is $\varepsilon>0$,
the penalization parameter is $\gamma>0$,
and the viscosity parameter is $0<\eta\ll\gamma$.
Furthermore, $\C$ denotes the elasticity tensor
and $e(u)$ the symmetric gradient.
Then, we have
\[
  \C e(u) = \sigma(u) = 2 \mu e(u) + \lambda \operatorname{tr}(e(u)) I,
\]
where $\mu,\lambda>0$ are the Lam\'e parameters and $I$ is the identity matrix.

\begin{Remark}
The above weak formulation differs slightly from many other phase-field fracture
formulations found in the literature since the crack irreversibility constraint
$\partial_t \varphi \leq 0$ is kept on the time-continuous level in order to apply
a Galerkin discretization in time.
\end{Remark}

\begin{Remark}[Initial condition $u_0$]
\label{rem_ini_u0}
Note that we are concerned with quasi-static brittle fracture
without explicit time derivative in the displacement equation.
Nonetheless, we introduce for formal reasons $u_0$.
First, we can develop in an analogous fashion time discretization schemes
for the overall forward model.
Second, it facilitates the extension to problems
in which the displacement equation does have a time derivative,
such as dynamic fracture \cite{BouLarRi11,BoVeScoHuLa12}.
Third, having $u_0$ allows for a monolithic implementation structure,
and the system matrix for the initial condition is regular.
\end{Remark}
\begin{Remark}[Convexification]
  We notice that strict positivity $\eta > 0$ 
  improves the numerical solution process of \eqref{model}.
  In fact, one can show for the quasi-static case
  that for sufficiently large values of $\eta$ the control-to-state mapping
  is single valued
  due to strict convexity of the energy corresponding to the equation.
  However, the convexification term $\eta (\dphi, \Phiphi)_I$
  also penalizes crack growth.
  To ensure the dominance of the physically motivated term
  $\gamma (\dphi, \Phiphi)\subdphiI{I}$
  we have to choose $\gamma \gg \eta$.
\end{Remark}

\subsection{Space-time finite element discretization}

\subsubsection{Temporal discretization}

Given $T > 0$, we define the time grid $0 = t_0 < \dotsb < t_M = T$
to partition the interval $I$ into $M$ left-open subintervals
$I_m = (t_{m-1}, t_m]$,
\[
  I = \set{0} \cup I_1 \cup \dotsb \cup I_M.
\]
By using the discontinuous Galerkin method, here dG(0),
we seek a solution $\boldu$ in the space $X^0_k$
of piecewise polynomials of degree $0$,
\begin{equation}
  \label{eq_X_k^0}
  X^0_k \coloneqq \defset{\boldv \in X}
  {\boldv(0) \in V  \text{ and }
    \boldv|_{I_m} \in \polynomials_0(I_m, V), \, m = 1, \dots, M}.
\end{equation}
Here, the subindex $k$ indicates the time-discretized function space
in order to distinguish it from the continuous space $X$.
For the jump terms arising in $X^0_k$ we use the standard notation
\[
  \boldv^+_m \coloneqq \boldv(t_m+), \qquad
  \boldv^-_m \coloneqq \boldv(t_m-) = \boldv(t_m), \qquad
  [\boldv]_m \coloneqq \boldv^+_m - \boldv^-_m.
\]

\begin{Remark}
\label{rem:constant_time}
Since we work with dG(0), i.e., piece-wise constant functions in time, we have
\[
  \partial_t \boldv = \boldv_m^- - \boldv_{m-1}^+ = 0
  \qfor \boldv \in X_k^0 \tfor m = 1, \dots, M.
\]
\end{Remark}
The discretized state equation combines the two equations of~\eqref{model}.
For a concise formulation, the energy-related terms are expressed
as a semi-linear form $a\colon Q \times V \times V \to \R$,
\begin{equation*}
  \label{State_sum_4}
  \begin{aligned}
    a(q,\boldu)(\boldPhi)
    &\coloneqq g(\varphi) \cdot (\C e(u), e(\Phiu)) \\
    &\, + G_c \varepsilon (\nabla \varphi, \nabla \Phiphi)
    - \frac{G_c}{\varepsilon} (1 - \varphi, \Phiphi) \\
    &\, + (1 - \kappa) (\varphi \cdot \C e(u) : e(u), \Phiphi)
    - (q, \Phi_{u:\perp})_{\Gamma_N}
    .
  \end{aligned}
\end{equation*}
Now the fully discretized state equation determines
a function $\boldu \in X^0_k$ for a given
initial value $\boldu_0 = (u_0, \varphi_0) \in V$ and a given
control $q \in W$
such that for every $\boldPhi \in X^0_k$
\begin{subequations}
  \label{State_dG}
  \begin{align}
    \label{State_dG_a}
    0
    &= \sum_{m=1}^M \bigl[
      \gamma (\dphi, \Phiphi)\subdphiI{I_m}
      + \eta (\dphi, \Phiphi)_{I_m} \bigr] \\
    \label{State_dG_b}
    &+ \sum_{m=0}^{M-1} \bigl[
      \gamma ([\varphi]_m, \Phiphi[m]^+)\subphim{m+1}{m}
      + \eta ([\varphi]_m, \Phiphi[m]^+) \bigr] \\
    \label{State_dG_c}
    &+ \sum_{m=1}^M a(q(t_m),\boldu(t_m))(\boldPhi(t_m)) \dtm \\
    \label{State_dG_d}
    &+ (u^-_0 - u_0, \Phiu[0]^-) + (\varphi^-_0 - \varphi_0, \Phiphi[0]^-)
    .
  \end{align}
\end{subequations}
The time integral in \eqref{State_dG_c} has been approximated
by the right-sided box rule, where $\dtm \coloneqq t_m - t_{m-1}$.
Discontinuities of the functions in $X^0_k$ are captured
by the jump terms in \eqref{State_dG_b} in the typical dG(0) manner.
These jump terms can be rewritten as
\begin{equation}
  \label{State_dG_tilde}
  \sum_{m=1}^M \bigl[
  \gamma (\varphi^+_{m-1} - \varphi^-_{m-1}, \Phiphi[m-1]^+)\subphim{m}{m-1}
  + \eta (\varphi^+_{m-1} - \varphi^-_{m-1}, \Phiphi[m-1]^+) \bigr].
\end{equation}
Moreover, since we are employing a dG(0) scheme,
our test functions satisfy
\[
  \boldPhi^+_{m-1} = \boldPhi^-_m \qfor m = 1, \dots, M.
\]
Thus the first sum \eqref{State_dG_a}
vanishes entirely by \cref{rem:constant_time},
and the two terms containing $\varphi_{m-1}^+$
in \eqref{State_dG_tilde} become
$(\varphi^-_m, \Phiphi[m]^-)\subphim{m}{m-1}$ and
$(\varphi^-_m, \Phiphi[m]^-)$, respectively.
Together with \eqref{State_dG_b} and \eqref{State_dG_d},
the discrete state equation \eqref{State_dG} is finally written as
\begin{equation}
  \label{State_sum_5}
  \begin{aligned}
    0
    &= \sum_{m=1}^M \Bigl( \gamma \bigl[
    (\varphi^-_m, \Phiphi[m]^-)\subphim{m}{m-1} -
    (\varphi^-_{m-1}, \Phiphi[m]^-)\subphim{m}{m-1} \bigr] \\[-2\jot]
    &\mkern70mu + \eta \bigl[
    (\varphi^-_m, \Phiphi[m]^-) -
    (\varphi^-_{m-1}, \Phiphi[m]^- ) \bigr] \\
    &\mkern70mu + a(q(t_m),\boldu(t_m))(\boldPhi(t_m)) \dtm \Bigr) \\
    &+ (u^-_0 - u_0, \Phiu[0]^-) + (\varphi^-_0 - \varphi_0, \Phiphi[0]^-)
    \qfor \boldPhi \in X^0_k
    .
  \end{aligned}
\end{equation}
To solve \eqref{State_sum_5}, we first obtain $\boldu_0^- = \boldu(0)$
from the initial condition
\begin{equation}
  \label{state_t=0}
  (\boldu(0), \boldPhi_0^-) = (\boldu_0, \boldPhi_0^-)
  \qfor \boldPhi_0^- \in V.
\end{equation}
Then we compute $\boldu(t_m)$ for $m = 1, \dots, M$ from
\begin{equation}
  \begin{aligned}\label{state_t>0}
    0
    &= \gamma (\varphi(t_m), \Phiphi(t_m))\subphi{t_m}{t_{m-1}}
    + \eta (\varphi(t_m), \Phiphi(t_m)) \\
    &- \gamma (\varphi(t_{m-1}), \Phiphi(t_{m}))\subphi{t_m}{t_{m-1}}
    - \eta (\varphi(t_{m-1}), \Phiphi(t_{m})) \\
    &+ a(q(t_m),\boldu(t_m))(\boldPhi(t_m)) \dtm
    \qfor \boldPhi \in X^0_k
    .
  \end{aligned}
\end{equation}

\subsubsection{Spatial discretization}

For the spatial discretization, we employ again a Galerkin
finite element scheme by introducing $H^1$ conforming discrete spaces.
We consider two-dimensional shape-regular meshes
with quadrilateral elements $K$ forming  the mesh
$\mathcal{T}_h = \set{K}$; see \cite{Cia87}.
The spatial discretization parameter is the
diameter $h_K$ of the element $K$.
On the mesh $\mathcal{T}_h$ we construct a finite element space
$V_h \coloneqq V_{uh} \times V_{\varphi h}$ as usual:
\begin{align*}
  V_{uh} &\coloneqq
  \defset{v \in H^1_D(\Omega; \R^2)}
  {v|_K \in Q_s(K) \text{ for } K \in \mathcal{T}_h}, \\
  V_{\varphi h} &\coloneqq
  \defset{v \in H^1(\Omega)}
  {v|_K \in Q_s(K) \text{ for } K \in \mathcal{T}_h}.
\end{align*}
Herein $Q_s(K)$ consists of shape functions
that are obtained as bilinear transformations of functions
defined on the master element $\hat{K} = (0, 1)^2$,
where $\hat Q_s(\hat K)$ is the space
of tensor product polynomials up to degree $s$ in dimension $d$, defined as
\[
  \hat Q_s (\hat K) \coloneqq \operatorname{span}
  \Defset{\prod_{i=1}^d \hat x_i^{\alpha_i}}{\alpha_i \in \set{0, 1 \dots, s}}.
\]
Specifically, for $s = 1$ and $d = 2$ we have
\[
  \hat Q_1(\hat K)
  =
  \operatorname{span}\set{1, \hat x_1, \hat x_2, \hat x_1 \hat x_2}.
\]
With these preparations, based on \eqref{eq_X_k^0},
we now design the fully discrete function space
\begin{equation*}
  X^0_{hk} \coloneqq \defset{\boldv \in X}
  {\boldvh(0) \in V_h \text{ and }
    \boldv|_{I_m} \in \polynomials_0(I_m, V_h), \, m = 1, \dots, M}.
\end{equation*}
The discrete control space $Q_h$ is constructed like $X^0_{hk}$
using $Q_1(K)$ (again $s=1$)
elements, but restricted to the Neumann boundary $\Gamma_N$.
Then, the fully discrete system consists of the initial condition
\begin{equation}
  \label{state_th=0}
  (\boldu_h(0), \boldPhi_{h,0}^-) = (\boldu_{h,0}, \boldPhi_{h,0}^-)
  \qfor \boldPhi_{h,0}^- \in V_h
\end{equation}
and for $m = 1, \dots, M$ of the local system
\begin{equation}
  \begin{aligned}\label{state_th>0}
    0
    &= \gamma (\varphi_h(t_m), \Phiphih(t_m))\subphih{t_m}{t_{m-1}}
    + \eta (\varphi_h(t_m), \Phiphih(t_m)) \\
    &- \gamma (\varphi_h(t_{m-1}), \Phiphih(t_{m}))\subphih{t_m}{t_{m-1}}
    - \eta (\varphi_h(t_{m-1}), \Phiphih(t_{m})) \\
    &+ a(q_h(t_m),\boldu_h(t_m))(\boldPhi_h(t_m)) \dtm
    \qfor \boldPhi_h \in X^0_{hk}.
  \end{aligned}
\end{equation}

\section{Optimization with phase-field fracture}
\label{sec_opt}

In this section, we state the phase-field optimal control problem and
introduce the reduced solution approach. Therein, the primal forward problem
plus three additional equations must be solved. Their combination
yields the final solution algorithm.

\subsection{Optimization problem}

We consider a separable NLP (Non-Linear Program)
with a cost functional of tracking type.
In this tracking type functional, the objective is to approximate
a given phase-field fracture pattern $\varphi_d$ by determining a suitable
control $q$.
The corresponding minimization problem is given by:
\begin{equation}
\label{EQ:NLPgamma}
\begin{aligned}
  \min_{q,\boldu} \quad & \J(q,\boldu) \coloneqq
  \frac{1}{2}\sum_{m=1}^M  \norm{\varphi(t_m) - \varphi_d(t_m)}^2
  + \frac{\alpha}{2}\sum_{m=1}^M  \norm{q(t_m) - q_d(t_m)}_{\Gamma_N}^2 \\
  \textq{s.t.} &\text{ \eqref{state_t=0} and \eqref{state_t>0} for
 $m=1,\dots,M$, with $(u_0,\varphi_0) \in V$ and $(q,\boldu) \in \CIQ \times X$,}
\end{aligned}
\end{equation}

where $\varphi_d \in L^{\infty}(\Omega)$ is some desired phase-field
and $q_d$ is a suitable nominal control
that we use for numerical stabilization.
The second sum represents a common Tikhonov regularization
with parameter $\alpha$.
The existence of a global solution of \eqref{EQ:NLPgamma}
in $L^2(I,Q) \times X$
has been shown in \cite[Theorem 4.3]{NeiWiWo17}
for functions that are non-negative and weakly semi-continuous.

\begin{Remark}
The fully discrete version of \eqref{EQ:NLPgamma} is obtained by working
with the equations \eqref{state_th=0} and \eqref{state_th>0}. In what follows,
in order to keep the notation comfortable, we omit the index $h$ indicating
the spatial discretization.
\end{Remark}

\subsection{Reduced optimization problem}

In order to handle \eqref{EQ:NLPgamma} by the reduced approach,
we assume that a solution operator
$S\colon \CIQ \to X$ exists for the PDE \eqref{model}.
The cost functional $\J(q, \boldu)$ then reduces to
$j\colon \CIQ \to \R$, $j(q) \coloneqq \J(q, S(q))$,
and we replace \eqref{EQ:NLPgamma} by the unconstrained optimization problem
\begin{equation}
  \label{EQ:NLPgammared}
  \min_q \ j(q).
\end{equation}
To solve $j'(q) = 0$ by Newton's method,
we compute representations of $j'$ and $j''$
using the established approach in \cite{BeMeVe07}.
It requires the solution of four equations (given below)
for derivatives of the Lagrangian
$\Lagr\colon \CIQ \times X^0_k \times X^0_k \to \R$,
which is defined within the dG($r$) setting as
\begin{equation}
  \label{Lagrange_dG}
  \begin{aligned}
    \Lagr(q,\boldu,\boldz)
    \coloneqq \J(q,\boldu)
    &-\sum_{m=1}^M \bigg( \gamma (\dphi, \zphi)\subdphiI{I_m} + \eta (\dphi, \zphi)_{I_m} \bigg)\\
    &-\sum_{m=0}^{M-1} \bigg( \gamma ([\varphi]_m, z^+_{\varphi,m})\subphim{m+1}{m}
      + \eta ([\varphi]_m, z^+_{\varphi,m}) \bigg) \\
    &- \int_I a(q(t),\boldu(t))(\boldz(t)) \dt \\
    &- \eta_0 (u(0) - u_0, \zu(0)) -\eta (\varphi(0) - \varphi_0, \zphi(0))
    ,
  \end{aligned}
\end{equation}
and for the time continuous case as

\begin{equation}
  \label{Lagrange}
  \begin{aligned}
    \Lagr(q,\boldu,\boldz)
    \coloneqq \J(q,\boldu)
    &- \gamma (\dphi, \zphi)\subdphiI{I} - \eta (\dphi, \zphi)_I \\
    &- \int_I a(q(t),\boldu(t))(\boldz(t)) \dt \\
    &- \eta_0 (u(0) - u_0, \zu(0)) -\eta (\varphi(0) - \varphi_0, \zphi(0))
    .
  \end{aligned}
\end{equation}

\begin{Remark}
We notice that starting with \eqref{Lagrange_dG} and deriving 
the state, adjoint, tangent, adjoint Hessian equations, exhibits 
the property that discretization and optimization interchange, i.e.,
the discretize-then-optimize and the optimize-then-discretize 
approaches are equal; see \cite{BeMeVe07,Me08} for parabolic 
optimization problems. However in what follows, 
we start from the time-continuous 
formulation \eqref{Lagrange} for the ease of presentation (which nonetheless 
becomes difficult enough) and we add only afterwards the dG-in-time
representations, which yields in the end the same result if we had started 
with \eqref{Lagrange_dG}.
\end{Remark}

\subsection{State, adjoint, tangent, adjoint Hessian}

In this section we state the four equations
to be solved for computing $j'$ and $j''$.
\begin{enumerate}
\item \emph{State equation:}
  given $q \in \CIQ$,
  find $\boldu {= S(q)} \in X$ such that the PDE \eqref{model} holds:
  \begin{equation*}
    \label{Stateaux}
    \Lagr'_{\boldz}(q,\boldu,\boldz)(\boldPhi) = 0
    \qfor \boldPhi \in X.
  \end{equation*}
\item \emph{Adjoint equation:} given $q \in \CIQ$ and $\boldu = S(q)$,
  find $\boldz \in X$ such that
  \begin{equation}
    \label{Adjointaux}
    \Lagr'_{\boldu}(q,\boldu,\boldz)(\boldPhi) = 0
    \qfor \boldPhi \in X.
  \end{equation}
\item \emph{Tangent equation:} given $q \in \CIQ$, $\boldu = S(q)$,
  and a direction $\delta q \in \CIQ$,
  find $\bolddeltau \in X$ such that
  \begin{equation}
    \label{Tangentaux}
    \Lagr''_{q\boldz}(q,\boldu,\boldz)(\delta q, \boldPhi) +
    \Lagr''_{\boldu\boldz}(q,\boldu,\boldz)(\bolddeltau,\boldPhi) = 0
    \qfor \boldPhi \in X.
  \end{equation}
\item \emph{Adjoint Hessian equation:} given $q \in \CIQ$, $\boldu = S(q)$,
  $\boldz \in X$ from \eqref{Adjointaux},
  a direction $\delta q \in \CIQ$,
  and $\bolddeltau \in X$ from \eqref{Tangentaux},
  find $\bolddeltaz \in X$ such that
  \begin{equation}
    \label{AdjointHessianaux}
    \Lagr''_{q\boldu}(q,\boldu,\boldz)(\delta q,\boldPhi) +
    \Lagr''_{\bolduu}(q,\boldu,\boldz)(\bolddeltau,\boldPhi) +
    \Lagr''_{\boldzu}(q,\boldu,\boldz)(\bolddeltaz,\boldPhi) = 0
    \qfor \boldPhi \in X.
  \end{equation}
\end{enumerate}
Solving these equations in a specific order
(see for instance \cite{BeMeVe07,Me08})
leads to the following representations
of the derivatives that we need for Newton's method:
\begin{equation}
\label{Representations}
\begin{aligned}
  j'(q)(\delta q)
  &= \Lagr'_q(q,\boldu,\boldz)(\delta q) \qfor \delta q \in \CIQ, \\
  j''(q)(\delta q_1, \delta q_2)
  &= \Lagr''_{qq}(q,\boldu,\boldz)(\delta q_1, \delta q_2) +
    \Lagr''_{\boldu q}(q,\boldu,\boldz)(\bolddeltau, \delta q_2) \\
  &+ \Lagr''_{\boldz q}(q,\boldu,\boldz)(\bolddeltaz,\delta q_2)
    \qfor \delta q_1, \delta q_2 \in \CIQ.
\end{aligned}
\end{equation}

\section{Auxiliary equations}
\label{sec_Lag_aux}
Starting from the Lagrangian \eqref{Lagrange}, we derive in detail the three auxiliary equations
\eqref{Adjointaux}--\eqref{AdjointHessianaux}.
Specific emphasis is on the regularization terms
for the crack irreversibility and the convexification.

\subsection{Adjoint}
In the adjoint for dG(0) we seek
$\boldz = (\zu, \zphi) \in X_k^0$ such that
\[
  \Lagr'_{\boldu}(q,\boldu,\boldz)(\boldPhi) = 0 \qfor \boldPhi \in X^0_k.
\]
The first interesting part is the calculation of the derivative of $\Lagr$.
We formulate it directly in the weak form
\begin{equation}
  \label{Lagrange_derivative_1}
  \begin{aligned}
    \Lagr'_{\boldu}(q,\boldu,\boldz) (\boldPhi)
    &= \J'_{\boldu}(q,\boldu)(\boldPhi) \\
    &- \gamma (\partial_t \Phiphi, \zphi)\subdphiI{I}
      - \eta (\partial_t \Phiphi, \zphi)_I \\
    &\,{}- \int_I a'_{\boldu}(q(t),\boldu(t))(\boldPhi(t),\boldz(t)) \dt \\
    &-\eta_0 (\Phiu(0), \zu(0)) -\eta (\Phiphi(0), \zphi(0))
    .
  \end{aligned}
\end{equation}
\begin{Remark}
We notice that $\gamma (\partial_t \Phiphi, \zphi)\subdphiI{I}$
is a suitable numerical approximation to the derivative 
of $\gamma (\partial_t \varphi, \zphi)\subdphiI{I}$, since formally 
a characteristic function must be differentiated; see also 
\cite[Section 5]{MiWheWi19} 
for a similar numerical approximation in the context 
of a related forward problem. The same procedure as numerical 
approximation of the derivative 
is utlized in the other three auxiliary problems, 
namely the adjoint, tangent, and adjoint Hessian.
\end{Remark}
The partial derivative of $a$ in \eqref{Lagrange_derivative_1} reads
\begin{equation}
  \label{a'_u}
  \begin{aligned}
    a'_{\boldu}(q,\boldu)(\boldPhi,\boldz)
    &= ((1-\kappa) \varphi^2 + \kappa) \cdot (\C e(\Phiu), e(\zu)) \\
    &+ 2 \varphi (1-\kappa) \Phiphi (\C e(u), e(\zu)) \\
    &+ G_c \varepsilon (\nabla \Phiphi, \nabla \zphi)
    + \frac{G_c}{\varepsilon} (\Phiphi, \zphi) \\
    &+ (1-\kappa) (\Phiphi \cdot \C e(u) : e(u), \zphi) \\
    &+ 2 \varphi (1-\kappa) (\C e(\Phiu) : e(u), \zphi)
    .
  \end{aligned}
\end{equation}
Now the main problem is that the time derivatives are applied
to the test function~$\boldPhi$ as usual in the adjoint.
Therefore we use integration by parts
to shift the time derivatives over to $\boldz$.
Then the second line
in \eqref{Lagrange_derivative_1} becomes
\begin{equation}
  \label{Lagrange_derivative_2}
  \begin{aligned}
    &\quad\gamma (\Phiphi, \partial_t \zphi)\subdphiI{I}
    + \eta (\Phiphi, \partial_t \zphi)_I \\
    &+ \gamma (\Phiphi(0), \zphi(0))\subdphi{0} + \eta (\Phiphi(0), \zphi(0)) \\
    &- \gamma (\Phiphi(T), \zphi(T))\subdphi{T} - \eta (\Phiphi(T), \zphi(T)) .
  \end{aligned}
\end{equation}
At this point we have to decide how to approximate
the time derivative $\dphi(0)$.
While $\dphi(t_m)$ for $m = 1, \dots, M$ is easily approximated
by the backward difference
\[
  \dphi(t_m) \approx \frac{\varphi(t_m) - \varphi(t_{m-1})}{t_m - t_{m-1}},
\]
this procedure will not work for the first mesh point $t_0 = 0$.
The forward difference
\[
  \dphi(0) \approx \frac{\varphi(t_1) - \varphi(t_0)}{t_1 - t_0}
\]
is a good choice because it simplifies the condition $\dphi(t_0) > 0$
to $\varphi(t_1) > \varphi(t_0)$ and leads to desired
cancelations in \eqref{Lagrange_derivative_dG}.
Now we will repeat the procedure that we applied to the state equation.
We approximate the time derivatives and add the jump terms
(with shifted index) as we did in \eqref{State_dG},
obtaining expressions similar to \eqref{State_dG_tilde}:
\begin{equation}
  \label{Lagrange_derivative_dG}
  \begin{aligned}
    \Lagr'_{\boldu}(q,\boldu,\boldz) (\boldPhi)
    &= \J'_{\boldu}(q,\boldu)(\boldPhi) \\
    &+ \sum_{m=1}^M \bigl[ \gamma
      (\Phiphi[m]^-, \zphi[m]^- - \zphi[m-1]^+)\subphim{m}{m-1} \\[-2\jot]
    &\mkern60mu + \eta (\Phiphi[m]^-, \zphi[m]^- - \zphi[m-1]^+) \bigr] \\
    &- \gamma (\Phiphi[M]^-, \zphi[M]^-)\subphi{t_M}{t_{M-1}} -
    \eta (\Phiphi[M]^-, \zphi[M]^-) \\
    &+ \gamma (\Phiphi[0]^-, \zphi[0]^-)\subphi{t_1}{t_0} +
    \eta (\Phiphi[0]^-, \zphi[0]^-) \\
    &+ \sum_{m=1}^M \bigl[ \gamma
      (\Phiphi[m-1]^-, \zphi[m-1]^+ - \zphi[m-1]^-)\subphim{m}{m-1} \\[-2\jot]
    &\mkern60mu + \eta (\Phiphi[m-1]^-, \zphi[m-1]^+ - \zphi[m-1]^-) \bigr] \\
    &- \sum_{m=1}^M
    a'_{\boldu}(q(t_m),\boldu(t_m))(\boldPhi(t_m),\boldz(t_m)) \dtm \\
    &-\eta_0 (\Phiu[0]^-, \zu[0]^-) -\eta (\Phiphi[0]^-, \zphi[0]^-)
    .
  \end{aligned}
\end{equation}
Since $\zphi\in X^0_k$, we have $\zphi[m]^- = \zphi[m-1]^+$
and see that the first sum vanishes entirely.
We also see that the terms $\pm \eta (\Phiphi[0]^-, \zphi[0]^-)$
in the fifth and the last line of \eqref{Lagrange_derivative_dG} cancel.
Moreover, we assume that $\varphi(t_1) \leq \varphi(t_0)$ in the initial step,
and hence the term $\gamma (\Phiphi[0]^-, \zphi[0]^-)\subphi{t_1}{t_0}$
in the fifth line vanishes as well.
\begin{Remark}[Projection of the initial solution] \label{Remark_Projection_of_initial_solution}
  The assumption $\varphi(t_1) \leq \varphi(t_0)$ is numerically
  justified since at $t_0$ some initial phase-field solution is prescribed.
  From $t_0$ to $t_1$ an $L^2$ projection of the initial conditions
  is employed that conserves the crack irreversibility constraint.
\end{Remark}
By the above arguments we eliminate the second, third and fifth line
of \eqref{Lagrange_derivative_dG} and the second term of the last line,
whereas the initial values for $\zu$ are still present:
\begin{equation}
  \label{eq:Lu-final}
  \begin{aligned}
    \Lagr'_{\boldu}(q,\boldu,\boldz)(\boldPhi)
    &= \J'_{\boldu}(q,\boldu)(\boldPhi) \\
    &- \gamma (\Phiphi[M]^-, \zphi[M]^-)\subphi{t_M}{t_{M-1}}
      - \eta (\Phiphi[M]^-, \zphi[M]^-) \\
    &+ \sum_{m=1}^M \bigl[ \gamma
    (\Phiphi[m-1]^-, \zphi[m-1]^+ - \zphi[m-1]^-)\subphim{m}{m-1} \\[-2\jot]
    &\mkern60mu + \eta (\Phiphi[m-1]^-, \zphi[m-1]^+ - \zphi[m-1]^-) \bigr] \\
    &- \sum_{m=1}^M
      a'_{\boldu}(q(t_m),\boldu(t_m))(\boldPhi(t_m),\boldz(t_m)) \dtm \\
    &- \eta_0 (\Phiu[0]^-, \zu[0]^-)
    .
  \end{aligned}
\end{equation}

\subsection{Adjoint time-stepping scheme}
From here on we exploit the separable structure of
$\J(q,\boldu) = \sum_m J(q(t_m),\boldu(t_m))$.
We start the solution process by pulling out from \eqref{eq:Lu-final}
every term associated with the last time point $t_M$:
\begin{equation}
\label{Adjoint_t_M}
  \begin{aligned}
    a_{\boldu}'&(q(t_M)\boldu(t_M))(\boldPhi(t_M),\boldz(t_M)) \dtm[M] \\
    &+ \gamma (\Phiphi[M]^-, \zphi[M]^-)\subphim{m}{m-1}
    + \eta (\Phiphi[M]^-, \zphi[M]^-) \\
    &= J_{\boldu}'(q(t_M),\boldu(t_M))(\boldPhi(t_M))
    \qfor \boldPhi \in X^0_k
    .
  \end{aligned}
\end{equation}
Now we collect what is left,
multiply by $-1$ and use the $X_k^0$ property ($\zphi[m-1]^+ = \zphi[m]^-$):
\begin{equation}
  \begin{aligned}
    0
    &= \sum_{m=1}^M \bigl[ \gamma
    (\Phiphi[m-1]^-, \zphi[m-1]^- - \zphi[m]^-)\subphim{m}{m-1}
    + \eta (\Phiphi[m-1]^-, \zphi[m-1]^- - \zphi[m]^-) \bigr] \\
    &+ \sum_{m=1}^{M-1}
    a'_{\boldu}(q(t_m),\boldu(t_m))(\boldPhi(t_m),\boldz(t_m)) \dtm \\
    &- \sum_{m=1}^{M-1}
    J'_{\boldu}(q(t_m),\boldu(t_m))(\boldPhi(t_m)) \\
    &+ \eta_0 (\Phiu[0]^-, \zu[0]^-)
    \qfor \boldPhi \in X_k^0.
 \end{aligned}
\end{equation}
To formulate the equations that are actually solved in every time step
we want to rewrite the entire equation as a
single sum.
Therefore we shift down the index of the first sum (the jump terms),
take out the terms for $m = 0$, and obtain
\begin{equation*}
  \begin{aligned}
    0
    &=\sum_{m=1}^{M-1} \Bigl( \bigl[
    \gamma (\Phiphi[m]^-, \zphi[m]^- - \zphi[m+1]^-)\subphim{m+1}{m}
    + \eta (\Phiphi[m]^-, \zphi[m]^- - \zphi[m+1]^-) \bigr] \\[-2\jot]
    &\mkern70mu +
    a'_{\boldu}(q(t_m),\boldu(t_m))(\boldPhi(t_m),\boldz(t_m))
    \dtm \notag \\[\jot]
    &\mkern70mu -
    J'_{\boldu}(q(t_m),\boldu(t_m))(\boldPhi(t_m)) \Bigr) \\
    &+ \gamma (\Phiphi[0]^-, \zphi[0]^- - \zphi[1]^-)\subphim{1}{0}
    + \eta (\Phiphi[0]^-, \zphi[0]^- - \zphi[1]^-) \\
    &+ \eta_0 (\Phiu[0]^-, \zu[0]^-)
    .
  \end{aligned}
\end{equation*}
Now we solve for $m = M-1, M-2, \dots, 1$ the equation
\begin{equation}
\label{Adjoint_ste}
\begin{aligned}
  a_{\boldu}'&(q(t_m),\boldu(t_m))(\boldPhi(t_m),\boldz(t_m)) \dtm \\
  &+ \gamma (\Phiphi[m]^-, \zphi[m]^- - \zphi[m+1]^-)\subphim{m+1}{m}
    + \eta (\Phiphi[m]^-, \zphi[m]^- - \zphi[m+1]^-) \\
  &= J_{\boldu}'(q(t_m),\boldu(t_m))(\boldPhi(t_m))
  \qfor \boldPhi \in X^0_k
  .
\end{aligned}
\end{equation}
Finally three terms are left for $m = 0$,
\begin{equation}
  \label{Adjoint_t_0}
  \gamma (\Phiphi[0]^-, \zphi[0]^- - \zphi[1]^-)\subphim{1}{0}
  + \eta (\Phiphi[0]^-, \zphi[0]^- - \zphi[1]^-)
  + \eta_0 (\Phiu[0]^-,\zu[0]^-) = 0
  .
\end{equation}
For $\eta_0 \ll \eta$ small enough the last term of \eqref{Adjoint_t_0}
can be dropped and the following equation can be solved instead:
\begin{equation}\label{Adj_initial}
  (\Phiphi[0]^-, \zphi[1]^-) = (\Phiphi[0]^-, \zphi[0]^-).
\end{equation}
\begin{Remark}[Algorithmic realization]\label{Remark_Algorithmic_realization}
  To avoid singular matrices that would lead to a loss of convergence
  in the linear solvers,
  we have to add an intial condition for $\zu[0]^-$:
  $(\Phiu[0]^-, \zu[1]^-) = (\Phiu[0]^-, \zu[0]^-)$.
  In total we replace \eqref{Adj_initial} by
  $(\boldPhi_0^-, \boldz^-_1) =  (\boldPhi_0^-, \boldz^-_0)$.
  We also refer the reader to the third reason outlined in \cref{rem_ini_u0}.
\end{Remark}

\subsection{Tangent equation}
The second auxiliary equation is the tangent equation.
In this equation we seek
$\bolddeltau = (\delta u, \delta\varphi) \in X^0_k$ such that
\[
  \Lagr''_{q\boldz}(q,\boldu,\boldz)(\delta q,\boldPhi) +
  \Lagr''_{\bolduz}(q,\boldu,\boldz)(\bolddeltau,\boldPhi) = 0
  \qfor \boldPhi \in X_k^0.
\]
Here we will apply the same procedure as for the state equation.
Recall that $\Lagr(q,\boldu,\boldz)$ contains the integrand
$a(q(t),\boldu(t))(\boldz(t))$ with $\boldz(t)$ entering linearly.
Hence the partial derivative required for
$\Lagr''_{\bolduz}(q,\boldu,\boldz)(\bolddeltau,\boldPhi)$
is simply $a'_{\boldu}(q,\boldu)(\bolddeltau,\boldPhi)$,
and the partial derivative required for
$\Lagr''_{q\boldz}(q,\boldu,\boldz)(\delta q,\boldPhi)$
can be derived from \eqref{State_sum_4} as
\begin{equation}
  \label{a'_q}
  a'_q(q,\boldu)(\delta q, \boldPhi) = -(\delta q, \Phiuy)_{\Gamma_N}.
\end{equation}
Furthermore, $\J(q,\boldu)$ does not depend on $\boldz$,
hence $\J''_{q\boldz}$ and $\J''_{\boldu\boldz}$ vanish.
Using the right-sided box rule again,
we thus obtain the discretized tangent equation
\begin{equation}
  \label{Tangent}
  \begin{aligned}
    0
    &= \sum_{m=1}^M \bigl[ \gamma
    (\delta\varphi^-_m - \delta\varphi^+_{m-1}, \Phiphi[m]^-)\subphim{m}{m-1}
    + \eta (\delta\varphi^-_m - \delta\varphi^+_{m-1}, \Phiphi[m]^-) \bigr] \\
    &+ \sum_{m=1}^M
    a'_{\boldu}(q(t_m),\boldu(t_m))(\bolddeltau(t_m),\boldPhi(t_m)) \dtm \\
    &+ \sum_{m=0}^{M-1} \bigl[ \gamma
    (\delta\varphi^+_m - \delta\varphi^-_m, \Phiphi[m]^+)\subphim{m+1}{m}
    + \eta (\delta\varphi^+_m - \delta\varphi^-_m, \Phiphi[m]^+) \bigr] \\
    &{} + \eta_0 (\delta u^-_0, \Phiu[0]^-)
      + \eta(\delta\varphi^-_0, \Phiphi[0]^-) \\
    &+ \sum_{m=1}^M
    a'_q(q(t_m),\boldu(t_m))(\delta q(t_m), \boldPhi(t_m)) \dtm
    \qfor \boldPhi \in X^0_k
    .
  \end{aligned}
\end{equation}
It is clear that the first sum is zero due to the dG(0) property.
By shifting the index of the third sum in \eqref{Tangent}
and applying the dG(0) property to $\Phiphi[m-1]^+$
we can combine the last three sums and rewrite \eqref{Tangent} as
\begin{equation}
  \label{Tangent_reform}
  \begin{aligned}
    0
    &= \sum_{m=1}^M \Bigl(
    a'_{\boldu}(q(t_m),\boldu(t_m))(\bolddeltau(t_m),\boldPhi(t_m))
    \dtm \\[-2\jot]
    &\mkern70mu+ \gamma (\delta\varphi^+_{m-1} - \delta\varphi^-_{m-1},
    \Phiphi[m]^-)\subphim{m}{m-1} \\
    &\mkern70mu+
    \eta (\delta\varphi^+_{m-1} - \delta\varphi^-_{m-1}, \Phiphi[m]^-) \\
    &\mkern70mu+
    a'_q(q(t_m),\boldu(t_m))(\delta q(t_m), \boldPhi(t_m)) \dtm \Bigr) \\
    &{} + \eta_0 (\delta u^-_0, \Phiu[0]^-)
      + \eta(\delta\varphi^-_0, \Phiphi[0]^-)
    \qfor \boldPhi \in X^0_k
    .
  \end{aligned}
\end{equation}

\subsection{Tangent time-stepping schemes}
As in the state equation we first solve the initial conditions,
\[
\begin{aligned}
  (\delta u^-_0, \Phiu[0]^-) &= 0, \\
  (\delta\varphi^-_0, \Phiphi[0]^-) &= 0,
\end{aligned}
\]
in short
\begin{equation}
  \label{Tangent_t_0}
  (\bolddeltau(t_0), \boldPhi_0^-) = 0 \qfor \boldPhi_0^- \in V.
\end{equation}
Applying the $X_k^0$ property to $\delta \varphi_{m-1}^{+}$ we can finally solve for $m=1,\dots,M$
the following equation
\begin{equation}
  \label{Tangent_ste}
  \begin{aligned}
    \gamma(\delta\varphi^-_m, &\Phiphi[m]^-)\subphim{m}{m-1}
    + \eta (\delta\varphi^-_m, \Phiphi[m]^-) \\
    &+ a'_{\boldu}(q(t_m),\boldu(t_m))(\bolddeltau(t_m),\boldPhi(t_m)) \dtm \\
    &= \eta (\delta\varphi^-_{m-1}, \Phiphi[m]^-) +
    \gamma (\delta\varphi^-_{m-1}, \Phiphi[m]^-)\subphim{m}{m-1} \\
    &- a'_q(q(t_m),\boldu(t_m))(\delta q(t_m),\boldPhi(t_m)) \dtm
    \qfor \boldPhi \in X^0_k
    .
  \end{aligned}
\end{equation}

\subsection{Adjoint Hessian equation}
The third and last auxiliary equation is the adjoint Hessian equation.
In this equation we seek $\bolddeltaz = (\delta\zu, \delta\zphi) \in X^0_k$
such that for all $\boldPhi \in X_k^0$ the following equation holds true:
\begin{equation}
  \label{Adj_Hessian}
  \Lagr''_{q\boldu}(q,\boldu,\boldz)(\delta q,\boldPhi) +
  \Lagr''_{\bolduu}(q,\boldu,\boldz)(\bolddeltau,\boldPhi) +
  \Lagr''_{\boldzu}(q,\boldu,\boldz)(\bolddeltaz, \boldPhi) = 0
  .
\end{equation}
First we see that $\Lagr''_{q\boldu}(q,\boldu,\boldz)(\delta q,\boldPhi) = 0$
since $q$ and $\boldu$ are decoupled. The derivative of $a$ in
$\Lagr''_{\boldzu}(q,\boldu,\boldz)(\bolddeltaz, \boldPhi)$
is given by $a'_{\boldu}(q,\boldu)(\boldPhi,\bolddeltaz)$
due to the linearity of $\boldz$ in $a$.
However, a genuine second-order derivative of $a$ arises in
$\Lagr''_{\bolduu}(q,\boldu,\boldz)(\bolddeltau,\boldPhi)$:
\begin{equation}
  \label{a'_uu}
  \begin{aligned}
    a''_{\bolduu}(q,\boldu)(\bolddeltau,\boldPhi,\boldz)
    &= 2 \varphi \cdot (1-\kappa) \Phiphi \cdot (\C e(\delta u), e(\zu)) \\
    &+ 2 \delta\varphi \cdot (1-\kappa) (\C e(u), e(\zu)) \cdot \Phiphi \\
    &+ 2 \varphi \cdot (1-\kappa) (\C e(u), e(\zu)) \delta\varphi \\
    &+ 2 \varphi \cdot (1-\kappa) (\C e(\Phiu) : e(\delta u), \zphi) \\
    &+ 2 \delta\varphi \cdot (1-\kappa) (\C e(\Phiu) : e(u), \zphi) \\
    &+ 2 (\C e(\delta u) : e(u), \zphi) \cdot \Phiphi
    .
  \end{aligned}
\end{equation}
Now we can rewrite \eqref{Adj_Hessian} in a dG(0) setting:
\begin{equation}
  \begin{aligned}
    0
    &= \sum_{m=1}^M
    J''_{\bolduu}(q(t_m),\boldu(t_m))(\bolddeltau(t_m),\boldPhi(t_m)) \\
    &- \sum_{m=1}^M a''_{\bolduu}(q(t_m),\boldu(t_m))
    (\bolddeltau(t_m),\boldPhi(t_m),\boldz(t_m)) \dtm \\
    &+ \sum_{m=1}^M \bigl[ \gamma
    (\Phiphi[m]^-, \delta\zphi[m]^- - \delta\zphi[m-1]^+)\subphim{m}{m-1}
    + \eta (\Phiphi[m]^-, \zphi[m]^- - \zphi[m-1]^+) \bigr] \\
    &- \gamma (\Phiphi[M]^-, \delta\zphi[M]^-)\subphim{M}{M-1}
    - \eta (\Phiphi[M]^-, \delta\zphi[M]^-) \\
    &+ \gamma (\Phiphi[0]^-, \delta\zphi[0]^-)\subphim{1}{0}
    + \eta(\Phiphi[0]^-, \delta\zphi[0]^-) \\
    &- \sum_{m=1}^M
    a'_{\boldu}(q(t_m),\boldu(t_m))(\boldPhi(t_m),\bolddeltaz(t_m)) \dtm \\
    &+\sum_{m=0}^{M-1} \gamma
    (\Phiphi[m]^-, \delta\zphi[m]^+ - \delta\zphi[m]^-)\subphim{m+1}{m}
    + \eta (\Phiphi[m]^-, \delta\zphi[m]^+ - \delta\zphi[m]^-) \\
    &- \eta_0 (\Phiu[0]^-, \delta\zu[0]^-)
    - \eta (\Phiphi[0]^-, \delta\zphi[0]^-)
    \qfor \boldPhi \in X^0_k
    .
  \end{aligned}
\end{equation}
Note that the same scaling of initial data was applied
that we already used for the adjoint equation.
By the $X_k^0$ property the third sum vanishes entirely.
Due to \cref{Remark_Projection_of_initial_solution} and the cancelation of
$\pm \eta(\Phiphi[0]^-, \delta\zphi[0]^-)$
the fifth line vanishes as well.
By shifting the index of the jump terms we can rewrite the equation as:
\begin{equation}
  \label{Adj_Hessian_1}
  \begin{aligned}
    0
    &= \sum_{m=1}^M \Bigl(
    J''_{\bolduu}(q(t_m),\boldu(t_m))(\bolddeltau(t_m),\boldPhi(t_m)) \\[-2\jot]
    &\mkern70mu- a''_{\bolduu}(q(t_m),\boldu(t_m))
    (\bolddeltau(t_m),\boldPhi(t_m),\boldz(t_m)) \dtm \\
    &\mkern70mu-
    a'_{\boldu}(q(t_m),\boldu(t_m))(\boldPhi(t_m),\bolddeltaz(t_m)) \dtm \\
    &\mkern70mu+ \gamma
    (\Phiphi[m-1]^-, \delta\zphi[m-1]^+ - \delta\zphi[m-1]^-)\subphim{m}{m-1} \\
    &\mkern70mu+ \eta
    (\Phiphi[m-1]^-, \delta\zphi[m-1]^+ - \delta\zphi[m-1]^-) \Bigr) \\
    &- \gamma (\Phiphi[M]^-, \delta\zphi[M]^-)\subphim{M}{M-1}
    - \eta(\Phiphi[M]^-, \delta\zphi[M]^-) \\
    &- \eta_0 (\Phiu[0]^-, \delta\zu[0]^-)
    \qfor \boldPhi \in X^0_k
    .
  \end{aligned}
\end{equation}
\subsection{{Adjoint Hessian time-stepping schemes}}
As in the adjoint time-step\-ping scheme
we first collect all terms that contain the last time point $t_M$ and solve
\begin{equation}
\label{Adj_Hessian_t_M}
  \begin{aligned}
    0
    &= J''_{\bolduu}(q(t_M)\boldu(t_M))(\bolddeltau(t_M), \boldPhi(t_M)) \\
    &- a'_{\boldu}(q(t_M)\boldu(t_M))(\boldPhi(t_M),\bolddeltaz(t_M)) \dtm[M] \\
    &- a''_{\bolduu}(q(t_M)\boldu(t_M))
    (\bolddeltau(t_M),\boldPhi(t_M),\boldz(t_M)) \dtm[M] \\
    &- \gamma (\Phiphi[M]^-, \delta\zphi[M]^-)\subphim{M}{M-1}
    - \eta (\Phiphi[M]^-, \delta\zphi[M]^-)
    \qfor \boldPhi \in X^0_k
    .
  \end{aligned}
\end{equation}
Then \eqref{Adj_Hessian_1} becomes
\begin{equation}
  \label{Adj_Hessian_ste}
  \begin{aligned}
    0
    &= \sum_{m=1}^{M-1} \Bigl(
    J''_{\bolduu}(q(t_m),\boldu(t_m))(\bolddeltau(t_m),\boldPhi(t_m)) \\[-2\jot]
    &\mkern70mu- a''_{\bolduu}(q(t_m),\boldu(t_m))
    (\bolddeltau(t_m),\boldPhi(t_m),\boldz(t_m)) \dtm \\
    &\mkern70mu- a'_{\boldu}(q(t_m),\boldu(t_m))(\boldPhi(t_m),\bolddeltaz(t_m))
    \dtm \Bigl) \\
    &+ \sum_{m=1}^M \Bigl( \gamma
    (\Phiphi[m-1]^-, \delta\zphi[m-1]^+ - \delta\zphi[m-1]^-)\subphim{m}{m-1}
    \\[-2\jot]
    &\mkern70mu+ \eta
    (\Phiphi[m-1]^-, \delta\zphi[m-1]^+ - \delta\zphi[m-1]^-) \Bigr) \\
    &- \eta_0 (\Phiu[0]^-, \delta\zu[0]^-)
    \qfor \boldPhi \in X^0_k
    .
  \end{aligned}
\end{equation}
In the final reformulation we shift the index of the second sum (jump-terms)
and take out the terms corresponding to $m = 0$
\begin{equation}
  \label{Adj_Hessian_t_0}
  \begin{aligned}
    0
    &= \sum_{m=1}^{M-1} \Bigl(
    J''_{\bolduu}(q(t_m),\boldu(t_m))(\bolddeltau(t_m),\boldPhi(t_m)) \\[-2\jot]
    &\mkern70mu- a''_{\bolduu}(q(t_m),\boldu(t_m))
    (\bolddeltau(t_m),\boldPhi(t_m),\boldz(t_m)) \dtm \\
    &\mkern70mu-
    a'_{\boldu}(q(t_m),\boldu(t_m))(\boldPhi(t_m),\bolddeltaz(t_m)) \dtm \\
    &\mkern70mu+ \gamma
    (\Phiphi[m]^-, \delta\zphi[m]^+ - \delta\zphi[m]^-)\subphim{m+1}{m} \\
    &\mkern70mu+ \eta
    (\Phiphi[m]^-, \delta\zphi[m]^+ - \delta\zphi[m]^-) \Bigr) \\
    &+ \gamma (\Phiphi[0]^-, \delta\zphi[0]^+ - \delta\zphi[0]^-)\subphim{1}{0}
    + \eta (\Phiphi[0]^-, \delta\zphi[0]^+ - \delta\zphi[0]^-) \\
    &- \eta_0 (\Phiu[0]^-, \delta\zu[0]^-)
    \qfor \boldPhi \in X^0_k
    .
  \end{aligned}
\end{equation}
As already pointed out in the time-stepping scheme for the adjoint equation,
all dual equations have to be solved backwards in time.
Therefore, we solve the following equation for $m=M-1,M-2,\dots,1$
\begin{equation*}
  \begin{aligned}
    0
    &= J''_{\bolduu}(q(t_m),\boldu(t_m))(\bolddeltau(t_m),\boldPhi(t_m)) \\
    &- a''_{\bolduu}(q(t_m),\boldu(t_m))
    (\bolddeltau(t_m),\boldPhi(t_m),\boldz(t_m)) \dtm \\
    &- a'_{\boldu}(q(t_m),\boldu(t_m))(\boldPhi(t_m),\bolddeltaz(t_m)) \dtm \\
    &+ \gamma
    (\Phiphi[m]^-, \delta\zphi[m]^+ - \delta\zphi[m]^-)\subphim{m+1}{m} \\
    &+ \eta (\Phiphi[m]^-, \delta\zphi[m]^+ - \delta\zphi[m]^-)
    \qfor \boldPhi \in X^0_k
    .
  \end{aligned}
\end{equation*}
As a result, the only remaning terms in \eqref{Adj_Hessian_t_0} are
\begin{equation}\label{Adj_hessian_before_drop}
  \gamma (\Phiphi[0]^-, \delta\zphi[0]^+ - \delta\zphi[0]^-)\subphim{1}{0}
  + \eta (\Phiphi[0]^-, \delta\zphi[0]^+ - \delta\zphi[0]^-)
  - \eta_0 (\Phiu[0]^-, \delta\zu[0]^-)
  .
\end{equation}
Finally we can apply the assumption $\eta_0 \ll \eta$ once more and drop the last term
in \eqref{Adj_hessian_before_drop}. Consequently the following equations have to be solved
for all $ \boldPhi \in X^0_k$:
\begin{align*}
  (\Phiphi[0]^-, \delta\zphi[0]^-) &= (\Phiphi[0]^-, \delta\zphi[1]^-), \\
  (\Phiu[0]^-, \delta\zu[0]^-) &= (\Phiu[0]^-, \delta\zu[1]^-).
\end{align*}
Note that \cref{Remark_Algorithmic_realization} was applied
to \eqref{Adj_hessian_before_drop} as well.

\subsection{Final complete algorithm}
\label{SEC:Algorithm}

Gathering the optimization problem statement and the space-time discretizations from the previous sections and resulting time-stepping schemes for the
four equations yields the complete method given in Algorithm~\ref{alg:main}.

\newcommand\boldf{\bs f}
\newcommand\boldd{\bs d}
\newcommand\boldG{\bs G}
\newcommand\boldH{\bs H}
\newcommand\boldh{\bs h}
\begin{algorithm}[tp]
  \vspace{2pt}
  \KwData{Domain $\Omega$, mesh $\mathcal{T}_h$, number of time intervals $M$,
    parameters $\varepsilon$, $\kappa$, $G_c$, $\mu$, $\lambda$, $\gamma$,
    $\eta$, $\alpha$, initial value~$\boldu_0$, initial control guess $q^0$.}
  \KwResult{Optimal control $q$ and admissible solution $\boldu$.}
  \vspace{2pt}
  1: Set $k = 0$ and $q^k = q^0$ and solve the
  state equation for $\boldu$:
  $\Lagr'_{\boldz}(q^k,\boldu,\boldz)(\boldPhi) = 0 \tfor \boldPhi$.
  Specifically, obtain $\boldu(t_0)$ from \eqref{state_t=0} and then
  $\boldu(t_1),\dots,\boldu(t_M)$ from \eqref{state_t>0}\;
  2: Solve the adjoint equation for $\boldz$:
  $\Lagr'_{\boldu}(q^k,\boldu,\boldz)(\boldPhi) = 0 \tfor \boldPhi$.
  Obtain $\boldz(t_M)$ from \eqref{Adjoint_t_M},
  then $\boldz(t_{M-1}),\dots,\boldz(t_1)$ from \eqref{Adjoint_ste},
  and finally $\boldz(t_0)$ from \eqref{Adjoint_t_0}\;
  3: Construct the coefficient vector $\boldf \in \R^n$
  for the reduced gradient $\nabla j(q^k)$ by solving
  $\boldG \boldf = [j'(q^k)(q_i)]_{i=1}^n$.
  Here $q_i$ denotes the $i$-th basis function
  of the discrete control space $Q_h$
  and $\boldG_{ij} = (q_i,q_j)$ defines the mass matrix.
  The derivatives $j'(q^k)(q_i)$ for the right hand side
  are computed from the representation \eqref{Representations}\;
  \While{$\norm{\boldf}_2 > TOL$}
  {
    4: Obtain $\delta q$ from the Newton equation,
    $j''(q^k)(\delta q, q_i) = -j'(q^k)(q_i) \tfor q_i$,
    by minimizing $m(q^k,\boldd) = j(q^k) + \sprod{\boldf}{\boldd}
    + \frac12 \sprod{\boldH \boldd}{\boldd}$
    for a vector $\boldd \in \R^n$ using the CG-method (matrix free).
    Here $\boldH \in \R^{n \times n}$ denotes the coefficient matrix
    of $\nabla^2 j(q^k) \delta q$\;
    \For{every CG step}
    {
      5: Solve the tangent equation for $\bolddeltau$:
      $\Lagr''_{q\boldz}(q^k,\boldu,\boldz)(\delta q, \boldPhi) +
      \Lagr''_{\boldu\boldz}(q^k,\boldu,\boldz)(\bolddeltau,\boldPhi) = 0
      \tfor \boldPhi$.
      Obtain $\bolddeltau(t_0)$ from \eqref{Tangent_t_0} and then
      $\bolddeltau(t_1),\dots,\bolddeltau(t_M)$ from \eqref{Tangent_ste}\;
      6: Solve the adjoint Hessian equation for $\bolddeltaz$:
      $\Lagr''_{q\boldu}(q^k,\boldu,\boldz)(\delta q,\boldPhi) +
      \Lagr''_{\bolduu}(q^k,\boldu,\boldz)(\bolddeltau,\boldPhi) +
      \Lagr''_{\boldzu}(q^k,\boldu,\boldz)(\bolddeltaz,\boldPhi) = 0
      \tfor \boldPhi$.
      Obtain $\bolddeltaz(t_M)$ from \eqref{Adj_Hessian_t_M}, then
      $\bolddeltaz(t_{M-1}),\dots,\bolddeltaz(t_1)$ from \eqref{Adj_Hessian_ste},
      and finally $\bolddeltaz(t_0)$ from \eqref{Adj_Hessian_t_0}\;
      7: Construct the coefficient vector $\boldh \in \R^n$
      for $\nabla^2 j(q^k)\delta q$
      by solving $\boldG \boldh = j''(q^k)(\delta q, q_i)_{i=1}^n$,
      where $j''(q^k)(\delta q, q_i)$ is represented
      via \eqref{Representations}\;
    }
    8: Choose a step length $\nu$ by an Armijo backtracking method\;
    9: Set $q^{k+1} = q^k + \nu \delta q$\;
    10: Repeat steps 1, 2, 3 for the new control $q^{k+1}$
    to obtain $\boldf$ for $\nabla j(q^{k+1})$.\;
    11: Increment $k=k+1$\;
  }
  \caption{Overall space-time phase-field fracture control algorithm}
  \label{alg:main}
\end{algorithm}

\section{Numerical studies}
\label{SEC:Numtests}
In this section we present six numerical experiments.
In these experiments we use the tracking type functional of \eqref{EQ:NLPgamma}
to find an optimal control force
that approximately produces a desired phase-field.
All numerical computations are performed with the open source software libraries
\texttt{deal.II} \cite{dealII91,deal2020}
and \textsc{DOpElib} \cite{dope,DOpElib}.

Since a large body of the published literature deals with forward phase-field 
fracture in which cracks propagate through large parts of the domain or even until domain 
boundaries, we emphasize that on purpose, short fractures are considered 
in our optimal control settings only. Specifically, $\varphi_d$ is prescribed 
sufficiently small such that we clearly can distinguish 
between our optimal control final fractures 
and classical non-controlled fractures.

\subsection{Experiment 1: horizontal fracture in right half domain}
\label{SEC:Experiment1}

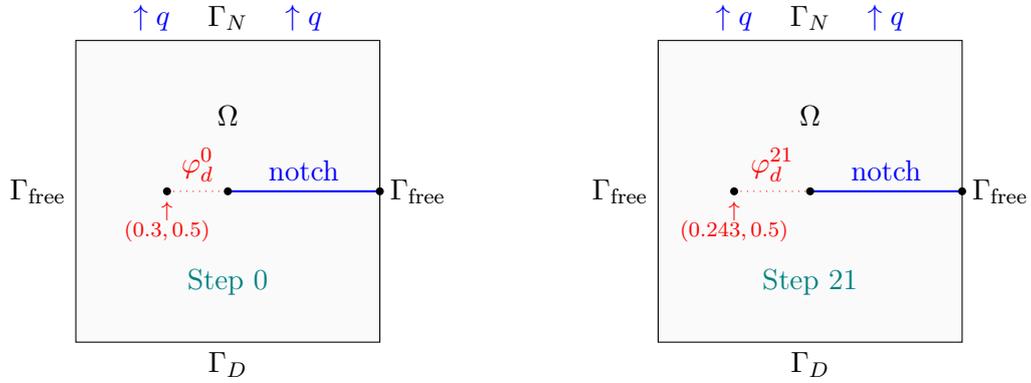
\begin{figure}[tp]
  \centering
  \begin{tikzpicture}
    \draw[fill=black!2] (0,0)
    -- (0,4) node[pos=0.5,left]  {$\Gamma\free$}
    -- (4,4) node[pos=0.5,above] {$\Gamma_N$}
    -- (4,0) node[pos=0.5,right] {$\Gamma\free$}
    -- cycle node[pos=0.5,below] {$\Gamma_D$};
    \node[teal] at (2,0.8) {Step 0};
    \node at (2,3) {$\Omega$};
    \node[blue,above] at (1,4) {$\uparrow q$};
    \node[blue,above] at (3,4) {$\uparrow q$};
    \draw[dotted,red] (1.2,2) -- (2,2) node[pos=0.5,above,red] {$\varphi^0_d$};
    \draw[thick,blue] (2,2) -- (4,2) node[pos=0.5,above,blue] {notch};
    \node[red,below] at (1.2,2) {\scriptsize{$\uparrow$}};
    \node[red,below] at (1.2,1.75) {\scriptsize{$(0.3,0.5)$}};
    \fill (1.2,2) circle (1.5pt) (2,2) circle (1.5pt) (4,2) circle (1.5pt);
  \end{tikzpicture}\hfil
  \begin{tikzpicture}
    \draw[fill=black!2] (0,0)
    -- (0,4) node[pos=0.5,left]  {$\Gamma\free$}
    -- (4,4) node[pos=0.5,above] {$\Gamma_N$}
    -- (4,0) node[pos=0.5,right] {$\Gamma\free$}
    -- cycle node[pos=0.5,below] {$\Gamma_D$};
    \node[teal] at (2,0.8) {Step 21};
    \node at (2,3) {$\Omega$};
    \node[blue,above] at (1,4) {$\uparrow q$};
    \node[blue,above] at (3,4) {$\uparrow q$};
    \draw[dotted,red]
    (1.0,2) -- (2,2) node[pos=0.5,above,red] {$\varphi^{21}_d$};
    \draw[thick,blue] (2,2) -- (4,2) node[pos=0.5,above,blue] {notch};
    \node[red,below] at (1.0,2) {\scriptsize{$\uparrow$}};
    \node[red,below] at (1.0,1.75) {\scriptsize{$(0.243,0.5)$}};
    \fill (1.0,2) circle (1.5pt) (2,2) circle (1.5pt) (4,2) circle (1.5pt);
  \end{tikzpicture}
  \caption{Experiment 1:
    domain $\Omega = (0, 1)^2$ with partitioned boundary $\partial \Omega$,
    initial notch and desired crack~$\varphi_d$
    for homotopy steps 0 and 21 in approach~(a).}
  \label{fig:EX_1_experiments}
\end{figure}

The first experiment is motivated by a standard problem:
the single edge notched tension test \cite{MieWelHof10a,MieWelHof10b}.
Here we consider the square domain $\Omega = (0, 1)^2$
with a horizontal notch, see \cref{fig:EX_1_experiments}.
The notch is in the middle of the right side of the domain,
defined as $(0.5, 1) \times \set{0.5}$.
The boundary $\partial \Omega$ is partitioned as
$\partial \Omega = \Gamma_N \cup \Gamma_D \cup \Gamma\free$,
where $\Gamma_N \coloneqq [0,1] \times \set{1}$,
$\Gamma_D \coloneqq [0, 1] \times \set{0}$, and
$\Gamma\free \coloneqq \set{0, 1} \times (0, 1)$.
On $\Gamma_N$ we apply the force $q$ in orthogonal direction to the domain,
on $\Gamma_D$ we enforce homogeneous Dirichlet boundary conditions for the
displacement $u=0$, and on $\Gamma\free$ we set homogeneous Neumann
boundary conditions.
We choose the time interval $[0, 1]$ with 41 equidistant time points $t_m$,
i.e., $T = 1$ and $M = 40$.
The discrete control space $Q_h$
is one-dimensional in the sense that the force is constant in time
and is only applied in $y$ direction.
The spatial mesh consists of $64 \times 64$ square elements,
hence the element diameter is $h = \sqrt2/64 \approx 0.0221$.
The initial values are given by $\boldu_0 = (u_0, \varphi_0)$
where $\varphi_0$ describes the horizontal notch:
\begin{equation*}
  \label{EQ:EX_1_phi0}
  \varphi_0(x, y) \coloneqq
  \begin{cases}
    0, & x \in (0.50, 1.00) \text{ and } y = 0.5 , \\
    1, & \text{else}.
  \end{cases}
\end{equation*}
The desired phase-field $\varphi_d$ is defined as a continuation
of the initial notch to the left hand side of the domain,
see $\varphi_d^0$ in \cref{fig:EX_1_experiments}.
In order to investigate the effect of $\varphi_d$ on the optimal solution,
we will use two different homotopy approaches.
In approach (a) we will successively increase
the length of the desired phase-field,
and in approach (b) we will successively reduce
the Tikhonov parameter $\alpha$.
In both cases the motivation is to increase the weight
of the physically motivated term $\frac{1}{2}\norm{\varphi-\varphi_d}^2$
in relation to the Tikhonov term.
We will perform as many homotopy steps as possible,
solving one NLP per step.
The common nominal parameters used in both approaches
are given in \cref{tab:EX_1_params}.
\newcommand\captionpar[1]{\caption{Experiment #1:
    regularization and penalty parameters (left),
    model and material parameters (right).}}
\begin{table}[tp]
  \captionpar{1}
  \centering
  \begin{tabular}{llS[table-format=1.2e+2]}
    \toprule
    Par. & Definition & {Value} \\
    \midrule
    $\varepsilon$ & Regul. (crack) $\approx 4 h$ & 0.0884 \\
    $\kappa$ & Regul. (crack) & 1.00e-10 \\
    $\eta$ & Regul. (viscosity) & 1.00e3 \\
    $\gamma$ & Penalty & 1.00e5  \\
    $\alpha$ & Tikhonov & 4.75e-10 \\
    \bottomrule
  \end{tabular}
  \hfil
  \begin{tabular}{llS[table-format=1.1e1]}
    \toprule
    Par. & Definition & {Value} \\
    \midrule
    $G_c$ & Fracture toughness & 1.0 \\
    $\nu_s$ & Poisson's ratio & 0.2 \\
    $E$ & Young's modulus & 1.0e6 \\
    $q_0$ & Initial control & 1.0 \\
    $q_d$ & Nominal control & 1.0e3 \\
    \bottomrule
  \end{tabular}
  \label{tab:EX_1_params}
\end{table}

\subsubsection{Approach (a): length increment}

Here we will solve the NLP \eqref{EQ:NLPgamma} several times
with different desired phase-fields $\varphi^k_d$.
Formally we define a sequence of desired phase-fields $\varphi^k_d$
with $\varphi^{k+1}_d < \varphi^k_d$
where $\varphi^k_d$ is defined as
\begin{equation*}
  \label{EQ:EX_1_phid}
  \varphi^k_d(x, y) \coloneqq
  \begin{cases}
    0, & x \in (0.3 \times 0.99^k, 0.5) \text{ and }
    y \in (0.5 - h, 0.5 + h), \\
    1, & \text{else}.
  \end{cases}
\end{equation*}
By this definition of $\varphi_d^k$ we extend
  the desired crack to the left so that it becomes gradually longer.
The number of homotopy steps performed in this experiment is 21:
in step 22, the iterative solution of the nonlinear state equation fails
because the Newton residuals do not decrease towards zero.
Probably this means that large numerical errors
  prevent finding a descent direction or that the initial estimate
  lies outside the area of convergence, but various other reasons
  might be possible as well.
Our results are presented in
\cref{tab:EX_1_results1}.
The first column (Step) counts the homotopy steps.
The second column (Iter) gives the number of Newton iterations
for solving the associated reduced problem \eqref{EQ:NLPgammared},
except that Iter 0 in Step 0 refers to the
initial guess from which the homotopy starts.
The remaining values are the absolute Newton residual,
the cost functional $\J$ and its tracking part
$\frac12 \sum_m \norm{\varphi(t_m) - \varphi_d(t_m)}^2$,
the maximal force $\abs{q_{\max}}$ applied on $\Gamma_N$,
and the Tikhonov regularization term,
$\frac{\alpha}{2}\sum_m \norm{q(t_m) - q_d(t_m)}_{\Gamma_N}^2$.
All values are rounded to three, five, or six significant digits.
For every NLP the Newton iteration terminates
when the residual falls below the tolerance \num{5e-11}.

\newcommand\captionperfhom[2][]{\caption{Experiment #2:
    number of Newton iterations, absolute residual, cost terms and maximal
    force during homotopy. Iter 0 in step 0 refers to initial state from
    which homotopy starts.\ifx#1\empty\else\ #1.\fi}}
\newcommand\captionhom[3]{\captionperf[Part #2: homotopy steps #3]{#1}}
\newcommand\Hstep[2][]{&&&&\multicolumn1c{\text{Step#1 #2}}\\[1pt]}
\begin{table}[tp]
  \centering
  \captionperfhom{1a}
  \sisetup{table-format=1.4e+1}
  \begin{tabular}{rrS[table-format=1.2e+2]SSSS[table-format=4.2]}
    \toprule
    Step&Iter&{Residual}&{Cost}&{Tracking}&{Tikhonov}&{Force}\\
    \midrule
    0 & 0 & 4.62e-07 & 4.1532e-3 & 3.9192e-3 & 2.3406e-4  & 1.0    \\
    0 & 9 & 2.62e-11 & 3.4863e-3 & 3.3648e-3 & 1.2150e-4  & 2379.02\\
    1 & 3 & 2.70e-11 & 3.5681e-3 & 3.4447e-3 & 1.2337e-4  & 2388.20\\
    2 & 4 & 1.83e-11 & 3.6503e-3 & 3.5254e-3 & 1.2489e-4  & 2395.58\\
    3 & 0 & 1.83e-11 & 3.6503e-3 & 3.5254e-3 & 1.2489e-4  & 2395.58\\
    4 & 4 & 2.76e-11 & 3.7823e-3 & 3.6552e-3 & 1.2708e-4  & 2405.78\\
    5 & 0 & 2.76e-11 & 3.7822e-3 & 3.6552e-3 & 1.2708e-4  & 2405.78\\
    6 & 2 & 2.37e-11 & 3.8645e-3 & 3.7357e-3 & 1.2883e-4  & 2414.00\\
    7 & 0 & 2.38e-11 & 3.8645e-3 & 3.7357e-3 & 1.2883e-4  & 2414.00\\
    8 & 2 & 1.29e-11 & 3.9466e-3 & 3.8156e-3 & 1.3099e-4  & 2424.06\\
    9 & 0 & 1.28e-11 & 3.9466e-3 & 3.8156e-3 & 1.3099e-4  & 2424.06\\
    10& 2 & 4.28e-11 & 4.0779e-3 & 3.9433e-3 & 1.3465e-4  & 2439.95\\
    11& 0 & 4.28e-11 & 4.0779e-3 & 3.9433e-3 & 1.3465e-4  & 2439.95\\
    12& 2 & 2.77e-11 & 4.1605e-3 & 4.0241e-3 & 1.3640e-4  & 2447.20\\
    13& 6 & 4.03e-11 & 4.2438e-3 & 4.1063e-3 & 1.3748e-4  & 2451.38\\
    14& 0 & 4.03e-11 & 4.2438e-3 & 4.1063e-3 & 1.3748e-3  & 2451.38\\
    15& 0 & 4.03e-11 & 4.2438e-3 & 4.1063e-3 & 1.3748e-3  & 2451.38\\
    16& 2 & 4.43e-11 & 4.3760e-3 & 4.2356e-3 & 1.4043e-4  & 2464.78\\
    17& 0 & 4.43e-11 & 4.3760e-3 & 4.2356e-3 & 1.4043e-4  & 2464.78\\
    18& 2 & 3.82e-11 & 4.4585e-3 & 4.3160e-3 & 1.4258e-4  & 2472.69\\
    19& 7 & 4.98e-11 & 4.5360e-3 & 4.3858e-3 & 1.5012e-4  & 2498.67\\
    20& 0 & 4.98e-11 & 4.5360e-3 & 4.3858e-3 & 1.5012e-4  & 2498.67\\
    21& 0 & 4.98e-11 & 4.5360e-3 & 4.3858e-3 & 1.5012e-4  & 2498.67\\
    \bottomrule
  \end{tabular}
  \label{tab:EX_1_results1}
\end{table}
\begin{figure}[tp]
  \footnotesize
  \centering
  \begin{tikzpicture}
    \begin{axis}[width=1.0\linewidth,height=50mm,grid=major,
      xmin=-2,xmax=68,ymin=3.2e-3,ymax=4.6e-3,ytick distance={2e-4}]
      \addplot[red!40] table[x expr=\coordindex,y index=5]
      {TXT11_Data.txt};
      \addplot[ycomb,red] table[x expr=\coordindex,y index=5]
      {TXT11_Data.txt};
      \addplot[blue!40] table[x expr=\coordindex,y index=6]
      {TXT11_Data.txt};
      \addplot[ycomb,blue] table[x expr=\coordindex,y index=6]
      {TXT11_Data.txt};
    \end{axis}
  \end{tikzpicture}
  \caption{Experiment 1a: cost functional of each NLP iteration in homotopy
    (blue: tracking part above \num{3.2e-3} + red: Tikhonov part).}
  \label{fig:EX1_1_Cost}
\end{figure}
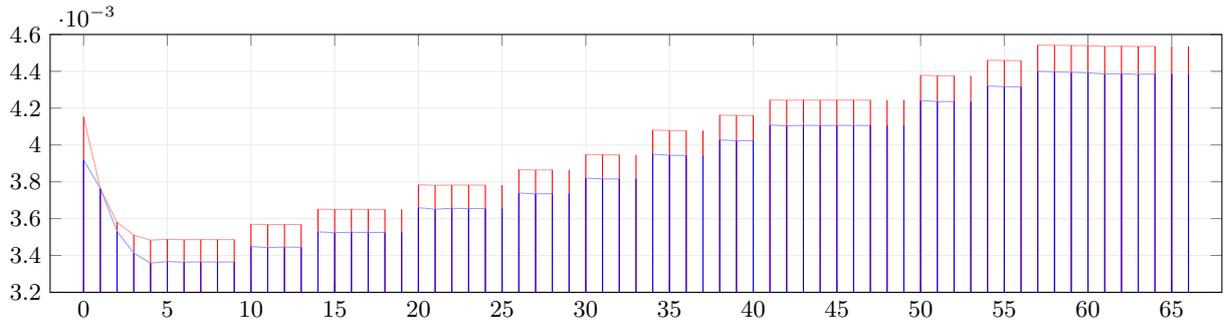
\begin{figure}[tp]
  \footnotesize
  \centering
  \begin{tikzpicture}
    \begin{axis}[width=1.0\linewidth,height=50mm,grid=major,
      xmin=-2,xmax=68,ymin=1.0e-11,ymax=1e-6,ytick distance={10},ymode=log]
      \addplot[blue,mark=*,mark size=1pt] table[x expr=\coordindex,y index=4]
      {TXT11_Data.txt};
    \end{axis}
  \end{tikzpicture}
  \caption{Experiment 1a: absolute residual of each NLP iteration in homotopy.}
  \label{fig:EX1_1_Res}
\end{figure}
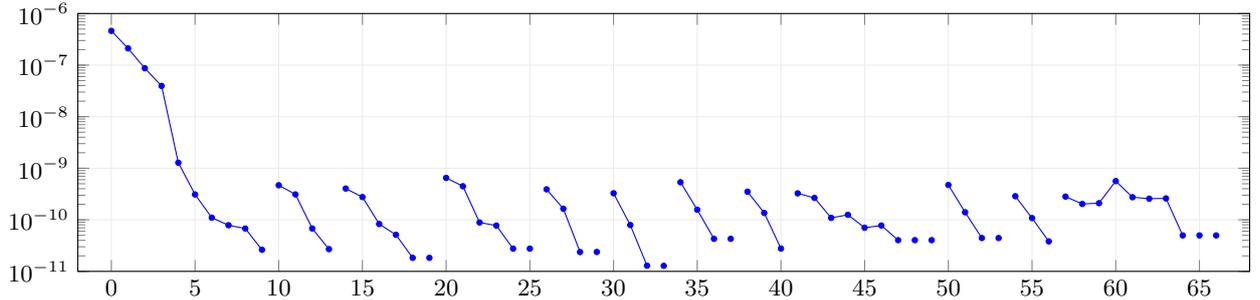

\subsubsection{Approach (b): Tikhonov iteration}

\begin{figure}[tp]
  \footnotesize
  \centering
  \begin{tikzpicture}
    \begin{axis}[width=1.0\linewidth,height=50mm,grid=major,
      xmin=-2,xmax=51,ymin=1.0e-11,ymax=1e-6,ytick distance={10},ymode=log]
      \addplot[blue,mark=*,mark size=1pt] table[x expr=\coordindex,y index=4]
      {TXT12_Data.txt};
    \end{axis}
  \end{tikzpicture}
  \caption{Experiment 1b: absolute residual of each NLP iteration in homotopy.}
  \label{fig:EX1_2_Res}
\end{figure}
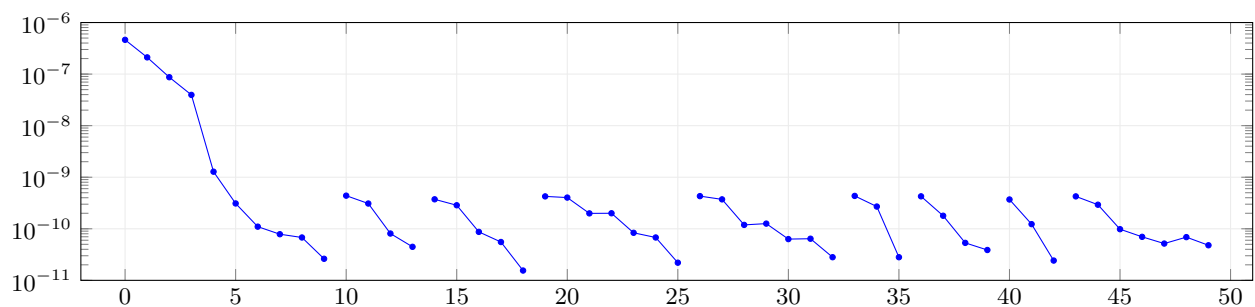

The second approach is a successive reduction of~$\alpha$,
a so called Tikhonov iteration.
In this case the length of the desired phase-field remains constant
at $\varphi^0_d$ for all homotopy steps
while the weight of the Tikhonov term in $\J$ is successively reduced.
Here we define the sequence by $\alpha_k = 0.99^k \alpha_0$
with $\alpha_0 = \num{4.75e-10}$.
The number of homotopy steps performed is 8:
in step 9, we have terminated the computation
because of very slow alternating convergence of the residual
as is often observed for very small values of the parameter~$\alpha$.
The results are presented in \cref{tab:EX_1_results2}.
First we notice the high sensitivity of our NLP
solution with respect to the control force.
In \cref{fig:EX1_sensitivity} we present the difference of
the controls on iterations 8 and~9 of the initial homotopy step.
A comparison of the corresponding residuals shows a reduction
from \num{6.78e-11} to \num{2.62e-11} (approximately 60\%),
even though the maximal difference between the applied controls
is only $2.1$ (or 0.1\%).
The values of the cost functional and the residual
on all iterations of both homotopies are presented in
\cref{fig:EX1_1_Cost,fig:EX1_1_Res,fig:EX1_2_Cost,fig:EX1_2_Res}.
Each dot in \cref{fig:EX1_1_Res} stands for one Newton step within
the corresponding homotopy iteration.
The behavior of the residual values in both approaches
is typical for homotopy methods:
in each homotopy step they are reduced below the tolerance,
and they increase slightly afterwards.
In the final homotopy step of each approach the reduction is non-monotonous
because the maximal number of line search iterations is reached;
this indicates the difficulty of the NLP.
In \cref{fig:EX1_1_Cost} we observe that the value of the
cost functional increases with each homotopy step.
This is a consequence of the increasing length of the desired
phase-field: $\varphi_d^{21} < \dots <{} \varphi_d^0$.
A closer look at the results reveals that the tracking part actually
increases non-linearly with the length of the desired phase-field,
which is not surprising as our overall problem is nonlinear.
Finally we observe that both approaches yield larger maximal
control forces when compared to the results without homotopy ansatz.
In approach (a) the maximal final control is $2498.67$,
and in approach (b) it is $2438.79$.
This corresponds to the different cracks being produced:
without any homotopy approach the crack has a total length of $0.063$,
with the Tikhonov iteration (approach b) we obtain $0.078$,
and with the crack length increment (approach a) we obtain $0.094$;
see \cref{fig:EX_1_PF}. We notice that we do not 
infer from these results any evaluation on which approach is better for 
this test case, but we can only say that both yield different findings.
\begin{figure}[tp]
  \footnotesize
  \centering
  \begin{tikzpicture}
    \begin{axis}[width=1.0\linewidth,height=50mm,grid=major,
      xmin=0,xmax=1.0,ymin=-0.5,ymax=2.5,ytick distance={0.5}]
      \addplot[blue] table
      {TXT12_Approach_2_Step_0_Iter_8_min_Iter_9.txt}
      node[pos=0.3,above] {$q_8 - q_9$};
    \end{axis}
  \end{tikzpicture}
  \caption{Experiment 1b:
    difference of control forces on iterations 8 and 9 in homotopy step 0.}
  \label{fig:EX1_sensitivity}
\end{figure}
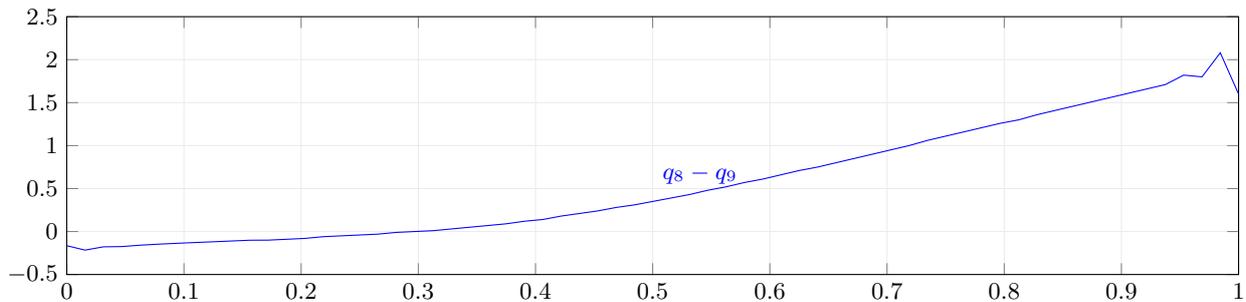
\begin{table}[tp]
  \centering
  \captionperfhom{1b}
  \sisetup{table-format=1.4e+1}
  \begin{tabular}{rrS[table-format=1.2e+2]SSSS[table-format=4.2]}
    \toprule
    Step&Iter&{Residual}&{Cost}&{Tracking}&{Tikhonov}&{Force}\\
    \midrule
    0 & 0 & 4.62e-07 & 4.1532e-3 & 3.9192e-3 & 2.3406e-4  & 1.0    \\
    0 & 9 & 2.62e-11 & 3.4863e-3 & 3.3648e-3 & 1.2150e-4  & 2379.02\\
    1 & 3 & 4.47e-11 & 3.4835e-3 & 3.3614e-3 & 1.2206e-4  & 2387.69\\
    2 & 4 & 1.55e-11 & 3.4812e-3 & 3.3590e-3 & 1.2212e-4  & 2393.56\\
    3 & 6 & 2.20e-11 & 3.4787e-3 & 3.3565e-3 & 1.2223e-4  & 2399.79\\
    4 & 6 & 2.81e-11 & 3.4764e-3 & 3.3541e-3 & 1.2226e-4  & 2405.43\\
    5 & 2 & 2.81e-11 & 3.4737e-3 & 3.3512e-3 & 1.2257e-4  & 2412.03\\
    6 & 3 & 3.88e-11 & 3.4706e-3 & 3.3469e-3 & 1.2362e-4  & 2423.04\\
    7 & 2 & 2.42e-11 & 3.4676e-3 & 3.3433e-3 & 1.2437e-4  & 2432.00\\
    8 & 6 & 4.79e-11 & 3.4648e-3 & 3.3310e-3 & 1.2487e-4  & 2438.79\\
    \bottomrule
  \end{tabular}
  \label{tab:EX_1_results2}
\end{table}
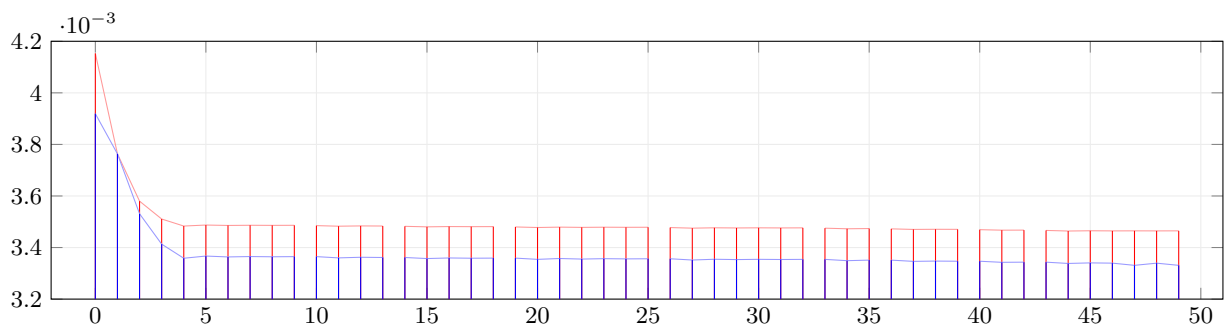
\begin{figure}[tp]
  \footnotesize
  \centering
  \begin{tikzpicture}
    \begin{axis}[width=1.0\linewidth,height=50mm,grid=major,
      xmin=-2,xmax=51,ymin=3.2e-3,ymax=4.2e-3,ytick distance={2e-4}]
      \addplot[red!40] table[x expr=\coordindex,y index=5]
      {TXT12_Data.txt};
      \addplot[ycomb,red] table[x expr=\coordindex,y index=5]
      {TXT12_Data.txt};
      \addplot[blue!40] table[x expr=\coordindex,y index=6]
      {TXT12_Data.txt};
      \addplot[ycomb,blue] table[x expr=\coordindex,y index=6]
      {TXT12_Data.txt};
    \end{axis}
  \end{tikzpicture}
  \caption{Experiment 1b: cost functional of each NLP iteration in homotopy
    (blue: tracking part above \num{3.2e-3} + red: Tikhonov part).}
  \label{fig:EX1_2_Cost}
\end{figure}

\newcommand\image[2][]{\bgroup\fboxsep0pt\fboxrule0pt
  \protect\fcolorbox{black}{black!10}{\protect\includegraphics[#1]{#2}}\egroup}

\newcommand\scale[5][0.1]{\parbox{#2\linewidth}{\vspace{0.3pc}%
    \makebox[\linewidth]{\num{#3}\hfill$#4$\hfill\num{#5}}\\[-#1pc]%
    \image[trim=225 1180 275 945,clip=true,width=1.0\linewidth]
    {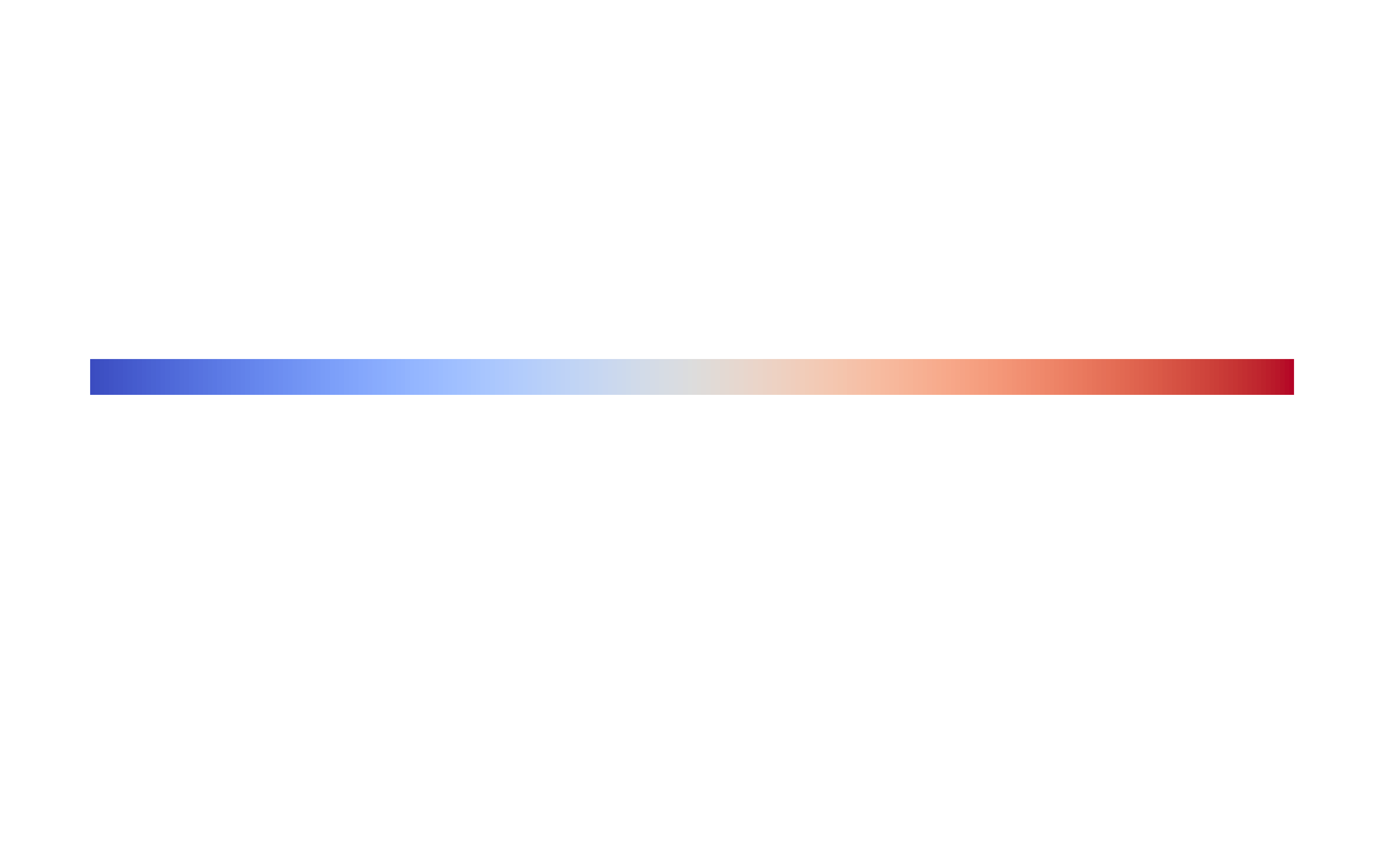}}%
}
\newcommand\img[3]{\image
  [trim=840 315 850 315,clip=true,width=.325\linewidth]
  {EX_1_AP_#1_UPD_#2_STATE_#3_TIME_40_PF.png}}
\begin{figure}[tp]
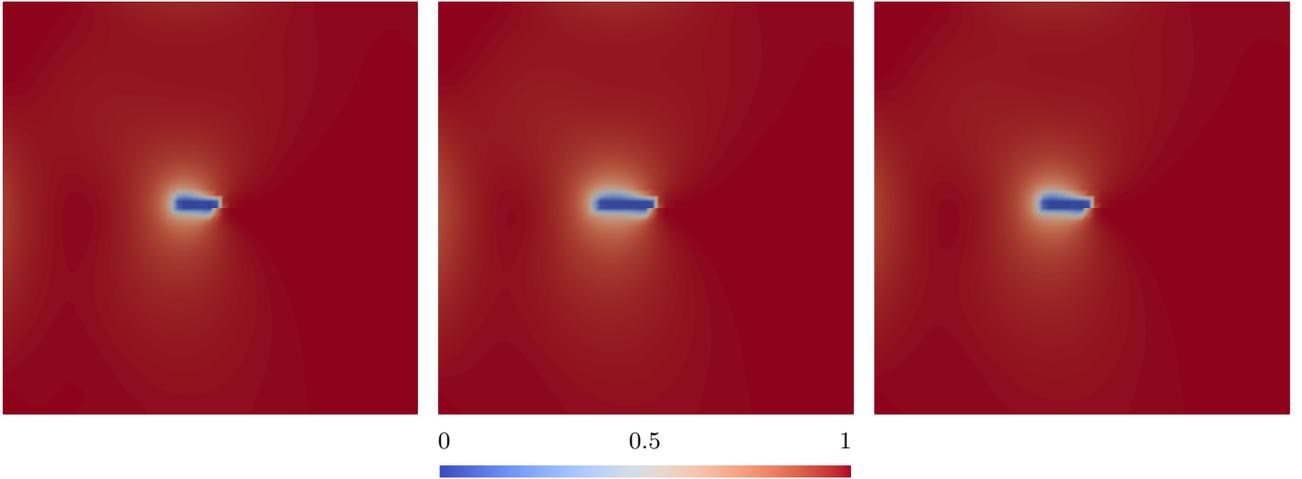

  \footnotesize
  \centering
  \img{1}{0}{9}\hfill
  \img{1}{19}{7}\hfill
  \img{2}{8}{6}\\%
  \scale[0.33]{0.32}{0}{0.5}{1}%
  \caption{Experiment 1: optimal phase-field $\varphi$ at time 40.
    Step 0 (left), step~19 of approach (a) (center, best result),
    step 8 of approach (b) (right).}
  \label{fig:EX_1_PF}
\end{figure}

\subsection{Experiment 2: two-sided control for diagonal crack}
\label{SEC:Experiment6}

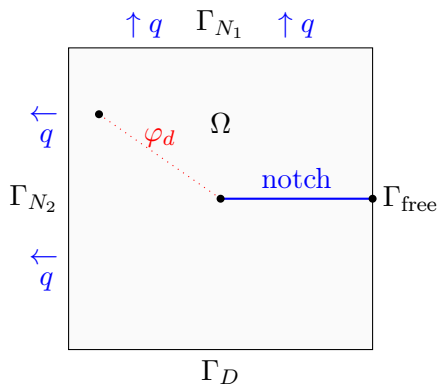
\begin{figure}
  \centering
 \begin{tikzpicture}
    \draw[fill=black!2] (0,0)
    -- (0,4) node[pos=0.5,left]  {$\Gamma_{N_2}$}
    -- (4,4) node[pos=0.5,above] {$\Gamma_{N_1}$}
    -- (4,0) node[pos=0.5,right] {$\Gamma\free$}
    -- cycle node[pos=0.5,below] {$\Gamma_D$};
    \node at (2,3) {$\Omega$};
    \node[blue,above] at (1,4) {$\uparrow q$};
    \node[blue,above] at (3,4) {$\uparrow q$};
    \node[blue,left] at (0,1.2) {$\leftarrow$};
    \node[blue,left] at (0,0.9) {$q\,$};
    \node[blue,left] at (0,3.1) {$\leftarrow$};
    \node[blue,left] at (0,2.8) {$q\,$};
   \draw[dotted,red] (0.4,3.12) -- (3-1,2) node[pos=0.5,above,red] {$\varphi_d$};
    \draw[thick,blue] (3-1,2) -- (5-1,2) node[pos=0.5,above,blue] {notch};
    \fill
    (0.4,3.12) circle (1.5pt) (3-1,2) circle (1.5pt) (5-1,2) circle (1.5pt);
  \end{tikzpicture}
  \caption{Experiment 2:
    domain $\Omega = (0,1)^2$ with partitioned boundary, initial notch
    and desired crack $\varphi_d$.}
  \label{fig:EX_6_experiments}
\end{figure}

\begin{figure}
  \footnotesize
  \centering
  \begin{tikzpicture}
    \begin{axis}[width=1.0\linewidth,height=8cm,
      xmin=0,xmax=1,ymin=1950,ymax=2650,ytick distance={100},grid=major]
      \draw[red,dotted] (0,2200) -- (1,2200) node[pos=0.52,above] {$q_d$};
      \addplot[blue] table
      {TXT6_EX_6_Force_TOP.txt} node[pos=0.534,above] {$q_1$};
      \addplot[blue] table
      {TXT6_EX_6_Force_LEFT.txt} node[pos=0.075,below] {$q_2$};
    \end{axis}
  \end{tikzpicture}
  \caption{Experiment 2: optimal control forces (solid:
    $q_1$ on upper boundary $\Gamma_{N_1}$,
    $q_2$ on left boundary $\Gamma_{N_2}$)
    and common nominal control force $q_d$ (dotted).}
  \label{fig:EX_6_Forces}
\end{figure}
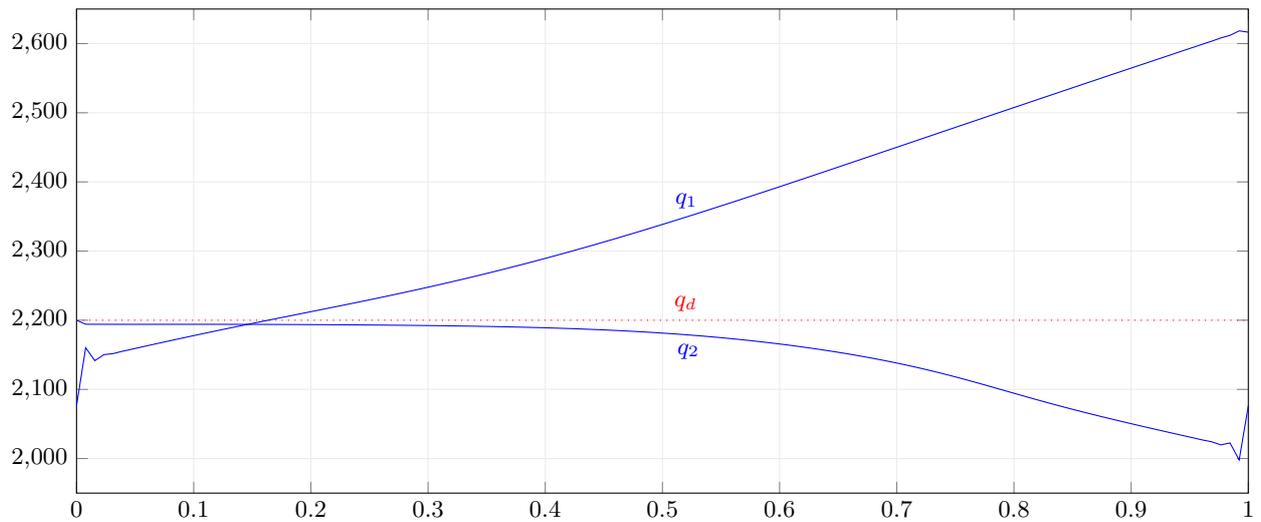

Our second experiment is an extension of the first one
with the aim to create a crack that grows diagonally,
in negative $x$ direction and positive $y$ direction.
We consider the same domain as before, $\Omega = (0,1)^2$,
but since cracks grow orthogonal to the maximum tensile
stress \cite[Chapter 4]{Zehnder2012},
the original control boundary $\Gamma_N = [0,1] \times \set{1}$
becomes $\Gamma_{N_1}$ and we extend the control
to a second boundary $\Gamma_{N_2} = \set{0} \times [0,1]$,
i.e., in the PDE constraint, we have
\[
  (q, \Phi_{u:\perp})_{\Gamma_N,I} =
  (q, \Phi_{u:y})_{\Gamma_{N_1},I} + (q, \Phi_{u:x})_{\Gamma_{N_2},I}.
\]
The overall setting is shown in \cref{fig:EX_6_experiments}.
Because of the second control boundary, the Tikhonov term in the cost functional
now becomes an integral over the union
$\Gamma_N \coloneqq \Gamma_{N_1} \cup\Gamma_{N_2}$.
The domain is partitioned into $128 \times 128$ square elements
with diameter $h = \sqrt2/128$.
The number of time steps is $M = 100$.
The desired phase-field is given as
\begin{equation*}
  \label{EQ:EX_6_phid}
  \varphi_d(x, y) \coloneqq
  \begin{cases}
    0, & x \in (0.1,0.5) \text{ and } \abs{y - (0.85 - 0.7 x)} \le 3 h, \\
    1, & \text{else}.
  \end{cases}
\end{equation*}
In short, the desired crack goes diagonally from
$(0.5,0.5)$ to $(0.1,0.78)$ with a vertical diameter of $6 h$.
The results are presented in
\cref{tab:results6,fig:EX_6_PF,fig:EX_6_ADJ,fig:EX_6_DISPL,fig:EX_6_Forces}.
From \cref{tab:results6} we can see that it takes 13 iterations
to solve the NLP with an absolute tolerance of \num{2.0e-10}.
Note that from now on the first two columns (Iter, CG)
give the iteration index of Newton's method on the reduced NLP
and the number of CG iterations required
for computing the Newton increment, respectively.
The Newton iteration terminates when either the relative
or the absolute residual falls below the requested tolerance.
The cost functional is reduced from \num{4.47e-2} to \num{1.28e-2},
by approximately 70\%.
The final phase-field is shown in \cref{fig:EX_6_PF}.
As one can clearly see, the desired diagonal crack propagation
has been produced successfully.
On the one hand, the crack has to propagate to the left,
therefore the control on the upper boundary $\Gamma_{N_1}$
has to increase from left to right.
On the other hand, the crack should propagate upwards,
therefore the control on the left boundary $\Gamma_{N_2}$
has to decrease from bottom to top.
In contrast to Experiment~1,
no symmetry in the displacement or adjoints fields can be expected
since here the notch is horizontal whereas
the desired phase-field is diagonal.
Note that \cref{fig:EX_6_Forces} shows a kink in each control.
This is a numerical artefact: at the cell in the top left corner,
the control acts on two adjacent boundaries simultaneously,
and the discretized quantities
interact within this single cell.
In \cref{fig:EX_6_PF} we finally notice a tiny crack propagation
starting from the bottom left edge $(0, 0)$.
This is due to the singularity caused by
the Dirichlet condition on the bottom boundary $\Gamma_D$
in combination with
the control acting as Neumann condition on $\Gamma_{N_2}$.
Similar observations are made in Section \ref{SEC:Experiment4}.

\newcommand\captionperf[2][]{\caption{Experiment #2:
    number of CG iterations, residuals, cost terms and maximal force
    during NLP iteration.\ifx#1\empty\else\ #1.\fi}}
\begin{table}[tp]
  \centering
  \captionperf{2}
  \sisetup{table-format=1.4e+1}
  \begin{tabular}{rr
    S[table-format=1.2e+1]S[table-format=1.2e+2]SSSS[table-format=4.2]}
    \toprule
    Iter&CG&{Relative}&{Absolute}&{Cost}&{Tracking}&{Tikhonov}&{Force}\\[-2pt]
        &  &{residual}&{residual}&      &          &       &          \\
    \midrule
     0 & -- & 1.0     & 1.99e-05 & 4.4742e-2 & 1.3723e-2 & 3.1019e-2 & 10     \\
     1 &  2 & 3.70e-3 & 7.36e-08 & 1.3075e-2 & 1.3075e-2 & 8.9902e-9 & 2202.36\\
     2 &  9 & 1.27e-3 & 2.52e-08 & 1.2879e-2 & 1.2824e-2 & 5.4707e-5 & 2443.45\\
     3 &  6 & 5.77e-4 & 1.15e-08 & 1.2818e-2 & 1.2719e-2 & 9.8469e-5 & 2525.62\\
     4 &  5 & 3.11e-4 & 6.20e-09 & 1.2790e-2 & 1.2667e-2 & 1.2261e-4 & 2565.05\\
     5 &  4 & 1.84e-4 & 3.67e-09 & 1.2775e-2 & 1.2638e-2 & 1.3678e-4 & 2584.23\\
     6 &  4 & 1.13e-4 & 2.26e-09 & 1.2766e-2 & 1.2620e-2 & 1.4553e-4 & 2596.34\\
     7 &  3 & 7.57e-5 & 1.51e-09 & 1.2760e-2 & 1.2609e-2 & 1.5102e-4 & 2605.50\\
     8 &  2 & 4.86e-5 & 9.67e-10 & 1.2756e-2 & 1.2602e-2 & 1.5474e-4 & 2609.38\\
     9 &  3 & 3.53e-5 & 7.03e-10 & 1.2754e-2 & 1.2597e-2 & 1.5710e-4 & 2613.69\\
     10&  2 & 2.42e-5 & 4.82e-10 & 1.2752e-2 & 1.2593e-2 & 1.5886e-4 & 2614.90\\
     11&  2 & 1.55e-5 & 3.09e-10 & 1.2751e-2 & 1.2591e-2 & 1.6008e-4 & 2616.43\\
     12&  2 & 1.01e-5 & 2.01e-10 & 1.2750e-2 & 1.2589e-2 & 1.6086e-4 & 2617.60\\
     13&  2 & 6.81e-6 & 1.36e-10 & 1.2750e-2 & 1.2588e-2 & 1.6138e-4 & 2618.43\\
    \bottomrule
  \end{tabular}
  \label{tab:results6}
\end{table}

\begin{table}
  \captionpar{2}
  \centering
  \begin{tabular}{llS[table-format=1.1e+2]}
    \toprule
    Par. & Definition & {Value} \\
    \midrule
    $\varepsilon$ & Regul. (crack) $\approx 4 h$ & 0.0442 \\
    $\kappa$ & Regul. (crack) & 1.0e-10 \\
    $\eta$ & Regul. (viscosity) & 1.0e3 \\
    $\gamma$ & Penalty & 1.0e5 \\
    $\alpha$ & Tikhonov & 6.5e-9 \\
    \bottomrule
  \end{tabular}
  \hfil
  \begin{tabular}{llS[table-format=2.1e1]}
    \toprule
    Par. & Definition & {Value} \\
    \midrule
    $G_c$ & Fracture toughness & 1.0 \\
    $\nu_s$ & Poisson's ratio & 0.2 \\
    $E$ & Young's modulus & 1.0e6 \\
    $q_0$ & Initial control & 10.0 \\
    $q_d$ & Nominal control & 2.2e3 \\
    \bottomrule
  \end{tabular}
  \label{tab:EX_6_params}
\end{table}%

\renewcommand\img[1]{\image
  [trim=840 315 840 315,clip=true,width=.325\linewidth]
  {EX_6_STATE_13_TIME_#1_PF.png}}
\begin{figure}[tp]
  \footnotesize
  \centering
  \img{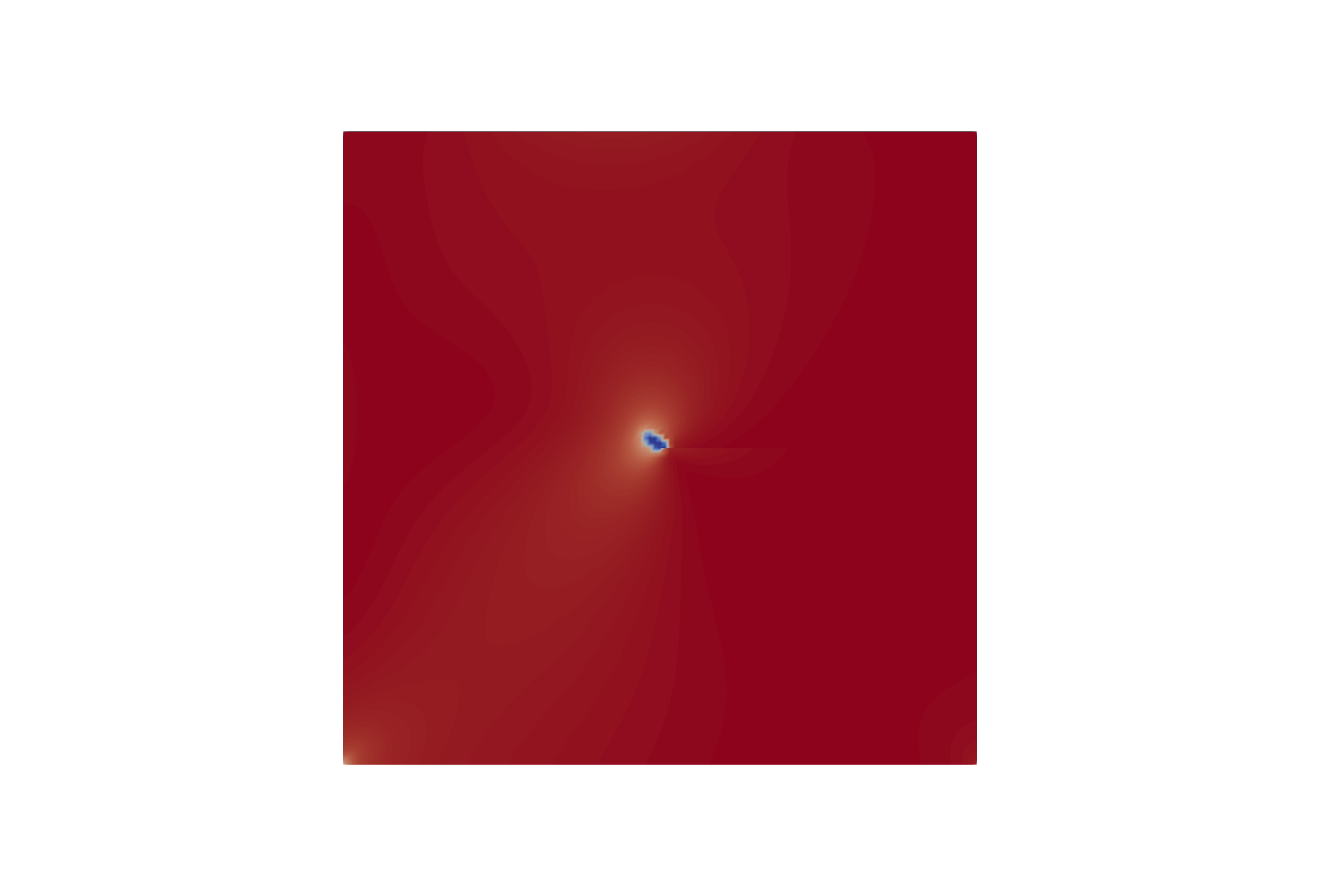}\hfill\img{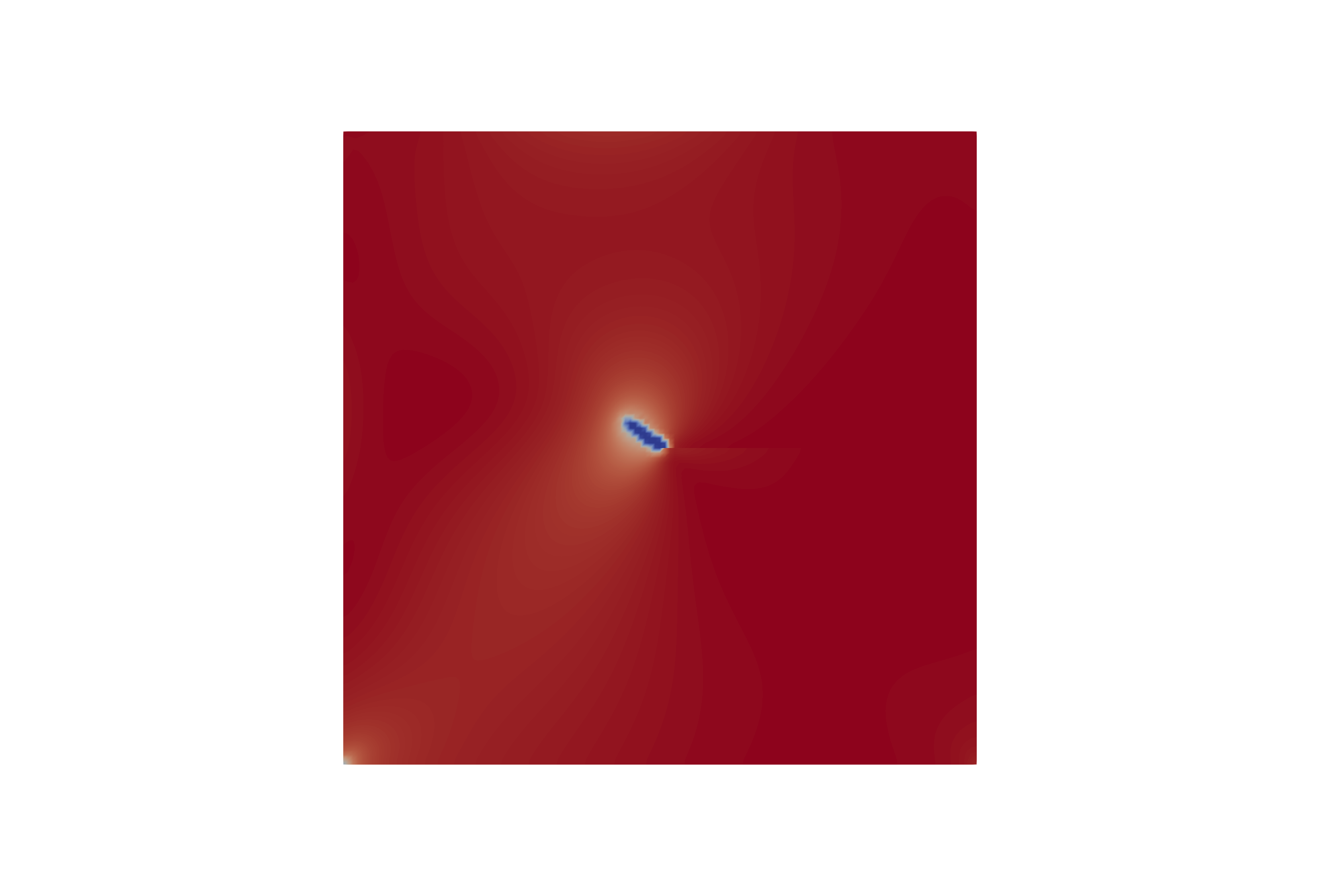}\hfill\img{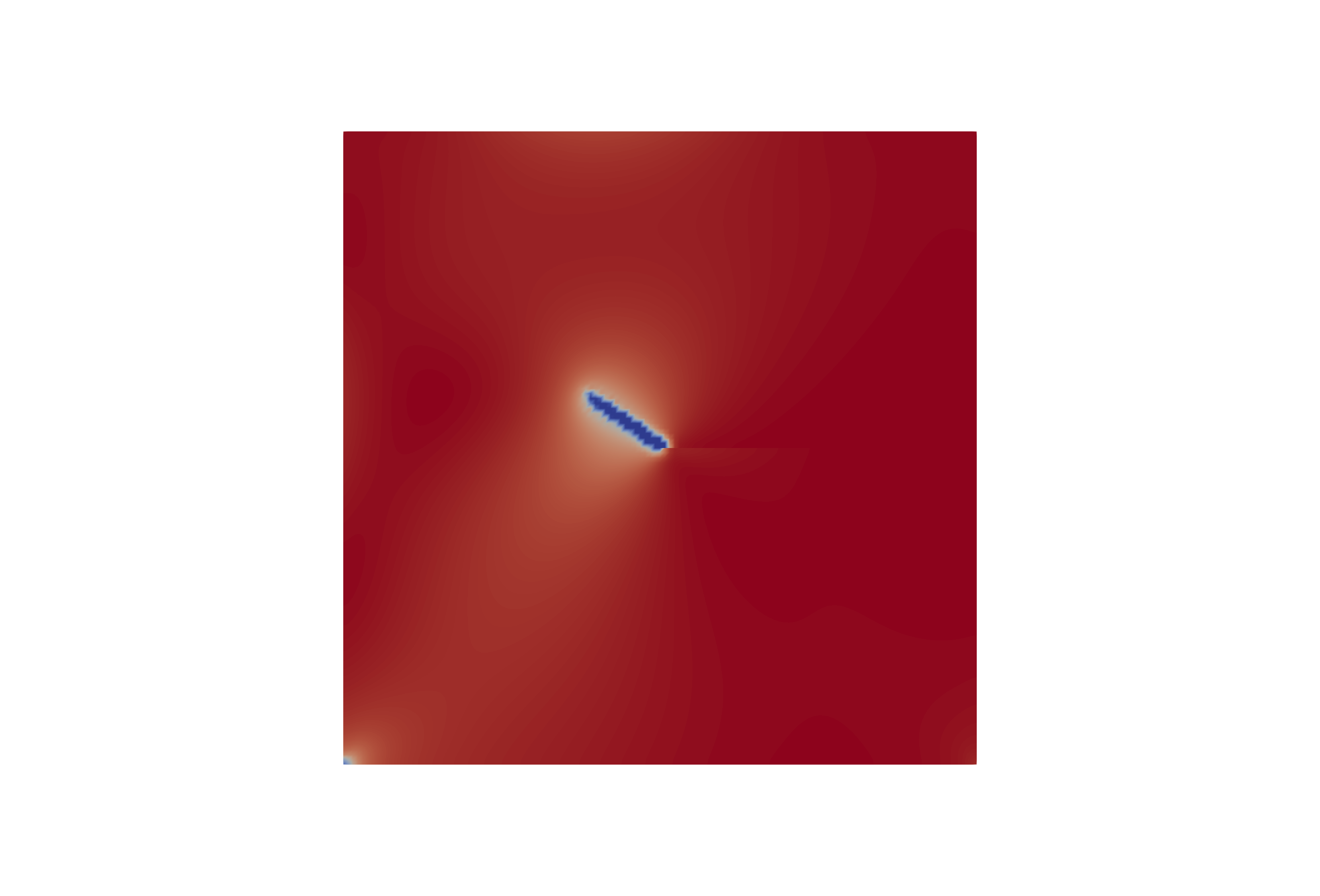}\\
  \scale[0.33]{0.32}{0}{0.5}{1}%
  \caption{Experiment 2:
    optimal phase-field $\varphi$ at times 50, 75, and 100.}
  \label{fig:EX_6_PF}
\end{figure}

\renewcommand\img[2]{\image
  [trim=845 320 850 320,clip=true,width=0.36\linewidth]
  {EX_6_#1_13_TIME_100_#2.png}}
\begin{figure}[tp]
  \footnotesize
  \centering
  \img{STATE}{U_X}\quad\img{STATE}{U_Y}\\
  \scale{0.36}{-2.2e-2}{}{1.5e-2}\quad
  \scale{0.36}{-3.2e-3}{}{5.8e-2}\\[3ex]
  \img{ADJ}{Z_X}\quad\img{ADJ}{Z_Y}\\
  \scale{0.36}{-2.1e-8}{}{3.8e-8}\quad
  \scale{0.36}{-8.5e-8}{}{1.4e-8}%
  \caption{Experiment 2:
    optimal displacement field $u$ (top: $x$ left, $y$ right)
    and adjoint field $z_u$ (bottom: $x$ left, $y$ right) at time 250.}
  \label{fig:EX_6_DISPL}
  \label{fig:EX_6_ADJ}
\end{figure}

\subsection{Experiment 3:
 connecting horizontal cracks for a sliced domain}
\label{SEC:Experiment2}

Our third experiment is motivated by a simple question:
Is it possible to connect some (but not all) notches in a given domain?
Here we consider the rectangule $\Omega = (0, 2.2) \times (0, 0.4)$
with four horizontal notches
$\notch_1 \coloneqq (0.3,0.5) \times \set{0.2}$,
$\notch_2 \coloneqq (0.7,0.9) \times \set{0.2}$,
$\notch_3 \coloneqq (1.3,1.5) \times \set{0.2}$,
$\notch_4 \coloneqq (1.7,1.9) \times \set{0.2}$,
see \cref{fig:EX_2_experiments}.
This yields the combined notch $\notch \coloneqq \bigcup_{i=1}^4 \notch_i$
with initial phase-field
\begin{equation*}
  \label{EQ:EX_2_phi0}
  \varphi_0(x, y) \coloneqq
  \begin{cases}
    0, & (x,y) \in \notch, \\
    1, & \text{else}.
  \end{cases}
\end{equation*}
The boundary $\partial \Omega$ is partitioned as in \cref{SEC:Experiment1}.
The time interval is again $[0, 1]$ but with 2001 equidistant time points,
i.e., $T = 1$ and $M = 2000$.
The spatial mesh now consists of $352 \times 64$ square elements
with diameter $h = \sqrt2 \times 0.4 / 64 \approx 0.00884$.
The desired phase-field $\varphi_d$ connects $\notch_1$ with $\notch_2$
and $\notch_3$ with $\notch_4$, hence it is defined as follows:
\begin{equation*}
  \label{EQ:EX_2_phid}
  \varphi_d(x, y) \coloneqq
  \begin{cases}
    0, & x \in (0.5,0.7) \cup (1.5,1.7) \text{ and }
    y \in (0.2 - 4h, 0.2 + 4h), \\
    1, & \text{else}.
  \end{cases}
\end{equation*}
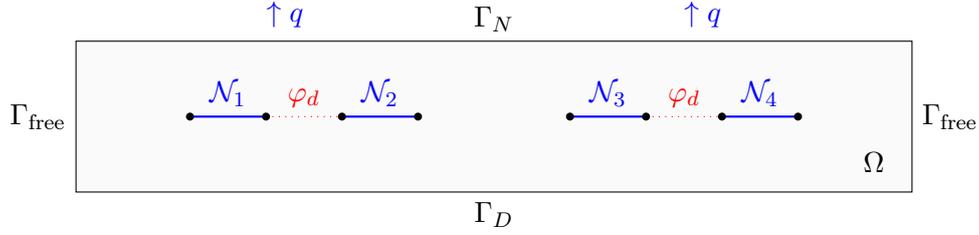
\begin{figure}[tp]
  \centering
  \begin{tikzpicture}[scale=5.0]
    \draw[fill=black!2] (0,0)
    -- (2.2,0) node[pos=0.5,below]  {$\Gamma_D$}
    -- (2.2,0.4) node[pos=0.5,right] {$\Gamma\free$}
    -- (0,0.4) node[pos=0.5,above] {$\Gamma_N$}
    -- cycle node[pos=0.5,left] {$\Gamma\free$};
    \node at (2.1,0.08) {$\Omega$};
    \node[blue,above] at (0.55,0.41) {$\uparrow q$};
    \node[blue,above] at (1.65,0.41) {$\uparrow q$};
    \draw[dotted,red]
    (0.5,0.2) -- (0.7,0.2) node[pos=0.5,above,red] {$\varphi_d$}
    (1.5,0.2) -- (1.7,0.2) node[pos=0.5,above,red] {$\varphi_d$};
    \draw[thick,blue]
    (0.3,0.2) -- (0.5,0.2) node[pos=0.5,above,blue] {$\notch_1$}
    (0.7,0.2) -- (0.9,0.2) node[pos=0.5,above,blue] {$\notch_2$}
    (1.3,0.2) -- (1.5,0.2) node[pos=0.5,above,blue] {$\notch_3$}
    (1.7,0.2) -- (1.9,0.2) node[pos=0.5,above,blue] {$\notch_4$};
    \fill
    (0.3,0.2) circle (0.3pt) (0.5,0.2) circle (0.3pt)
    (0.7,0.2) circle (0.3pt) (0.9,0.2) circle (0.3pt)
    (1.3,0.2) circle (0.3pt) (1.5,0.2) circle (0.3pt)
    (1.7,0.2) circle (0.3pt) (1.9,0.2) circle (0.3pt);
  \end{tikzpicture}
  \caption{Experiment 3:
    domain $\Omega = (0, 2.2) \times (0, 0.4)$
    with partitioned boundary $\partial \Omega$,
    initial notches $\notch_1,\dots,\notch_4$,
    and desired cracks $\varphi_d$.}
  \label{fig:EX_2_experiments}
\end{figure}%
All relevant parameters for this experiment are presented
in \cref{tab:EX_2_params}.
The results for the tolerance \num{2.0e-9} are shown in \cref{tab:results2}.
From the final optimal phase-field in \cref{fig:EX_2_PF} (bottom)
we see that the desired phase-field has indeed been reached,
since $\notch_1$ is connected with $\notch_2$ and $\notch_3$ with $\notch_4$.
The optimal control force shown in \cref{fig:EX_2_Force}
is rather strong and has two roughly parabolic maxima
right at the two sections where notches are to be connected,
which is to be expected from a mechanical point of view.
The four cracks propagating from both ends
of each pair of connected notches,
where no cracks are desired,
can be explained by the decreasing control at the end points.
That decreasing control generates
a different principal axis of tension
which in turn produces the non-horizontal crack growth.
In \cref{fig:EX_2_DISPL2000} (top) we present the optimal displacement fields
at time step $2000$.
They are both symmetric and reach their maxima right at
the two sections where notches are to be connected.
This is consistent with the behavior of the control forces
and again physically plausible.
For comparison, before the middle
cracks join, we also display the respective fields at time step
$1800$ in \cref{fig:EX_2_DISPL1800,fig:EX_2_ADJ1800}.

\begin{table}[tp]
  \captionpar{3}
  \centering
  \begin{tabular}{llS[table-format=1.1e+2]}
    \toprule
    Par. & Definition & {Value} \\
    \midrule
    $\varepsilon$ & Regul. (crack) $\approx 4 h$ & 0.035 \\
    $\kappa$ & Regul. (crack) & 1.0e-10 \\
    $\eta$ & Regul. (viscosity) & 1.0e3 \\
    $\gamma$ & Penalty & 1.0e5  \\
    $\alpha$ & Tikhonov & 2.1e-10 \\
    \bottomrule
  \end{tabular}
  \hfil
  \begin{tabular}{llS[table-format=1.2e1]}
    \toprule
    Par. & Definition & {Value} \\
    \midrule
    $G_c$ & Fracture toughness & 1.0 \\
    $\nu_s$ & Poisson's ratio & 0.2 \\
    $E$ & Young's modulus & 1.00e6 \\
    $q_0$ & Initial control & 1.0 \\
    $q_d$ & Nominal control & 6.53e3 \\
    \bottomrule
  \end{tabular}
  \label{tab:EX_2_params}
\end{table}%
\begin{table}[tp]
  \centering
  \captionperf{3}
  \sisetup{table-format=1.4e+1}
  \begin{tabular}{rr
    S[table-format=1.2e+1]S[table-format=1.2e+2]SSSS[table-format=4.2]}
    \toprule
    Iter&CG&{Relative}&{Absolute}&{Cost}&{Tracking}&{Tikhonov}&{Force}\\[-2pt]
        &  &{residual}&{residual}&      &          &       &          \\
    \midrule
    0 & -- & 1.0     & 2.00e-06 & 2.3850e-2 & 1.4302e-2 & 9.5483e-3  & 100.00  \\
    1 &  2 & 0.110   & 2.20e-07 & 1.2839e-2 & 1.2839e-2 & 3.4223e-8  & 6550.69 \\
    2 &  2 & 3.13e-2 & 6.27e-08 & 1.2428e-2 & 1.2309e-2 & 1.1897e-4  & 7899.32 \\
    3 &  2 & 1.22e-2 & 2.44e-08 & 1.2330e-2 & 1.2139e-2 & 1.9136e-4  & 8264.47 \\
    4 &  2 & 9.34e-3 & 1.87e-08 & 1.2292e-2 & 1.2066e-2 & 2.2519e-4  & 8403.51 \\
    5 &  2 & 1.44e-3 & 2.89e-09 & 1.2263e-2 & 1.2011e-2 & 2.5205e-4  & 8501.55 \\
    6 &  2 & 9.20e-4 & 1.84e-09 & 1.2261e-2 & 1.2006e-2 & 2.5478e-4  & 8514.26 \\
    \bottomrule
  \end{tabular}
  \label{tab:results2}
\end{table}
\renewcommand\img[1]{\image
  [trim=540 915 540 915,clip=true,width=0.75\linewidth]
  {EX_2_STATE_6_TIME_#1_PF.png}}
\begin{figure}[tp]
  \footnotesize
  \centering
  \img{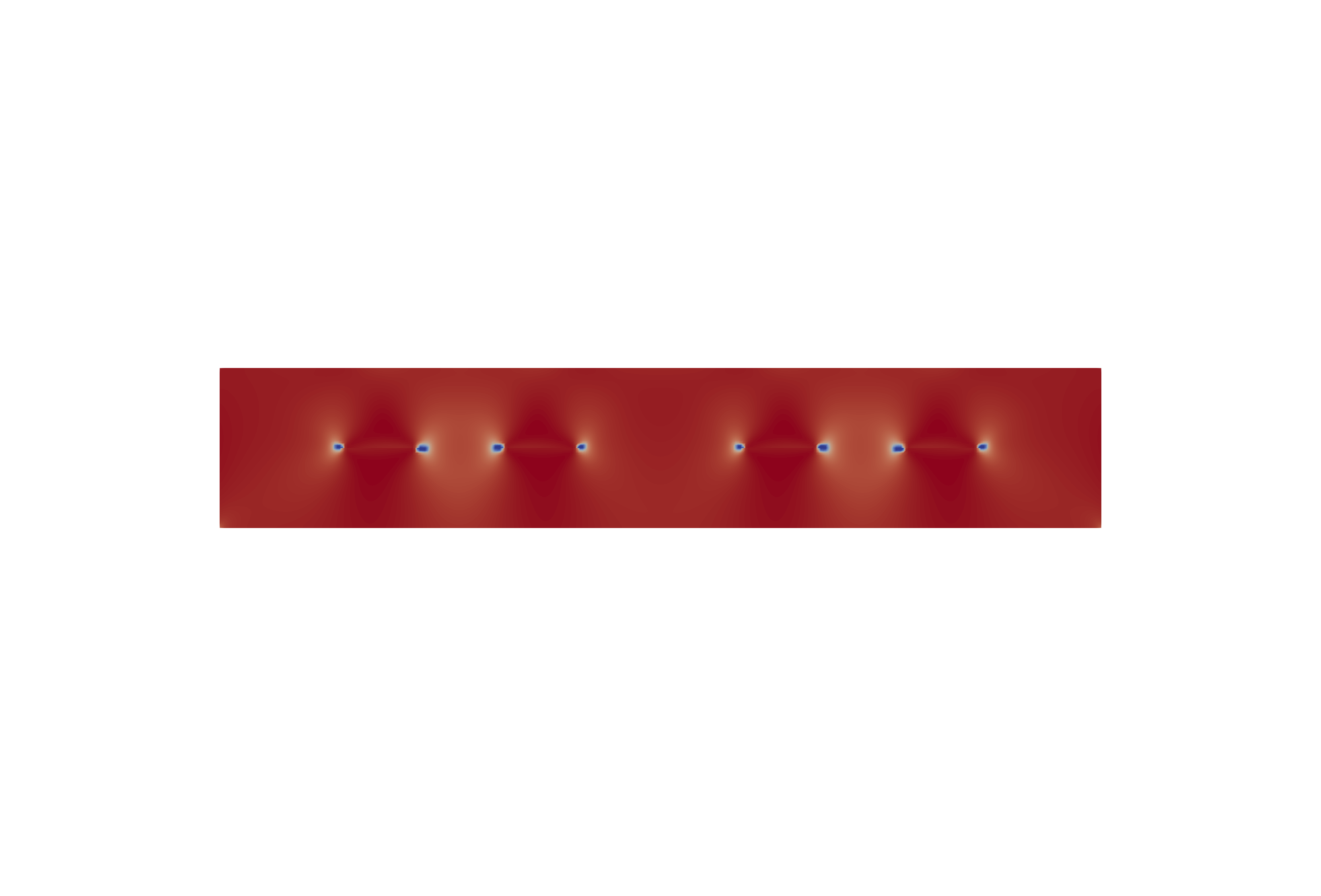}\\[3pt]
  \img{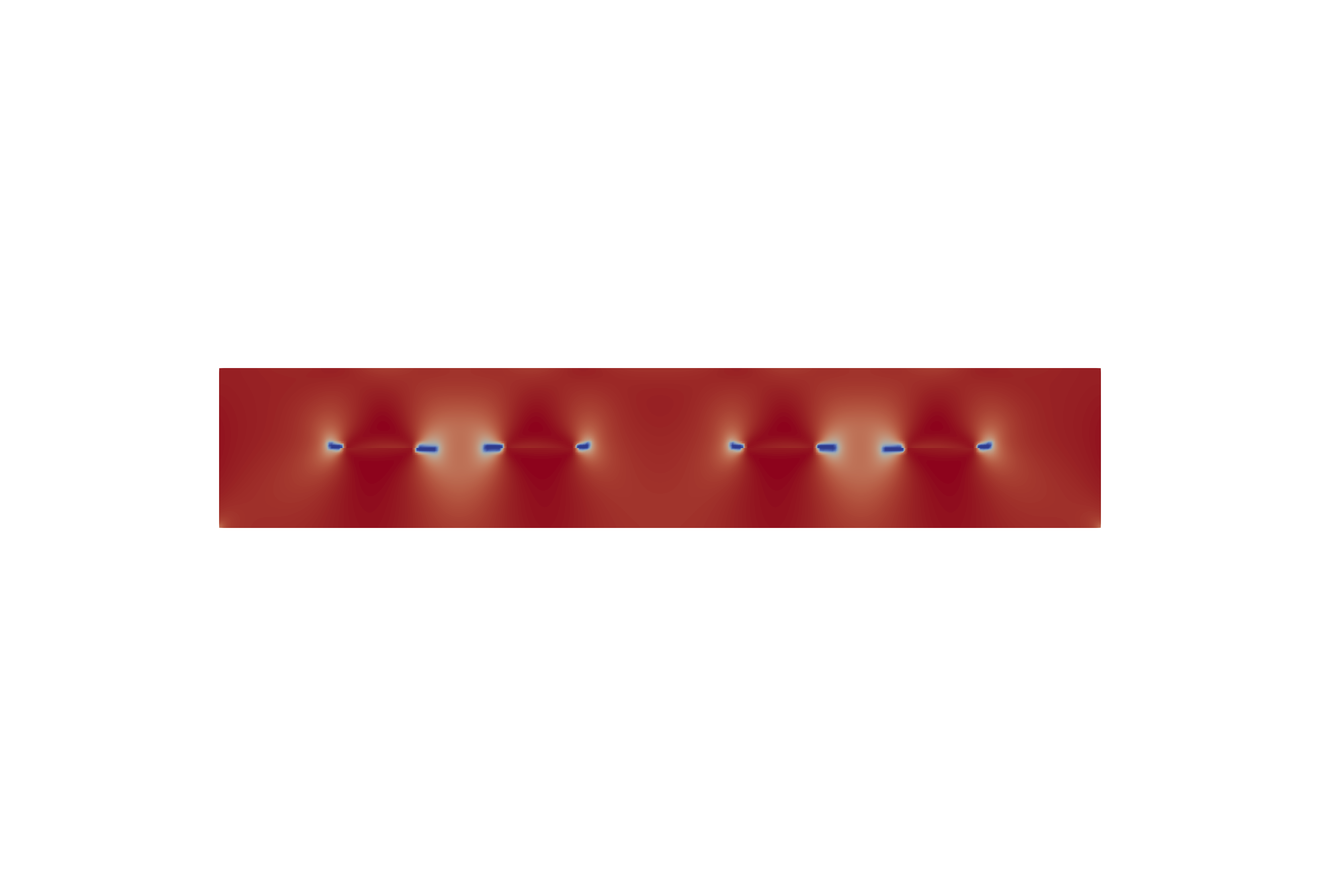}\\[3pt]
  \img{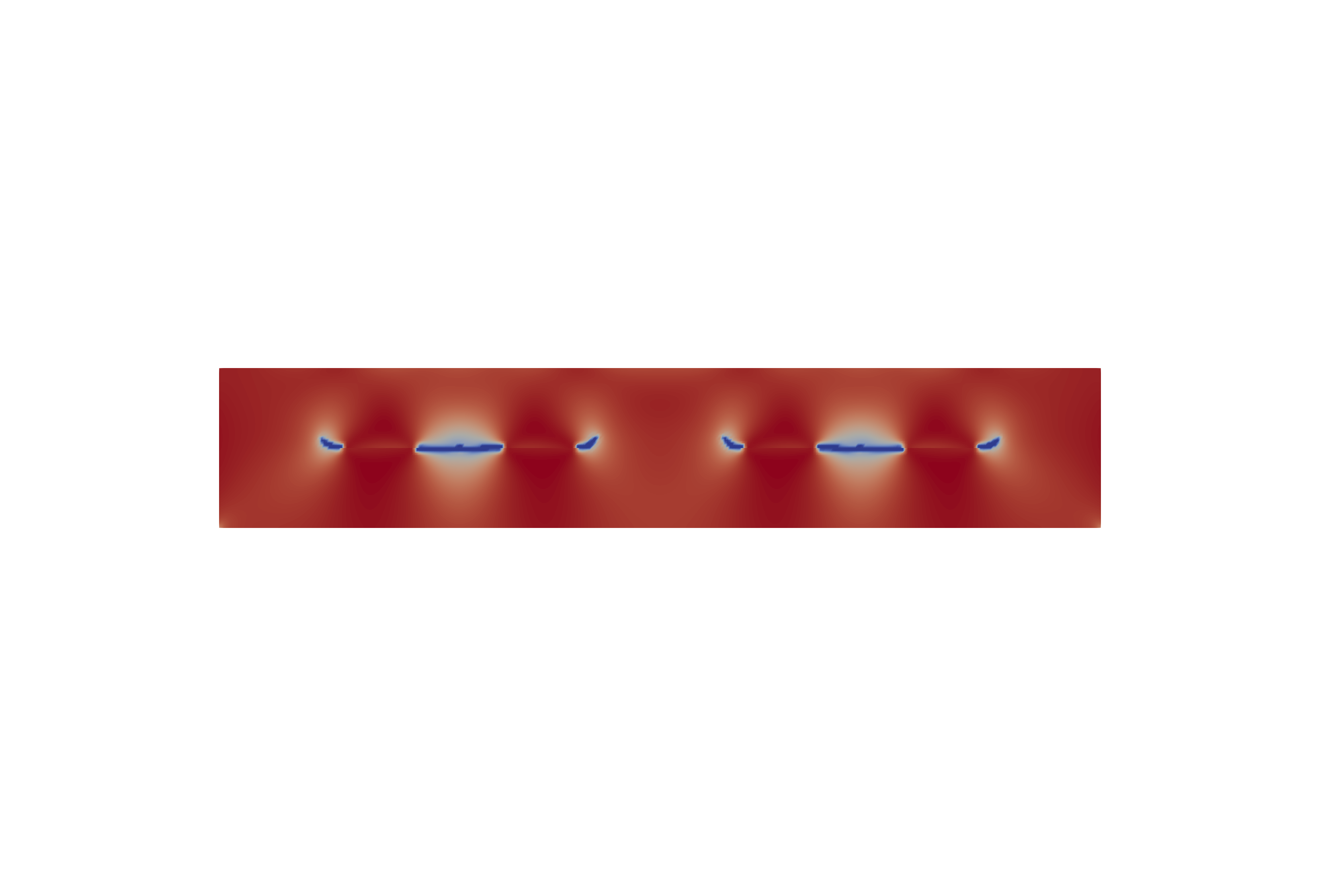}\\[3pt]
  \scale[0.33]{0.32}{0}{0.5}{1}%
  \caption{Experiment 3:
    optimal phase-field $\varphi$ at times 1400, 1800, and 2000.}
  \label{fig:EX_2_PF}
\end{figure}

\renewcommand\img[2]{\image
  [trim=625 910 625 910,clip=true,width=0.48\linewidth]
  {EX_2_#1_6_TIME_1800_#2.png}}
\begin{figure}[tp]
  \footnotesize
  \centering
  \img{STATE}{U_X}\hfill\img{STATE}{U_Y}\\%
  \scale{0.48}{-2.2e-3}{}{2.2e-3}\hfill
  \scale{0.48}{0}{}{9.6e-3}\\[2ex]%
  \img{ADJ}{Z_X}\hfill\img{ADJ}{Z_Y}\\%
  \scale{0.48}{-8.6e-8}{}{8.6e-8}\hfill
  \scale{0.48}{-2.3e-7}{}{1.0e-7}%
  \caption{Experiment 3:
    optimal displacement field $u$ (top: $x$ left, $y$ right)
    and adjoint field $z_u$ (bottom: $x$ left, $y$ right) at time 1800.}
  \label{fig:EX_2_DISPL1800}
  \label{fig:EX_2_ADJ1800}
\end{figure}

\renewcommand\img[2]{\image
  [trim=530 910 530 910,clip=true,width=0.48\linewidth]
  {EX_2_#1_6_TIME_2000_#2.png}}
\begin{figure}[tp]
  \footnotesize
  \centering
  \img{STATE}{U_X}\hfill\img{STATE}{U_Y}\\%
  \scale{0.48}{-1.5e-2}{}{1.5e-2}\hfill
  \scale{0.48}{-1.6e-3}{}{5.5e-2}\\[2ex]%
  \img{ADJ}{Z_X}\hfill\img{ADJ}{Z_Y}\\%
  \scale{0.48}{-2.1e-9}{}{2.1e-9}\hfill
  \scale{0.48}{-4.6e-10}{}{4.9e-9}%
  \caption{Experiment 3:
    optimal displacement field $u$ (top: $x$ left, $y$ right)
    and adjoint field $z_u$ (bottom: $x$ left, $y$ right) at time 2000.}
  \label{fig:EX_2_DISPL2000}
  \label{fig:EX_2_ADJ2000}
\end{figure}

\newcommand\captionctrl[2][{[0, 1]} \times \set{1}]
{\caption{Experiment #2: optimal control force (solid) and
    nominal control force (dotted) on upper boundary $\Gamma_N$.}}
\begin{figure}[tp]
  \footnotesize
  \centering
  \begin{tikzpicture}
    \begin{axis}[width=1.0\linewidth,height=60mm,
      xmin=0,xmax=2.2,ymin=6000,ymax=8700,ytick distance={500},grid=major]
      \draw[red,dotted] (0,6530) -- (2.2,6530) node[pos=0.5,above] {$q_d$};
      \addplot[blue] table
      {TXT2_EX_2_Force.txt} node[pos=0.5,below] {$q$};
    \end{axis}
  \end{tikzpicture}
  \captionctrl[{[0, 2.2]} \times \set{0.4}]{3}
  \label{fig:EX_2_Force}
\end{figure}

\subsection{Experiment 4:
  connecting two horizontal cracks for an entirely sliced domain}
\label{SEC:Experiment3}

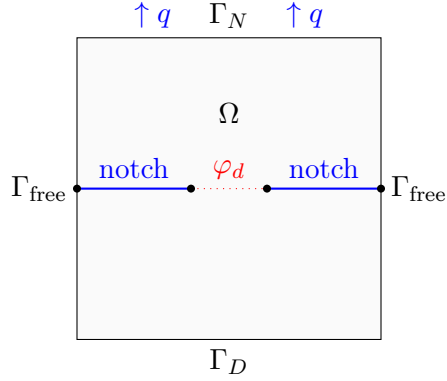
\begin{figure}[tp]
  \centering
  \begin{tikzpicture}
    \draw[fill=black!2] (0,0)
    -- (0,4) node[pos=0.5,left]  {$\Gamma\free$}
    -- (4,4) node[pos=0.5,above] {$\Gamma_N$}
    -- (4,0) node[pos=0.5,right] {$\Gamma\free$}
    -- cycle node[pos=0.5,below] {$\Gamma_D$};
    \node at (2,3) {$\Omega$};
    \node[blue,above] at (1,4) {$\uparrow q$};
    \node[blue,above] at (3,4) {$\uparrow q$};
    \draw[dotted,red] (1.5,2) -- (2.5,2) node[pos=0.5,above,red] {$\varphi_d$};
    \draw[thick,blue] (0,2) -- (1.5,2) node[pos=0.5,above,blue] {notch};
    \draw[thick,blue] (2.5,2) -- (4,2) node[pos=0.5,above,blue] {notch};
    \fill
    (0,2) circle (1.5pt) (1.5,2) circle (1.5pt)
    (2.5,2) circle (1.5pt) (4,2) circle (1.5pt);
  \end{tikzpicture}
  \caption{Experiment 4:
    domain $\Omega = (0, 1)^2$ with partitioned boundary $\partial \Omega$,
    initial notches, and desired crack~$\varphi_d$.}
  \label{fig:EX_3_experiment}
\end{figure}

The fourth experiment is motivated by the question
whether it is possible to connect two horizontal notches
to achieve an entirely sliced domain.
Here we consider again the square domain $\Omega = (0, 1)^2$,
but now with two horizontal notches, see \cref{fig:EX_3_experiment}.
The left notch is defined as $(0.0, 0.375) \times \set{0.5}$,
the right notch is defined as $(0.625,1.0) \times \set{0.5}$.
The boundary $\partial \Omega$ is partitioned as in \cref{SEC:Experiment1}.
We choose the time interval $[0, 1]$
with $251$ equidistant time points,
i.e., $T = 1$ and $M = 250$.
The spatial mesh now consists of $128 \times 128$ square elements
with diameter $h = \sqrt2/128 \approx 0.011$.
The desired phase-field $\varphi_d$ connects
the left notch with the right notch
and is defined as follows:
\begin{equation*}
  \label{EQ:EX_3_phid}
  \varphi_d(x, y) \coloneqq
  \begin{cases}
    0, & x \in (0.375, 0.625) \text{ and }
    y \in (0.5 - 2h, 0.5 + 2h), \\
    1, & \text{else}.
  \end{cases}
\end{equation*}

Our goal in this experiment is rather peculiar. Analytically,
the PDE constraint becomes singular once the domain is entirely sliced.
In the phase-field model this happens
when the left and right boundaries of $\Omega$
are connected by a path along which the phase-field $\varphi$ vanishes.
Numerical difficulties are to be expected even before such a path exists:
the PDE becomes increasingly ill-conditioned
when the transition zones with $0 < \varphi < 1$ come into contact.
Nevertheless it is possible to create a domain-splitting crack
with a pure forward model,
see for instance the related single edge notched tension test
\cite{MieWelHof10a,MieWelHof10b,AmGeraLoren15,BrWiBeNoRa20}.
Yet this experiment remains numerically difficult
and becomes even more challenging within our optimization setting.
Since we have to expect that the solution of the forward problem
might be close to singularities,
it is not clear what will happen when we insert this solution
into the optimization algorithm.
With regard to this challenge we have observed
that in many experiments the Tikhonov term acts against extreme forces
and improves the solvability of the PDE for the resulting controls.
In the experiment under consideration we set $\alpha$ to $\num{2.0e-10}$.
By this the Tikhonov term is not the driving factor of the optimization process,
but still large enough to avoid extreme forces.
The choice of the other parameters is shown in \cref{tab:params_Ex_3},
and our results are presented in
\cref{tab:results3,fig:EX_3_PF,fig:EX_3_DISPL,fig:EX_3_Force}.

In \cref{tab:results3} we observe that the residual value
is decreasing, except for the last iteration.
After iteration~6 the PDE forward problem becomes unsolvable.
Therefore we regard iteration~5 as the optimal solution:
it has the lowest absolute residual value, \num{5.13e-9},
and also the lowest relative residual value, \num{1.57e-2}.
The results presented in
\cref{tab:results3,fig:EX_3_PF,fig:EX_3_DISPL,fig:EX_3_Force}
refer to iteration~5.
The optimal phase-field presented in \cref{fig:EX_3_PF}
does not connect the two notches but reaches approximately two thirds
of the length of the desired phase-field.
In \cref{fig:EX_3_Force} we see that the optimal control force
is nearly twice as large as the nominal control force $q_d$.
This means that the optimization is primarily driven by
the physical term $\norm{\varphi - \varphi_d}^2$.

\begin{table}
  \captionpar{4}
  \centering
  \begin{tabular}{llS[table-format=1.1e+2]}
    \toprule
    Par. & Definition & {Value} \\
    \midrule
    $\varepsilon$ & Regul. (crack) $\approx 2 h$ & 0.0221 \\
    $\kappa$ & Regul. (crack) & 1.0e-10 \\
    $\eta$ & Regul. (viscosity) & 1.0e3 \\
    $\gamma$ & Penalty & 1.0e5  \\
    $\alpha$ & Tikhonov & 2.0e-10 \\
    \bottomrule
  \end{tabular}
  \hfil
  \begin{tabular}{llS[table-format=1.2e1]}
    \toprule
    Par. & Definition & {Value} \\
    \midrule
    $G_c$ & Fracture toughness & 1.0 \\
    $\nu_s$ & Poisson's ratio & 0.2 \\
    $E$ & Young's modulus & 1.00e6 \\
    $q_0$ & Initial control & 1.0 \\
    $q_d$ & Nominal control & 1.85e3 \\
    \bottomrule
  \end{tabular}
  \label{tab:params_Ex_3}
\end{table}%
\begin{table}[tp]
  \centering
  \captionperf{4}
  \sisetup{table-format=1.4e+1}
  \begin{tabular}{rr
    S[table-format=1.2e+1]S[table-format=1.2e+2]SSSS[table-format=4.2]}
    \toprule
    Iter&CG&{Relative}&{Absolute}&{Cost}&{Tracking}&{Tikhonov}&{Force}\\[-2pt]
        &  &{residual}&{residual}&      &          &       &          \\
    \midrule
    0 & -- & 1.0     & 3.26e-07 & 5.5093e-3 & 5.2936e-3 & 2.1566e-4  & 380.0  \\
    1 &  2 & 0.424   & 1.38e-07 & 4.9612e-3 & 4.9584e-3 & 2.7903e-6  & 2043.40\\
    2 &  3 & 0.238   & 7.75e-08 & 4.7486e-3 & 4.6748e-3 & 7.3883e-5  & 2845.85\\
    3 &  3 & 0.145   & 4.72e-08 & 4.6179e-3 & 4.4620e-3 & 1.5592e-4  & 3300.68\\
    4 &  3 & 4.54e-2 & 1.48e-08 & 4.5391e-3 & 4.3184e-3 & 2.2067e-4  & 3573.22\\
    5 &  2 & 1.57e-2 & 5.13e-09 & 4.5160e-3 & 4.2733e-3 & 2.4274e-4  & 3629.90\\
    6 &  2 & 0.26    & 8.36e-08 & 4.5088e-3 & 4.2591e-3 & 2.4976e-4  & 3632.82\\
    \bottomrule
  \end{tabular}
  \label{tab:results3}
\end{table}
\renewcommand\img[1]{\image
  [trim=1005 320 1005 320,clip=true,width=.325\linewidth]
  {EX_3_STATE_5_TIME_#1_PF.png}}
\begin{figure}[tp]
  \footnotesize
  \centering
  \img{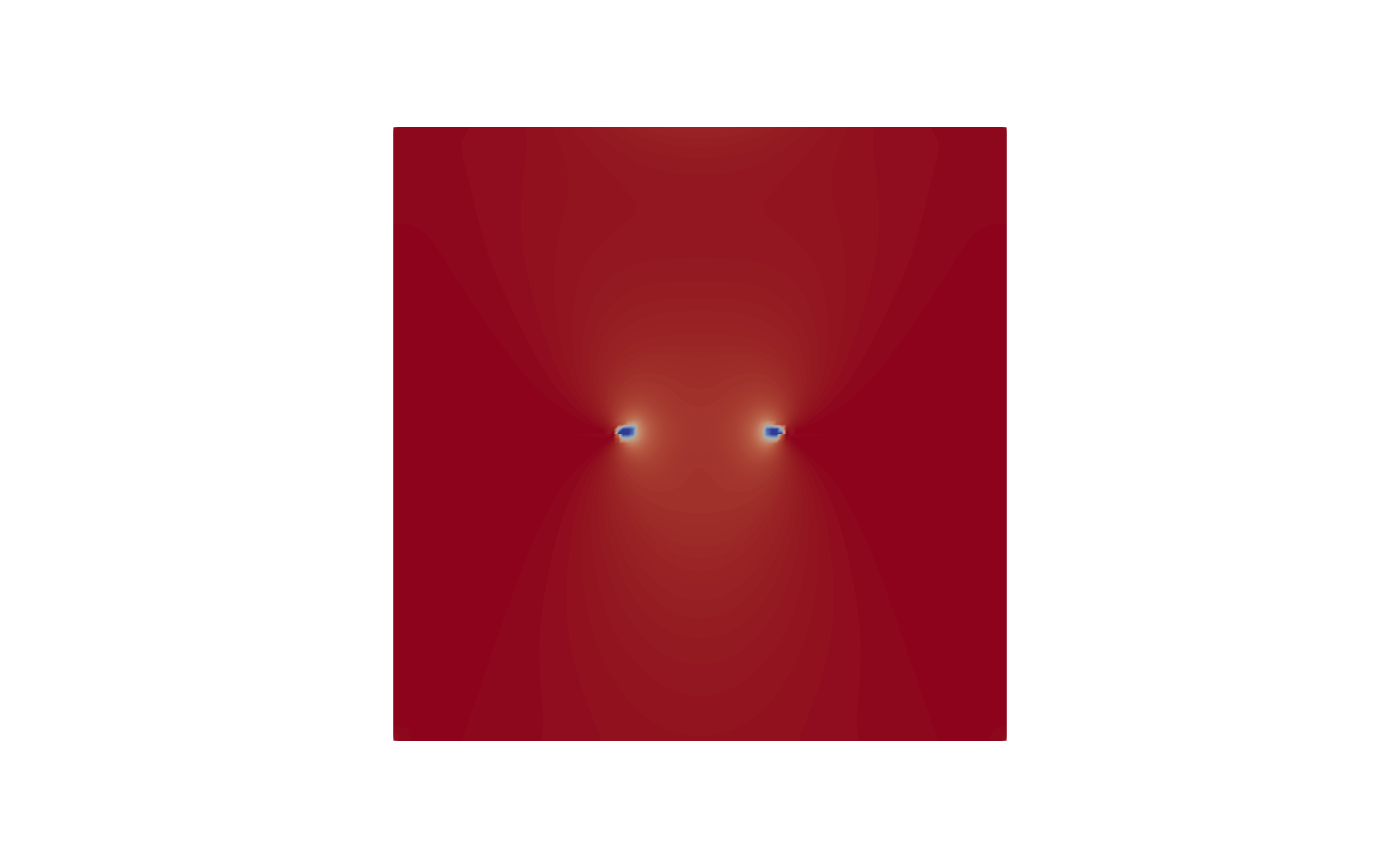}\hfill
  \img{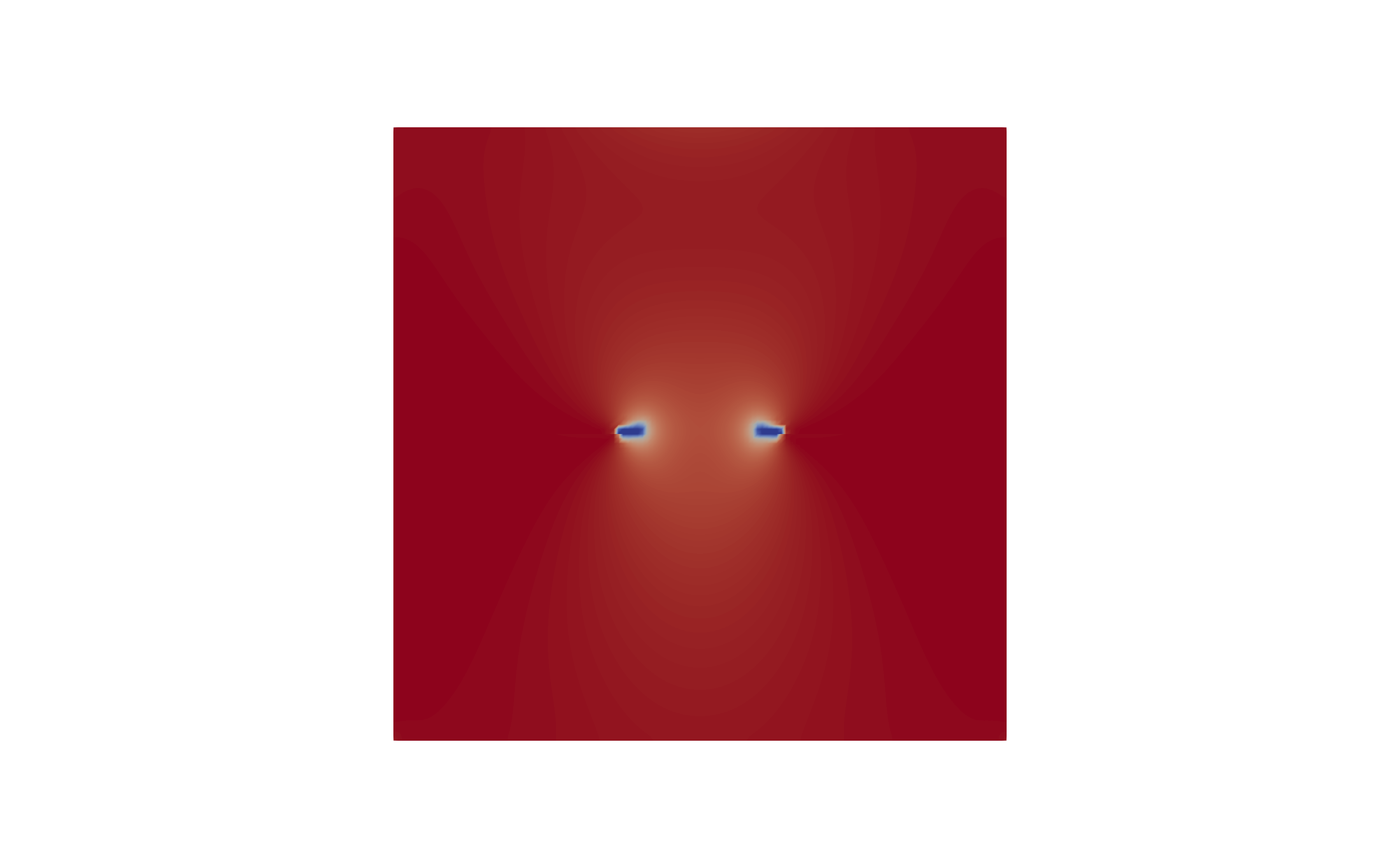}\hfill
  \img{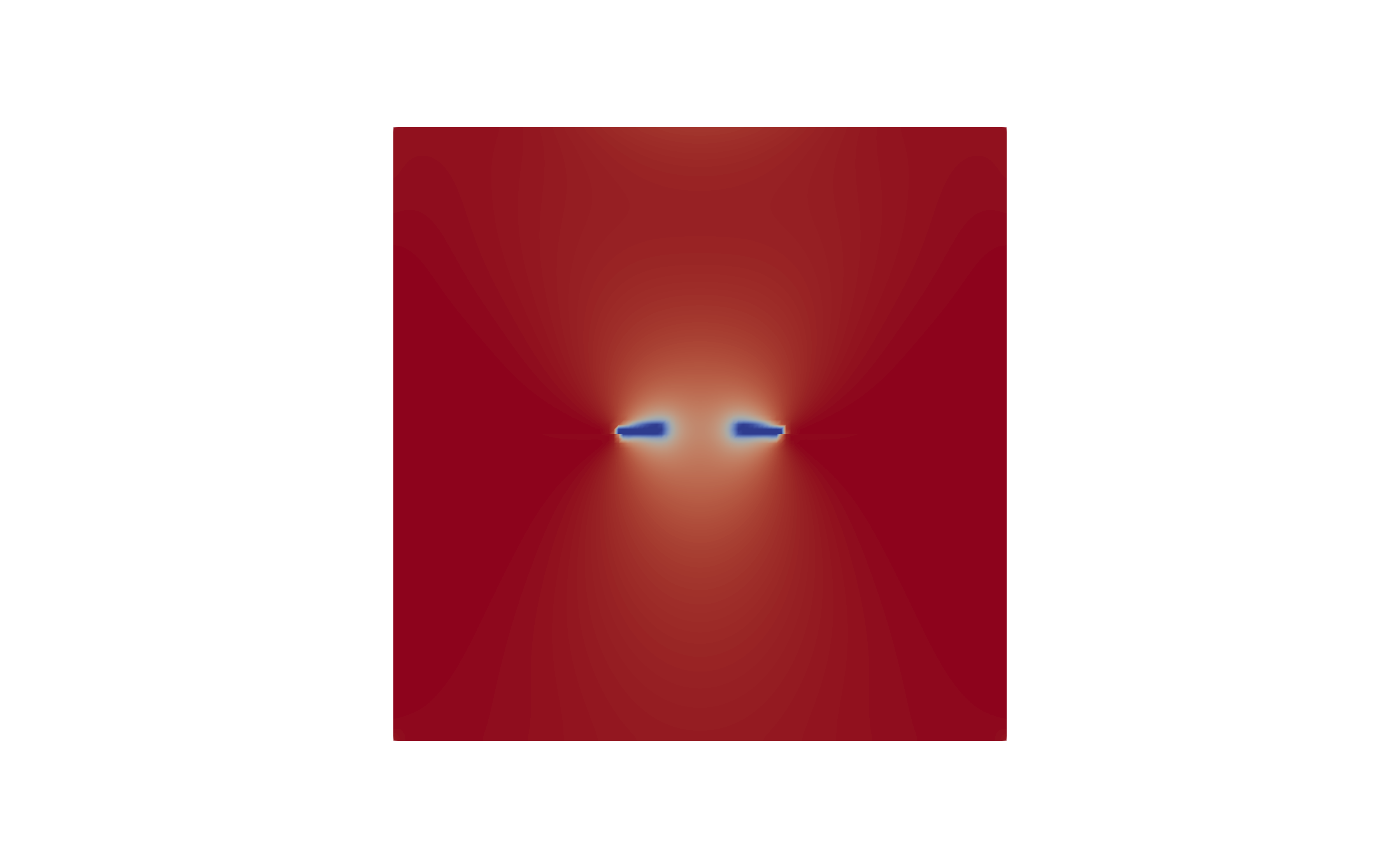}\\%
  \scale[0.33]{0.32}{0}{0.5}{1}%
  \caption{Experiment 4:
    optimal phase-field $\varphi$ at times 150, 200, and 250.}
  \label{fig:EX_3_PF}
\end{figure}
\renewcommand\img[2]{\image
  [trim=1005 320 1005 320,clip=true,width=0.37\linewidth]
  {EX_3_#1_5_TIME_250_#2.png}}
\begin{figure}[tp]
  \footnotesize
  \centering
  \img{STATE}{U_X}\quad\img{STATE}{U_Y}\\%
  \scale{0.37}{-2.5e-3}{}{2.5e-3}\quad
  \scale{0.37}{-1.0e-4}{}{2.1e-2}\\[3ex]%
  \img{ADJ}{Z_X}\quad\img{ADJ}{Z_Y}\\%
  \scale{0.37}{-7.1e-9}{}{7.1e-9}\quad
  \scale{0.37}{-1.3e-8}{}{1.8e-8}%
  \caption{Experiment 4:
    optimal displacement field $u$ (top: $x$ left, $y$ right)
    and adjoint field $z_u$ (bottom: $x$ left, $y$ right) at time 250.}
  \label{fig:EX_3_DISPL}
  \label{fig:EX_3_ADJ}
\end{figure}
\begin{figure}
  \footnotesize
  \centering
  \begin{tikzpicture}
    \begin{axis}[width=1.0\linewidth,height=60mm,
      xmin=0,xmax=1,ymin=1500,ymax=4000,ytick distance={500},grid=major]
      \draw[red,dotted] (0,1850) -- (1,1850) node[pos=0.52,below] {$q_d$};
      \addplot[blue]
      table {TXT3_EX_3_Force.txt} node[pos=0.501,above] {$q$};
    \end{axis}
  \end{tikzpicture}
  \captionctrl{4}
  \label{fig:EX_3_Force}
\end{figure}

\subsection{Experiment 5: L-shaped domain}
\label{SEC:Experiment4}

In our fifth experiment we study a modification of
the L-shaped panel test within an optimization context.
The L-shaped panel test was originally
developed by Winkler \cite{winkler2001}
and extensively studied in
\cite{AmGeraLoren15,MesBouKhon15,Wi17_SISC,Mang_2020}.
In the original test the applied force pushes upwards
against a small left-most section of the upper part of the domain.
In our experiment we apply a pulling force
on the top boundary $\Gamma_N$ instead.
We do this in order to have a complete control boundary
within the optimization context.
The L-shaped domain $\Omega = (0,1)^2 \setminus (0.5,0.5)^2$
and its partitioning of the boundary $\partial \Omega$
are shown in \cref{fig:EX_4_experiments}.
We choose the time interval $[0, 1]$
with 301 equidistant time points, i.e., $T = 1$ and $M = 300$.
Each of the $3 \times 80 \times 80$ square spatial mesh elements
has a diameter of $h = \sqrt2/160 \approx 0.00884$.
All other parameters are shown in \cref{tab:EX_4_params}.
From \cite{AmGeraLoren15,MesBouKhon15,Wi17_SISC,Mang_2020}
and \cite{winkler2001} we already know that the crack will grow
slightly above the horizontal line $[0.5,1] \times \set{0.5}$.
Therefore we place the desired phase-field $\varphi_d$
also slightly above that line,
\begin{equation*}
  \label{EQ:EX_4_phid}
  \varphi_d(x, y) \coloneqq
  \begin{cases}
    0, & x \in (0.5,1.0)\text{ and }
    y \in (0.53 - 4h, 0.53 + 4h), \\
    1, & \text{else}.
  \end{cases}
\end{equation*}
We are aware that a fracture with this phase-field
cannot be produced in our setting for two reasons.
First, a sharp crack along $[0.5,1.0] \times \set{0.53}$
is physically impossible
because the crack will always start to grow
from the singularity in $(0.5,0.5)$.
Second, a decomposition of the stress tensor is needed
in order to distinguish crack growth under tension and compression;
see extensive findings and discussions
for the L-shaped panel test in \cite{AmGeraLoren15}.
Since stress splitting laws introduce further nonlinearities
in the forward phase-field fracture model
and do not contribute to significant further insight in the current work,
we have not used them, despite implemented in our software,
e.g., \cite{Mang_2020}.
We also tried to define $\varphi_d$
on the horizontal line $[0.5,1] \times \set{0.5}$.
However, since the crack starts propagating
diagonally upwards from $(0.5,0.5)$,
the values of the residual and the cost functional
did not decrease, and as a consequence
the Newton iteration for the optimization
problem did not converge.
Our results for the tolerance \num{2.0e-10} are presented in
\cref{tab:results4,fig:EX_4_PF,fig:EX_4_DISPL,fig:EX_4_ADJ,fig:EX_4_Force}.

\begin{figure}[tp]
  \centering
  \begin{tikzpicture}
    \draw[fill=black!2] (0,2)
    -- (0,4) node[pos=0.5,left]  {$\Gamma\free$}
    -- (4,4) node[pos=0.5,above] {$\Gamma_N$}
    -- (4,0) node[pos=0.5,right] {$\Gamma\free$}
    -- (2,0) node[pos=0.5,below] {$\Gamma_D$}
    -- (2,2) node[below left] {$\Gamma\free$}
    -- (0,2);
    \node at (2,3) {$\Omega$};
    \node[blue,above] at (1,4) {$\uparrow q$};
    \node[blue,above] at (3,4) {$\uparrow q$};
    \draw[red,dotted] (2,2.12) -- (4,2.12) node[pos=0.5,above,red] {$\varphi_d$};
    \fill
    (4,2.12) circle (1.5pt) (2.0,2.12) circle (1.5pt);
  \end{tikzpicture}
  \caption{Experiment 5:
    L-shaped domain $\Omega = (0,1)^2 \setminus (0.5,0.5)^2$
    with partitioned boundary and desired crack $\varphi_d$.}
  \label{fig:EX_4_experiments}
\end{figure}%
\begin{table}[tp]
  \captionpar{5}
  \centering
  \begin{tabular}{llS[table-format=1.3e+2]}
    \toprule
    Par. & Definition & {Value} \\
    \midrule
    $\varepsilon$ & Regul. (crack) $\approx 4 h$ & 0.0354 \\
    $\kappa$ & Regul. (crack) & 1.000e-10 \\
    $\eta$ & Regul. (viscosity) & 1.000e3 \\
    $\gamma$ & Penalty & 1.000e5  \\
    $\alpha$ & Tikhonov & 2.625e-9 \\
    \bottomrule
  \end{tabular}
  \hfil
  \begin{tabular}{llS[table-format=1.1e1]}
    \toprule
    Par. & Definition & {Value} \\
    \midrule
    $G_c$ & Fracture toughness & 1.0 \\
    $\nu_s$ & Poisson's ratio & 0.2 \\
    $E$ & Young's modulus & 1.0e6 \\
    $q_0$ & Initial control & 1.0 \\
    $q_d$ & Nominal control & 1.6e3 \\
    \bottomrule
  \end{tabular}
  \label{tab:EX_4_params}
\end{table}%
\begin{table}[tp]
  \centering
  \captionperf{5}
  \sisetup{table-format=1.4e+1}
  \begin{tabular}{rr
    S[table-format=1.2e+1]S[table-format=1.2e+2]SSSS[table-format=4.2]}
    \toprule
    Iter&CG&{Relative}&{Absolute}&{Cost}&{Tracking}&{Tikhonov}&{Force}\\[-2pt]
        &  &{residual}&{residual}&      &          &       &          \\
    \midrule
    0 & -- & 1.0     & 4.34e-06 & 2.0637e-2 & 8.9261e-3 & 1.1711e-2 & 1.0\\
    1 &  2 & 0.117   & 5.09e-07 & 1.6652e-2 & 8.8773e-3 & 7.7750e-3 & 1600.27\\
    2 &  2 & 4.71e-2 & 2.04e-07 & 1.6473e-2 & 8.9092e-3 & 7.5637e-3 & 1982.44\\
    3 &  3 & 2.45e-2 & 1.06e-07 & 1.6404e-2 & 8.9210e-3 & 7.4829e-3 & 2143.55\\
    4 &  2 & 1.18e-3 & 5.11e-09 & 1.6367e-2 & 8.9284e-3 & 7.4387e-3 & 2226.13\\
    5 &  3 & 7.86e-4 & 3.41e-09 & 1.6344e-2 & 8.9340e-3 & 7.4103e-3 & 2279.40\\
    6 &  4 & 5.85e-4 & 2.54e-09 & 1.6329e-2 & 8.9382e-3 & 7.3905e-3 & 2304.63\\
    7 &  3 & 4.33e-4 & 1.88e-09 & 1.6317e-2 & 8.9419e-3 & 7.3750e-3 & 2334.29\\
    8 &  3 & 3.37e-4 & 1.46e-09 & 1.6308e-2 & 8.9450e-3 & 7.3630e-3 & 2350.53\\
    9 &  2 & 2.68e-4 & 1.16e-09 & 1.6301e-2 & 8.9478e-3 & 7.3534e-3 & 2360.71\\
    10&  3 & 2.20e-4 & 9.56e-10 & 1.6296e-2 & 8.9501e-3 & 7.3455e-3 & 2374.60\\
    11&  2 & 1.84e-4 & 7.96e-10 & 1.6291e-2 & 8.9521e-3 & 7.3388e-3 & 2380.56\\
    12&  2 & 1.58e-4 & 6.83e-10 & 1.6287e-2 & 8.9539e-3 & 7.3331e-3 & 2387.88\\
    13&  2 & 1.36e-4 & 5.92e-10 & 1.6284e-2 & 8.9555e-3 & 7.3281e-3 & 2394.83\\
    14&  3 & 1.20e-4 & 5.21e-10 & 1.6281e-2 & 8.9569e-3 & 7.3239e-3 & 2401.26\\
    15&  2 & 1.03e-4 & 4.46e-10 & 1.6278e-2 & 8.9583e-3 & 7.3200e-3 & 2403.82\\
    16&  2 & 9.21e-5 & 4.00e-10 & 1.6276e-2 & 8.9594e-3 & 7.3166e-3 & 2407.58\\
    17&  2 & 8.23e-5 & 3.57e-10 & 1.6274e-2 & 8.9605e-3 & 7.3135e-3 & 2411.42\\
    18&  2 & 7.45e-5 & 3.23e-10 & 1.6272e-2 & 8.9615e-3 & 7.3107e-3 & 2415.00\\
    19&  2 & 6.83e-5 & 2.96e-10 & 1.6271e-2 & 8.9624e-3 & 7.3082e-3 & 2418.35\\
    20&  2 & 6.17e-5 & 2.68e-10 & 1.6269e-2 & 8.9633e-3 & 7.3058e-3 & 2421.39\\
    21&  2 & 5.84e-5 & 2.53e-10 & 1.6268e-2 & 8.9640e-3 & 7.3037e-3 & 2424.17\\
    22&  2 & 5.31e-5 & 2.30e-10 & 1.6266e-2 & 8.9648e-3 & 7.3016e-3 & 2426.76\\
    23&  2 & 4.82e-5 & 2.09e-10 & 1.6265e-2 & 8.9658e-3 & 7.2997e-3 & 2429.09\\
    24&  2 & 4.35e-5 & 1.89e-10 & 1.6264e-2 & 8.9661e-3 & 7.2980e-3 & 2431.26\\
    \bottomrule
  \end{tabular}
  \label{tab:results4}
\end{table}
\renewcommand\img[1]{\image
  [trim=1005 320 1005 320,clip=true,width=.325\linewidth]
  {EX_4_STATE_24_TIME_#1_PF.png}}
\begin{figure}[tp]
  \footnotesize
  \centering
  \img{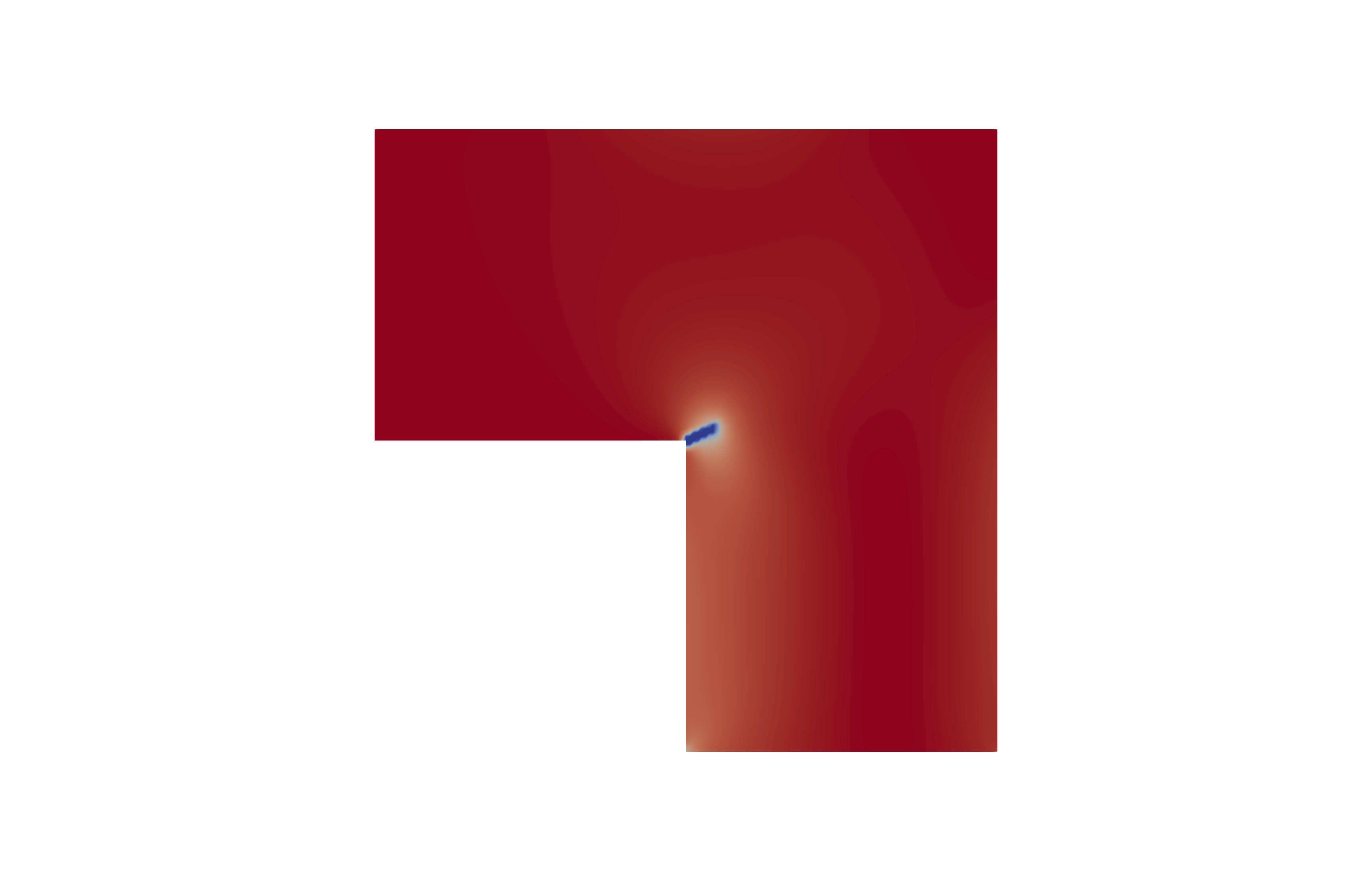}\hfill
  \img{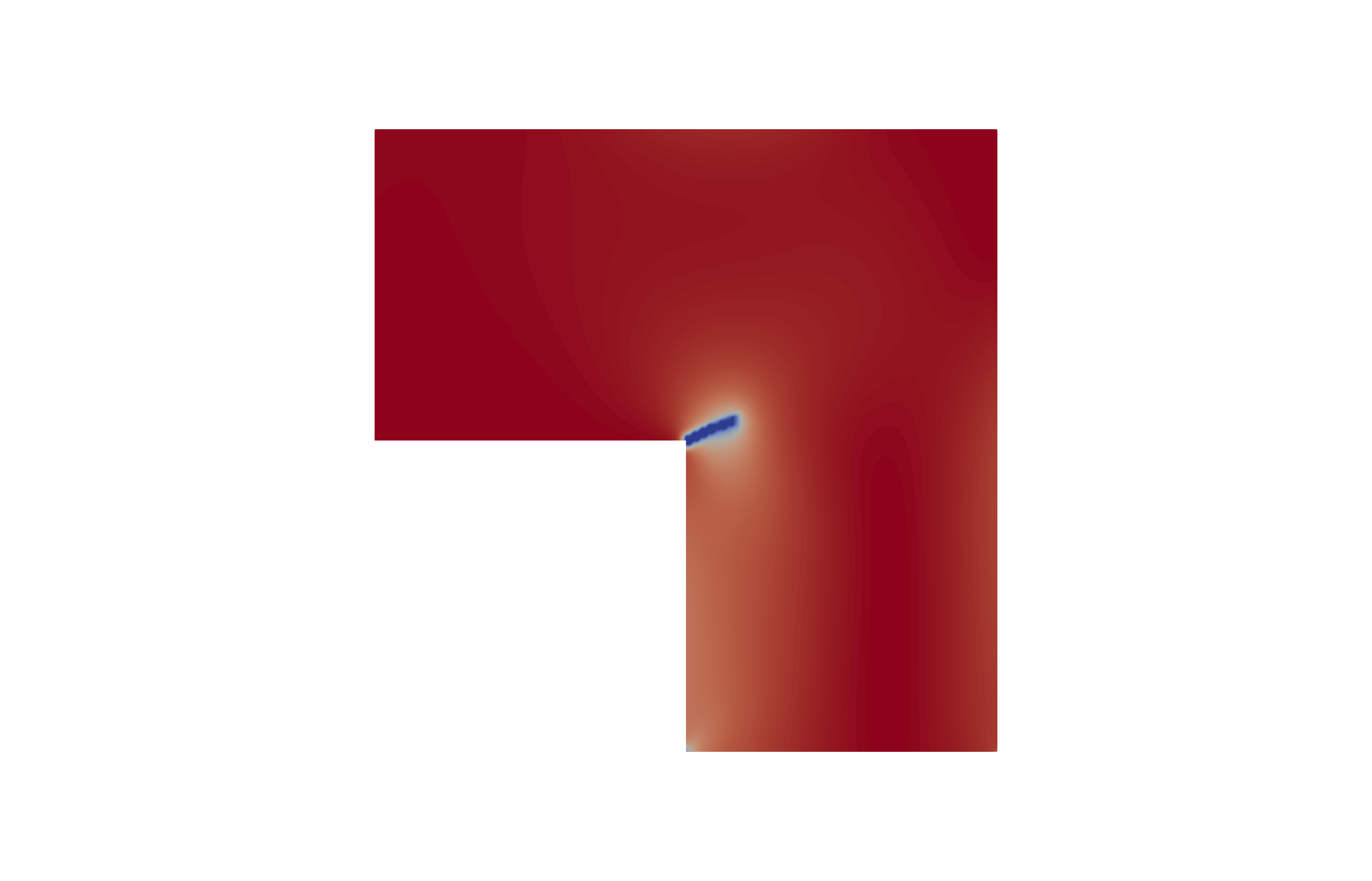}\hfill
  \img{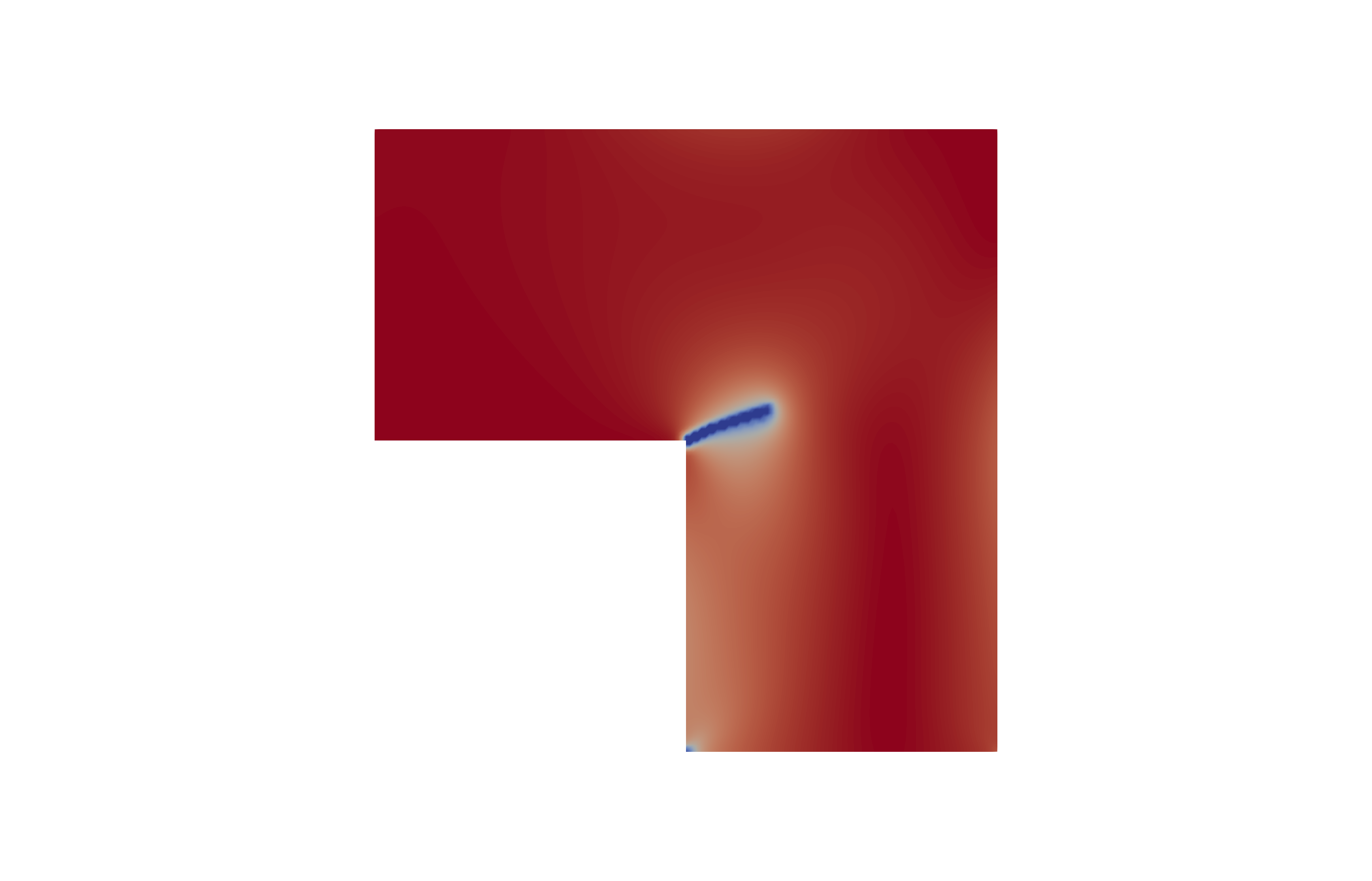}\\%
  \scale[0.33]{0.32}{0}{0.5}{1}%
  \caption{Experiment 5:
    optimal phase-field $\varphi$ at times 200, 250, and 300.}
  \label{fig:EX_4_PF}
\end{figure}
\renewcommand\img[2]{\image
  [trim=1005 320 1005 320,clip=true,width=0.37\linewidth]
  {EX_4_#1_24_TIME_300_#2.png}}
\begin{figure}[tp]
  \footnotesize
  \centering
  \img{STATE}{U_X}\quad\img{STATE}{U_Y}\\%
  \scale{0.37}{-3.8e-5}{}{6.3e-2}\quad
  \scale{0.37}{-1.5e-2}{}{1.0e-1}\\[3ex]%
  \img{ADJ}{Z_X}\quad\img{ADJ}{Z_Y}\\%
  \scale{0.37}{-4.1e-9}{}{1.1e-8}\quad
  \scale{0.37}{-2.1e-8}{}{7.9e-9}%
  \caption{Experiment 5:
    optimal displacement field $u$ (top: $x$ left, $y$ right)
    and adjoint field $z_u$ (bottom: $x$ left, $y$ right) at time 300.}
  \label{fig:EX_4_DISPL}
  \label{fig:EX_4_ADJ}
\end{figure}

\begin{figure}
  \footnotesize
  \centering
  \begin{tikzpicture}
    \begin{axis}[width=1.0\linewidth,height=60mm,
      xmin=0,xmax=1,ymin=1400,ymax=2500,ytick distance={200},grid=major]
      \draw[red,dotted] (0,1600) -- (1,1600) node[pos=0.52,above] {$q_d$};
      \addplot[blue] table
      {TXT4_EX_4_Force.txt} node[pos=0.65,above] {$q$};
    \end{axis}
  \end{tikzpicture}
  \captionctrl{5}
  \label{fig:EX_4_Force}
\end{figure}

In \cref{tab:results4} we see that 24 iterations were
required to reach the final residual value \num{1.89e-10}.
The propagating crack, shown in \cref{fig:EX_4_PF},
is very similar to the results from \cite{Mang_2020}.
The corresponding optimal control is presented in \cref{fig:EX_4_Force}.
It decreases almost linearly, approximately from 2400 to 1500,
which is plausible since this experiment
has similarities to \cref{SEC:Experiment1}.
Similar to \cref{SEC:Experiment6}
we notice a small crack propagation
starting from the lower left corner $(0.5,0)$.
This is due to the singularity caused by
the Dirichlet condition on $\Gamma_D$ in combination with
the Neumann condition on $\set{0.5} \times [0,0.5]$.

\subsection{Experiment 6: inhibiting horizontal crack growth}
\label{SEC:Experiment5}

In our final experiment we expose the domain $\Omega = (0,1)^2$
to a time-independent external force $q_c$
which creates a growing crack
(for the tiny initial control $q = 1$).
Then we seek an optimal control $q$ that counteracts
the external force $q_c$ to inhibit the crack growth.
We choose the same partitioning of $\partial \Omega$
and the same notch as in \cref{SEC:Experiment1}.
The initial phase-field is
\begin{equation*}
  \label{EQ:EX_5_phi0}
  \varphi_0(x, y) \coloneqq
  \begin{cases}
    0, & x \in (0.5,1)\text{ and } y=0.5 \\
    1, & \text{else}.
  \end{cases}
\end{equation*}
We define the external force as a linear function: $q_c(x) = 850 + 1800 x$.
The time interval is $[0, 1]$ with 101 equidistant time points,
i.e., $T = 1$ and $M = 100$.
The spatial mesh consists of $64 \times 64$ square elements
with diameter $h = \sqrt2/64$.
The desired phase-field $\varphi_d$ has the value one on the whole domain.
Our findings for the tolerance \num{2.0e-11} are presented in
\cref{tab:results5,fig:EX_5_PF,fig:EX_5_Force,fig:EX_5_Forces_1}.
\begin{figure}
  \centering
  \begin{tikzpicture}
    \draw[fill=black!2] (0,0)
    -- (0,4) node[pos=0.5,left]  {$\Gamma\free$}
    -- (4,4) node[pos=0.5,above] {$\Gamma_N$}
    -- (4,0) node[pos=0.5,right] {$\Gamma\free$}
    -- cycle node[pos=0.5,below] {$\Gamma_D$};
    \node at (2,3) {$\Omega$};
    \node[blue,above] at (1,4) {$\uparrow q$};
    \node[blue,above] at (3,4) {$\uparrow q$};
    \node[red,below] at (1,4) {$\,\uparrow q_c$};
    \node[red,below] at (3,4) {$\,\uparrow q_c$};
    \draw[dotted,red] (2-1,2) -- (3-1,2) node[pos=0.5,above] {$\neg\varphi_d$};
    \draw[thick,blue] (3-1,2) -- (5-1,2) node[pos=0.5,above] {notch};
    \fill
    (2-1,2) circle (1.5pt) (3-1,2) circle (1.5pt) (5-1,2) circle (1.5pt);
  \end{tikzpicture}
  \caption{Experiment 6:
    domain $\Omega = (0,1)^2$ with partitioned boundary, intial notch,
    undesired crack $\neg \varphi_d$ and constant pulling force $q_c$.}
\end{figure}
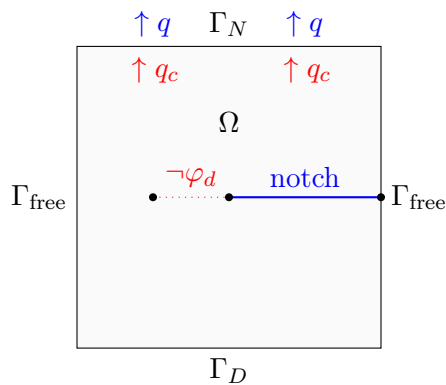
\begin{table}[tp]
  \centering
  \captionperf{6}
  \sisetup{table-format=1.4e+1}
  \begin{tabular}{rr
    S[table-format=1.2e+1]S[table-format=1.2e+2]SSSS[table-format=4.2]}
    \toprule
    Iter&CG&{Relative}&{Absolute}&{Cost}&{Tracking}&{Tikhonov}&{Force}\\[-2pt]
        &  &{residual}&{residual}&      &          &       &          \\
    \midrule
    0 & -- & 1.0     & 3.23e-06 & 1.7274e-3 & 1.4082e-3 & 3.1920e-4  & 1.0    \\
    1 &  3 & 0.558   & 1.80e-06 & 1.5649e-3 & 2.1034e-5 & 1.5438e-3  & 4262.61\\
    2 &  3 & 0.166   & 5.38e-07 & 1.9198e-4 & 1.8984e-4 & 2.1348e-6  & 824.97\\
    3 &  3 & 0.137   & 4.44e-07 & 1.2174e-4 & 7.6886e-6 & 1.1405e-4  & 1731.04\\
    4 &  3 & 5.26e-2 & 1.70e-07 & 9.3561e-5 & 9.0843e-5 & 2.7177e-6  & 956.04\\
    5 &  3 & 4.40e-2 & 1.42e-07 & 6.5026e-5 & 3.4897e-5 & 3.0129e-5  & 1280.33\\
    6 &  3 & 1.25e-3 & 4.05e-09 & 7.2949e-5 & 6.3052e-5 & 9.8963e-6  & 1079.64\\
    7 & 12 & 1.02e-3 & 3.31e-09 & 6.7607e-5 & 5.2538e-5 & 1.5069e-5  & 1142.82\\
    8 & 12 & 3.10e-4 & 1.00e-09 & 7.0443e-5 & 5.8585e-5 & 1.1858e-5  & 1104.25\\
    9 & 10 & 2.52e-4 & 8.16e-10 & 6.9141e-5 & 5.5993e-5 & 1.3148e-5  & 1120.23\\
    10&  9 & 7.94e-5 & 2.57e-10 & 6.9890e-5 & 5.7514e-5 & 1.2377e-5  & 1110.86\\
    11&  7 & 6.42e-5 & 2.07e-10 & 6.9560e-5 & 5.6853e-5 & 1.2707e-5  & 1114.79\\
    12&  7 & 2.03e-5 & 6.57e-11 & 6.9753e-5 & 5.7241e-5 & 1.2512e-5  & 1112.52\\
    13&  5 & 1.63e-5 & 5.28e-11 & 6.9669e-5 & 5.7073e-5 & 1.2596e-5  & 1113.45\\
    14&  4 & 5.26e-6 & 1.70e-11 & 6.9718e-5 & 5.7172e-5 & 1.2546e-5  & 1112.98\\
    15&  3 & 4.21e-6 & 1.36e-11 & 6.9697e-5 & 5.7128e-5 & 1.2568e-5  & 1113.14\\
    \bottomrule
  \end{tabular}
  \label{tab:results5}
\end{table}
\begin{table}
  \captionpar{6}
  \centering
  \begin{tabular}{llS[table-format=1.1e+2]}
    \toprule
    Par. & Definition & {Value} \\
    \midrule
    $\varepsilon$ & Regul. (crack) $\approx 2 h$ & 0.0442 \\
    $\kappa$ & Regul. (crack) & 1.0e-10 \\
    $\eta$ & Regul. (viscosity) & 1.0e3 \\
    $\gamma$ & Penalty & 1.0e5  \\
    $\alpha$ & Tikhonov & 1.0e-9 \\
    \bottomrule
  \end{tabular}
  \hfil
  \begin{tabular}{llS[table-format=+1.1e1]}
    \toprule
    Par. & Definition & {Value} \\
    \midrule
    $G_c$ & Fracture toughness & 1.0 \\
    $\nu_s$ & Poisson's ratio & 0.2 \\
    $E$ & Young's modulus & 1.0e6 \\
    $q_0$ & Initial control & 1.0 \\
    $q_d$ & Nominal control & -8.0e2 \\
    \bottomrule
  \end{tabular}
  \label{tab:EX_5_params}
\end{table}%
\renewcommand\img[1]{\image
  [trim=840 315 840 315,clip=true,width=.33\linewidth]
  {EX_5_STATE_#1_TIME_100_PF.png}}
\begin{figure}[tp]
  \footnotesize
  \centering
  \img{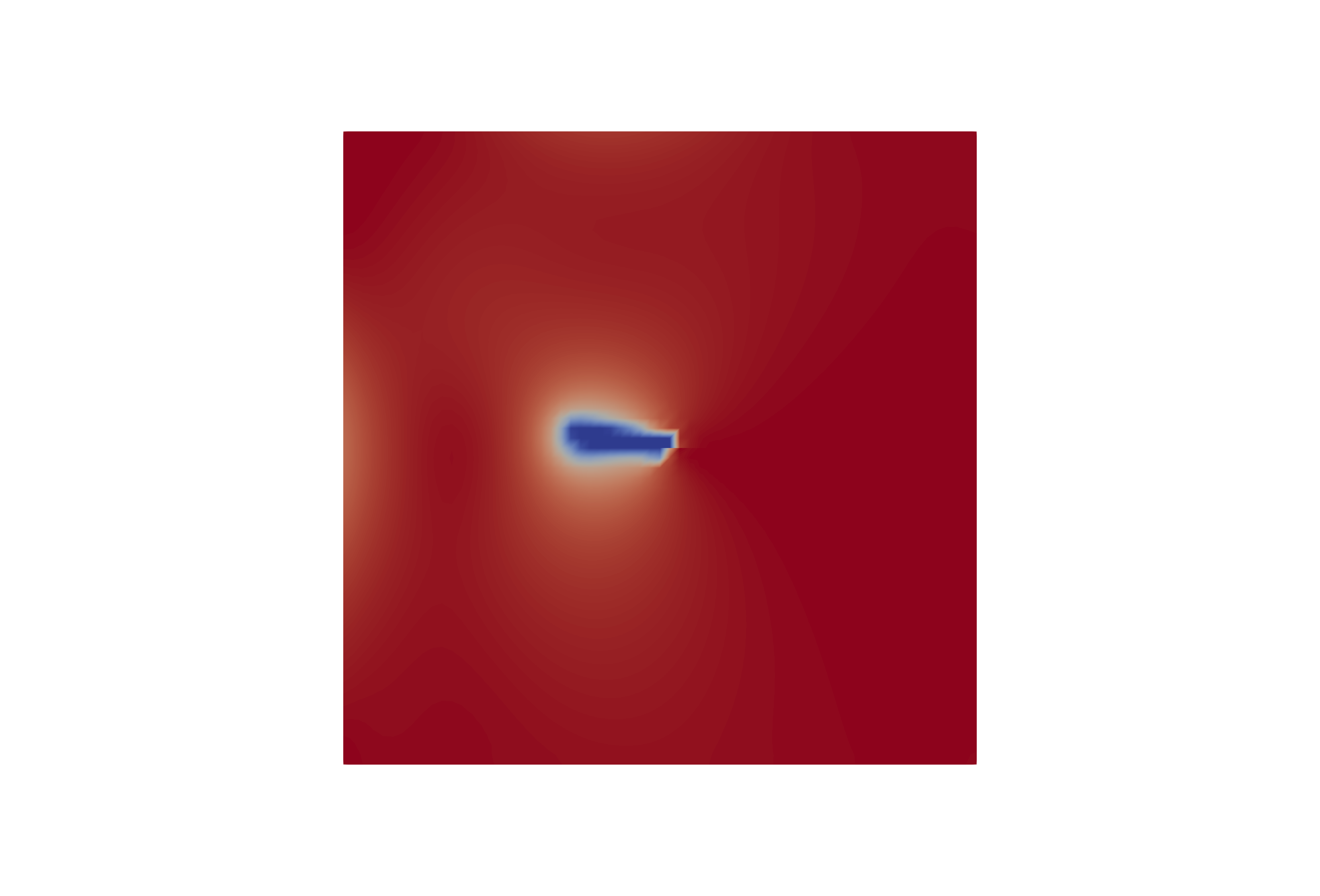}\enspace\img{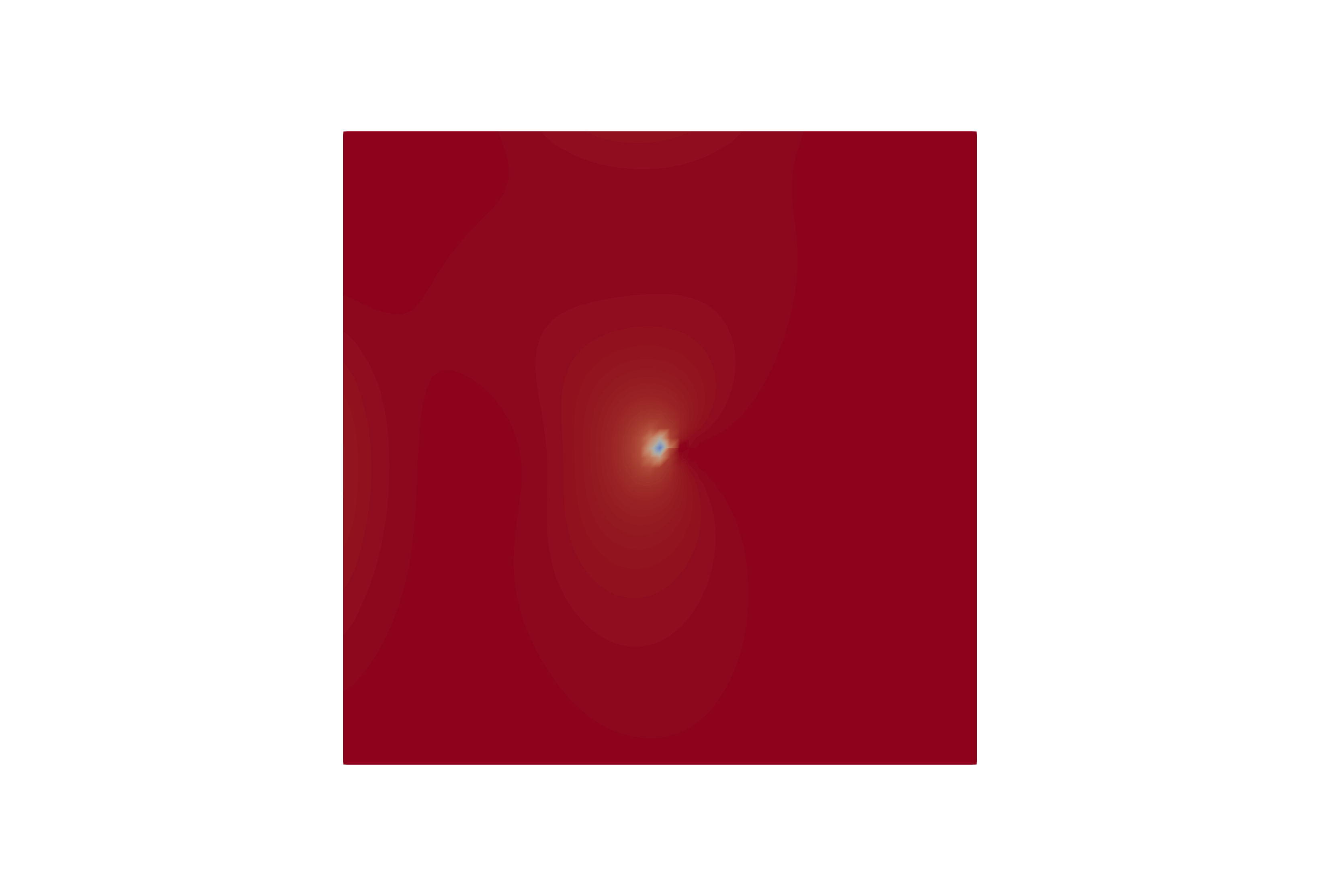}\\
  \scale[0.33]{0.32}{0}{0.5}{1}%
  \caption{Experiment 6:
    initial phase-field $\varphi$ (left, iteration 0) and
    optimal phase-field (right, iteration 15) at final time 100.}
  \label{fig:EX_5_PF}
\end{figure}
\definecolor{myForestGreen}{rgb}{0.133, 0.545, 0.133} 
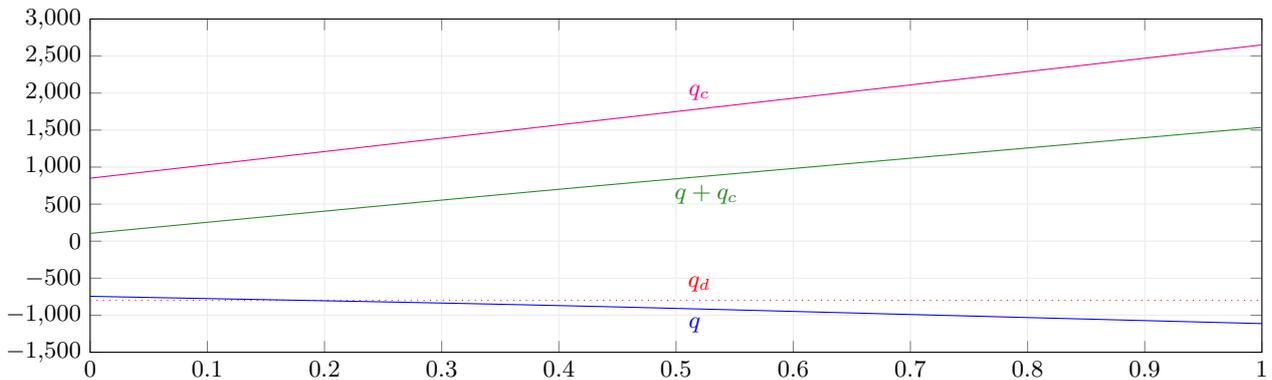
\begin{figure}
  \footnotesize
  \centering
  \begin{tikzpicture}
    \begin{axis}[width=1.0\linewidth,height=60mm,
      xmin=0,xmax=1,ymin=-1500,ymax=3000,ytick distance={500},grid=major]
      \draw[red,dotted] (0,-800) -- (1,-800) node[pos=0.52,above] {$q_d$};
      \draw[magenta] (0,850) -- (1,2650) node[pos=0.52,above] {$q_c$};
      \addplot[blue] table
      {TXT5_EX_5_Force.txt} node[pos=0.46,below] {$q$};
      \addplot[myForestGreen] table
      {TXT5_EX_5_Force_plus_q_c.txt} node[pos=0.54,below] {$q+q_c$};
    \end{axis}
  \end{tikzpicture}
  \caption{Experiment 6: optimal control force (blue),
    nominal control (red, dotted), constant control (magenta)
    and resulting total control $q + q_c$ (green) on upper boundary $\Gamma_N$.}
  \label{fig:EX_5_Force}
\end{figure}

In \cref{fig:EX_5_PF} we observe that no crack propagation
occurs with the computed optimal control.
As a result, the desired phase-field is successfully reproduced
and the initial value of the cost functional is reduced by $96\%$.
In comparison to all other experiments,
where we achieved a maximum reduction of $70\%$,
this is a remarkable result.
Although the sum of the optimal control $q$ and
the constant external force $q_c$ is positive everywhere,
see \cref{fig:EX_5_Force},
it is not large enough to create a propagating crack.
This is to be expected since the Tikhonov term
would penalize an unnecessarily strong control force.
In \cref{fig:EX_5_Forces_1} the control forces
of iterations 1 to~4 from \cref{tab:results5} are shown.
We observe that the first control $q_1$
decreases almost linearly to a minimal value of $-4262.61$,
which produces a relatively large Tikhonov term.
The second control $q_2$ is instead
a linearly increasing function that minimzes this term.
The third and fourth controls lie between $q_1$ and $q_2$,
and $q_4$ behaves similar to the optimal control.
When studying the value of the tracking part in \cref{tab:results5},
we observe that it becomes almost zero on the third iteration
where the Tikhonov term is more than 10 times larger.
Consequently the Tikhonov term must be reduced next.
Subsequently, the two terms oscillate until they are roughly balanced.
Some of the former experiments have been sensitive
to the choice of $\alpha$,
but none of them has been as sensitive as this experiment.

\begin{figure}
  \footnotesize
  \centering
  \begin{tikzpicture}
    \begin{axis}[width=1.0\linewidth,height=10cm,
      xmin=0,xmax=1,ymin=-4500,ymax=100,ytick distance={500},grid=major]
      \draw[red,dotted] (0,-800) -- (1,-800) node[pos=0.2,below] {$q_d$};
      \addplot[blue] table
      {TXT5_EX_5_Force_1.txt} node[pos=0.6,below] {$q_1$};
      \addplot[myForestGreen] table
      {TXT5_EX_5_Force_2.txt} node[pos=0.47,above] {$q_2$};
      \addplot[teal] table
      {TXT5_EX_5_Force_3.txt} node[pos=0.59,below] {$q_3$};
      \addplot[magenta] table
      {TXT5_EX_5_Force_4.txt} node[pos=0.555,below] {$q_4$};
    \end{axis}
  \end{tikzpicture}
  \caption{Experiment 6: control forces for iterations 1--4 (solid)
    and nominal control (dotted) on upper boundary $\Gamma_N$.}
  \label{fig:EX_5_Forces_1}
  \label{fig:EX_5_Forces_2}
  \label{fig:EX_5_Forces_3}
  \label{fig:EX_5_Forces_4}
\end{figure}
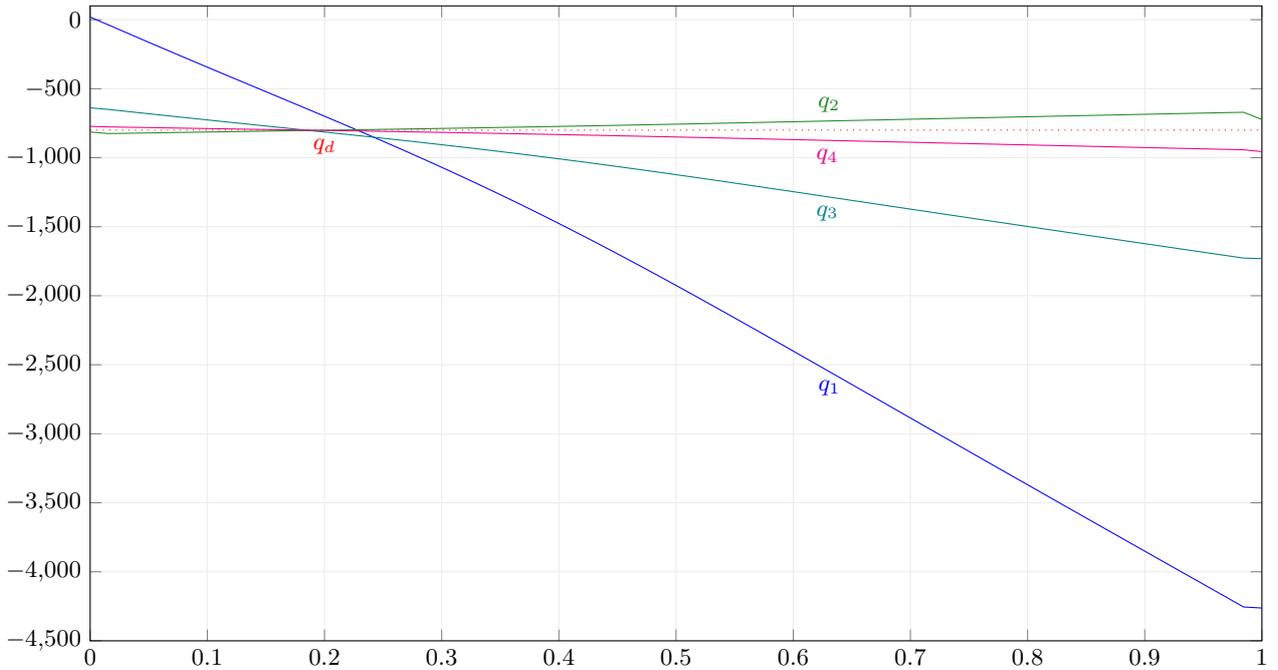

\section{Conclusions}
\label{sec_conclusions}

In this paper we derived a space-time Galerkin
formulation for a regularized phase-field fracture optimal control setting.
By introducing jump terms in time and with the help of a discontinuous
Galerkin discretization in time, specific time-stepping schemes were
obtained. A careful investigation of correct weighting
of two regularization terms and the initial conditions was necessary
for the forward phase-field fracture problem. The solution process
of the optimization problem was based on the reduced approach
in which the state variables are obtained from a solution operator
acting on the controls. The numerical solution algorithm is based
on Newton's method in which three auxiliary problems are required.
The first part of the paper was concerned with the detailed derivation
of these terms, which are to the best of our knowledge novel
in the published literature.

In Section \ref{SEC:Numtests}, we performed
several detailed computational performance studies
for space-time phase-field fracture optimal control problems.
The optimization problem was designed with the help of a reduced approach
in which the state variables are obtained from a solution operator
acting on the controls. Therein, a monolithic space-time respresentation
of the phase-field fracture problem was adopted. Moreover, the crack
irreversibility constraint was regularized using a penality approach.
To study the performance, we investigated six numerical experiments
with single (Experiments 1, 2, 5, 6) and multiple fractures
(Experiments 3, 4), single controls (Experiments 1, 3, 4, 5, 6) and two controls
(Experiments~2), propagating fractures (Experiments 1, 2, 3, 4, 5) and
inhibiting crack growth (Experiment~6). Therein, the performance
of the NLP solver (Algorithm~1) and the inner CG method as well as
the phase-field fracture PDE constraint were computationally analyzed
in great detail.
One main bottleneck is the computational cost
of the inner linear solver of the forward problem,
which is well-known and analogous in other PDE-constrained optimization problems.
In ongoing work, we plan to incorporate
parallel adaptive preconditioned iterative solvers 
\cite{JoLaWi20,JoLaWi20_parallel}, which, however,
is a major extension and was out of scope in this work.

\section{Acknowledgements}

The first and third author are partially funded by the
Deutsche Forschungsgemeinschaft (DFG, German Research Foundation)
Priority Program 1962 (DFG SPP 1962) within the subproject
\emph{Optimizing Fracture Propagation using a Phase-Field Approach}
with the project number 314067056.
The second author is funded by the DFG -- SFB1463 -- 434502799.

\bibliographystyle{abbrv}
\bibliography{lit}

\begin{thebibliography}{10}

\bibitem{Allaire20115010}
G.~Allaire, F.~Jouve, and N.~V. Goethem.
\newblock Damage and fracture evolution in brittle materials by shape
  optimization methods.
\newblock {\em Journal of Computational Physics}, 230(12):5010 -- 5044, 2011.

\bibitem{AmGeraLoren15}
M.~Ambati, T.~Gerasimov, and L.~De~Lorenzis.
\newblock A review on phase-field models of brittle fracture and a new fast
  hybrid formulation.
\newblock {\em Computational Mechanics}, 55(2):383--405, 2015.

\bibitem{AmTo90}
L.~Ambrosio and V.~Tortorelli.
\newblock Approximation of functionals depending on jumps by elliptic
  functionals via $\gamma$-convergence.
\newblock {\em Comm. Pure Appl. Math.}, 43:999--1036, 1990.

\bibitem{AmTo92}
L.~Ambrosio and V.~Tortorelli.
\newblock On the approximation of free discontinuity problems.
\newblock {\em Boll. Un. Mat. Ital. B}, 6:105--123, 1992.

\bibitem{dealII91}
D.~Arndt, W.~Bangerth, T.~C. Clevenger, D.~Davydov, M.~Fehling,
  D.~Garcia-Sanchez, G.~Harper, T.~Heister, L.~Heltai, M.~Kronbichler, R.~M.
  Kynch, M.~Maier, J.-P. Pelteret, B.~Turcksin, and D.~Wells.
\newblock The \texttt{deal.II} library, version 9.1.
\newblock {\em Journal of Numerical Mathematics}, 2019.

\bibitem{deal2020}
D.~Arndt, W.~Bangerth, D.~Davydov, T.~Heister, L.~Heltai, M.~Kronbichler,
  M.~Maier, J.-P. Pelteret, B.~Turcksin, and D.~Wells.
\newblock The \texttt{deal.II} finite element library: Design, features, and
  insights.
\newblock {\em Computers \& Mathematics with Applications}, 2020.

\bibitem{Barbu1984}
V.~P. Barbu.
\newblock {\em Optimal Control of Variational Inequalities}, volume 100.
\newblock Pitman Advanced Pub. Program, 1984.

\bibitem{BeMeVe07}
R.~Becker, D.~Meidner, and B.~Vexler.
\newblock Efficient numerical solution of parabolic optimization problems by
  finite element methods.
\newblock {\em Optim. Methods Softw.}, 22(5):813--833, 2007.

\bibitem{BeCoOhlWill15}
P.~Benner, A.~Cohen, M.~Ohlberger, and K.~Willcox.
\newblock {\em Model Reduction and Approximation: Theory and Algorithms}.
\newblock SIAM Philadelphia, 2015.

\bibitem{BoVeScoHuLa12}
M.~J. Borden, C.~V. Verhoosel, M.~A. Scott, T.~J.~R. Hughes, and C.~M. Landis.
\newblock A phase-field description of dynamic brittle fracture.
\newblock {\em Comput. Meth. Appl. Mech. Engrg.}, 217:77--95, 2012.

\bibitem{Bou99}
B.~Bourdin.
\newblock Image segmentation with a finite element method.
\newblock {\em Mathematical Modelling and Numerical Analysis}, 33(2):229--244,
  1999.

\bibitem{Bour07}
B.~Bourdin.
\newblock Numerical implementation of the variational formulation for
  quasi-static brittle fracture.
\newblock {\em Interfaces and free boundaries}, 9:411--430, 2007.

\bibitem{BourFraMar00}
B.~Bourdin, G.~Francfort, and J.-J. Marigo.
\newblock Numerical experiments in revisited brittle fracture.
\newblock {\em J. Mech. Phys. Solids}, 48(4):797--826, 2000.

\bibitem{BourFraMar08}
B.~Bourdin, G.~Francfort, and J.-J. Marigo.
\newblock The variational approach to fracture.
\newblock {\em J. Elasticity}, 91(1--3):1--148, 2008.

\bibitem{BouFra19}
B.~Bourdin and G.~A. Francfort.
\newblock Past and present of variational fracture.
\newblock {\em SIAM News}, 52(9), 2019.

\bibitem{BouLarRi11}
B.~Bourdin, C.~Larsen, and C.~Richardson.
\newblock A time-discrete model for dynamic fracture based on crack
  regularization.
\newblock {\em Int. J. Frac.}, 168(2):133--143, 2011.

\bibitem{Braides1998}
A.~Braides.
\newblock {\em Approximation of free-discontinuity problems}.
\newblock Springer Berlin Heidelberg, 1998.

\bibitem{BrWiBeNoRa20}
M.~K. Brun, T.~Wick, I.~Berre, J.~M. Nordbotten, and F.~A. Radu.
\newblock An iterative staggered scheme for phase field brittle fracture
  propagation with stabilizing parameters.
\newblock {\em Computer Methods in Applied Mechanics and Engineering},
  361:112752, 2020.

\bibitem{BuOrSue10}
S.~Burke, C.~Ortner, and E.~S\"uli.
\newblock An adaptive finite element approximation of a variational model of
  brittle fracture.
\newblock {\em SIAM J. Numer. Anal.}, 48(3):980--1012, 2010.

\bibitem{Cia87}
P.~G. Ciarlet.
\newblock {\em The Finite Element Method for Elliptic Problems}.
\newblock North-Holland, Amsterdam [u.a.], 2. pr. edition, 1987.

\bibitem{DESAI2022111048}
J.~Desai, G.~Allaire, and F.~Jouve.
\newblock Topology optimization of structures undergoing brittle fracture.
\newblock {\em Journal of Computational Physics}, 458:111048, 2022.

\bibitem{DiLiWiTy22}
P.~Diehl, R.~Lipton, T.~Wick, and M.~Tyagi.
\newblock {A comparative review of peridynamics and phase-field models for
  engineering fracture mechanics}.
\newblock {\em Computational Mechanics}, pages 1--35, 2022.

\bibitem{dope}
The {D}ifferential {E}quation and {O}ptimization {E}nvironment:
  \textsc{DOpElib}.

\bibitem{Fra21}
G.~Francfort.
\newblock Variational fracture: Twenty years after.
\newblock {\em International Journal of Fracture}, pages 1--11, 2021.

\bibitem{FraMar98}
G.~Francfort and J.-J. Marigo.
\newblock Revisiting brittle fracture as an energy minimization problem.
\newblock {\em J. Mech. Phys. Solids}, 46(8):1319--1342, 1998.

\bibitem{GeLo16}
T.~Gerasimov and L.~D. Lorenzis.
\newblock A line search assisted monolithic approach for phase-field computing
  of brittle fracture.
\newblock {\em Computer Methods in Applied Mechanics and Engineering},
  312:276--303, 2016.

\bibitem{GERASIMOV2020113353}
T.~Gerasimov, U.~Römer, J.~Vondřejc, H.~G. Matthies, and L.~{De Lorenzis}.
\newblock Stochastic phase-field modeling of brittle fracture: Computing
  multiple crack patterns and their probabilities.
\newblock {\em Computer Methods in Applied Mechanics and Engineering},
  372:113353, 2020.

\bibitem{DOpElib}
C.~Goll, T.~Wick, and W.~Wollner.
\newblock {DOpElib}: {D}ifferential equations and optimization environment; {A}
  goal oriented software library for solving pdes and optimization problems
  with pdes.
\newblock {\em Archive of Numerical Software}, 5(2):1--14, 2017.

\bibitem{HaKa09}
V.~Hakim and A.~Karma.
\newblock Laws of crack motion and phase-field models of fracture.
\newblock {\em J. Mech. Phys. Solids}, 57(2):342--368, 2009.

\bibitem{HeWheWi15}
T.~Heister, M.~F. Wheeler, and T.~Wick.
\newblock A primal-dual active set method and predictor-corrector mesh
  adaptivity for computing fracture propagation using a phase-field approach.
\newblock {\em Comp. Meth. Appl. Mech. Engrg.}, 290:466--495, 2015.

\bibitem{HiPiUlUl09}
M.~Hinze, R.~Pinnau, M.~Ulbrich, and S.~Ulbrich.
\newblock {\em Optimization with {PDE} Constraints}.
\newblock Number~23 in Mathematical modelling: theory and applications.
  Springer, Dordrecht u.a., 2009.

\bibitem{JoLaWi20}
D.~Jodlbauer, U.~Langer, and T.~Wick.
\newblock Matrix-free multigrid solvers for phase-field fracture problems.
\newblock {\em Computer Methods in Applied Mechanics and Engineering},
  372:113431, 2020.

\bibitem{JoLaWi20_parallel}
D.~Jodlbauer, U.~Langer, and T.~Wick.
\newblock Parallel matrix-free higher-order finite element solvers for
  phase-field fracture problems.
\newblock {\em Mathematical and Computational Applications}, 25(3):40, 2020.

\bibitem{KaKeLe01}
A.~Karma, D.~Kessler, and H.~Levine.
\newblock Phase-field model of mode iii dynamic fracture.
\newblock {\em Physical Review Letters}, 87(4):45501 618, 2001.

\bibitem{KhiSteiWi21}
D.~Khimin, M.~C. Steinbach, and T.~Wick.
\newblock Optimal control for phase-field fracture: Algorithmic concepts and
  computations.
\newblock In F.~Aldakheel, B.~Hudobivnik, M.~Soleimani, H.~Wessels,
  C.~Weißenfels, and M.~Marino, editors, {\em Current Trends and Open Problems
  in Computational Mechanics}. Springer, 2022.

\bibitem{KhoNoPaAbWiHei20}
A.~Khodadadian, N.~Noii, M.~Parvizi, M.~Abbaszadeh, T.~Wick, and C.~Heitzinger.
\newblock A {B}ayesian estimation method for variational phase-field fracture
  problems.
\newblock {\em Computational Mechanics}, 66:827--849, 2020.

\bibitem{KneRoZa13}
D.~Knees, R.~Rossi, and C.~Zanini.
\newblock A vanishing viscosity approach to a rate-independent damage model.
\newblock {\em Mathematical Models and Methods in Applied Sciences},
  23(04):565--616, 2013.

\bibitem{KOLDITZ2022100047}
L.~Kolditz and K.~Mang.
\newblock On the relation of gamma-convergence parameters for pressure-driven
  quasi-static phase-field fracture.
\newblock {\em Examples and Counterexamples}, 2:100047, 2022.

\bibitem{KoKr20}
A.~Kopanicakova and R.~Krause.
\newblock A recursive multilevel trust region method with application to fully
  monolithic phase-field models of brittle fracture.
\newblock {\em Computer Methods in Applied Mechanics and Engineering},
  360:112720, 2020.

\bibitem{KuMue10}
C.~Kuhn and R.~M{\"u}ller.
\newblock A continuum phase field model for fracture.
\newblock {\em Engineering Fracture Mechanics}, 77(18):3625 -- 3634, 2010.

\bibitem{Mang_2020}
K.~Mang, T.~Wick, and W.~Wollner.
\newblock A phase-field model for fractures in nearly incompressible solids.
\newblock {\em Computational Mechanics}, 65(1):61–78, 2020.

\bibitem{Me08}
D.~Meidner.
\newblock {\em Adaptive Space-Time Finite Element Methods for Optimization
  Problems Governed by Nonlinear Parabolic Systems}.
\newblock PhD thesis, University of Heidelberg, 2008.

\bibitem{MesBouKhon15}
A.~Mesgarnejad, B.~Bourdin, and M.~Khonsari.
\newblock Validation simulations for the variational approach to fracture.
\newblock {\em Computer Methods in Applied Mechanics and Engineering}, 290:420
  -- 437, 2015.

\bibitem{MieWelHof10b}
C.~Miehe, M.~Hofacker, and F.~Welschinger.
\newblock A phase field model for rate-independent crack propagation: {R}obust
  algorithmic implementation based on operator splits.
\newblock {\em Comput. Meth. Appl. Mech. Engrg.}, 199:2765--2778, 2010.

\bibitem{MieWelHof10a}
C.~Miehe, F.~Welschinger, and M.~Hofacker.
\newblock Thermodynamically consistent phase-field models of fracture:
  variational principles and multi-field fe implementations.
\newblock {\em Int. J. Numer. Methods Engrg.}, 83:1273--1311, 2010.

\bibitem{MIGNOT1976}
F.~Mignot.
\newblock Contrôle dans les inéquations variationelles elliptiques.
\newblock {\em Journal of Functional Analysis}, 22(2):130--185, 1976.

\bibitem{MIGNOT1984}
F.~Mignot and J.~P. Puel.
\newblock Optimal control in some variational inequalities.
\newblock 22(3):466--476, May 1984.

\bibitem{MiWheWi13a}
A.~Mikeli\'c, M.~Wheeler, and T.~Wick.
\newblock A phase-field approach to the fluid filled fracture surrounded by a
  poroelastic medium.
\newblock ICES Report 13-15, Jun 2013.

\bibitem{MiWheWi19}
A.~Mikeli{\'{c}}, M.~F. Wheeler, and T.~Wick.
\newblock Phase-field modeling through iterative splitting of hydraulic
  fractures in a poroelastic medium.
\newblock {\em GEM - International Journal on Geomathematics}, 10(1), Jan 2019.

\bibitem{MoWo20}
M.~Mohammadi and W.~Wollner.
\newblock Phase field modelling of fracture.
\newblock {\em Optimization and Engineering}, 2020.

\bibitem{NeiWiWo17}
I.~Neitzel, T.~Wick, and W.~Wollner.
\newblock An optimal control problem governed by a regularized phase-field
  fracture propagation model.
\newblock {\em SIAM Journal on Control and Optimization}, 55(4):2271--2288,
  2017.

\bibitem{NeiWiWo19}
I.~Neitzel, T.~Wick, and W.~Wollner.
\newblock An optimal control problem governed by a regularized phase-field
  fracture propagation model. {P}art {II}: The regularization limit.
\newblock {\em SIAM Journal on Control and Optimization}, 57(3):1672--1690,
  2019.

\bibitem{NoKhoUlAlWiFrWr21}
N.~Noii, A.~Khodadadian, J.~Ulloa, F.~Aldakheel, T.~Wick, S.~Francois, and
  P.~Wriggers.
\newblock Bayesian inversion for unified ductile phase-field fracture.
\newblock {\em Computational Mechanics}, 2021.

\bibitem{NoKhoWi21}
N.~Noii, A.~Khodadadian, and T.~Wick.
\newblock Bayesian inversion for anisotropic hydraulic phase-field fracture.
\newblock {\em Computer Methods in Applied Mechanics and Engineering},
  386:114--118, 2021.

\bibitem{Robinson1976}
S.~M. Robinson.
\newblock Stability theory for systems of inequalities, part ii: Differentiable
  nonlinear systems.
\newblock {\em SIAM Journal on Numerical Analysis}, 13(4):497--513, 1976.

\bibitem{SpBrKa11}
R.~Spatschek, E.~Brener, and A.~Karma.
\newblock Phase field modeling of crack propagation.
\newblock {\em Philosophical Magazine}, 91(1):75--95, 2011.

\bibitem{Troe09}
F.~Tr{\"o}ltzsch.
\newblock {\em {O}ptimale {S}teuerung partieller {D}ifferentialgleichungen -
  {T}heorie, {V}erfahren und {A}nwendungen}.
\newblock Vieweg und Teubner, Wiesbaden, 2nd edition, 2009.

\bibitem{wambacq2020interiorpoint}
J.~Wambacq, J.~Ulloa, G.~Lombaert, and S.~François.
\newblock Interior-point methods for the phase-field approach to brittle and
  ductile fracture, 2020.

\bibitem{Wi17_SISC}
T.~Wick.
\newblock An error-oriented {N}ewton/inexact augmented {L}agrangian approach
  for fully monolithic phase-field fracture propagation.
\newblock {\em SIAM Journal on Scientific Computing}, 39(4):B589--B617, 2017.

\bibitem{Wi17_CMAME}
T.~Wick.
\newblock Modified {N}ewton methods for solving fully monolithic phase-field
  quasi-static brittle fracture propagation.
\newblock {\em Computer Methods in Applied Mechanics and Engineering},
  325:577--611, 2017.

\bibitem{Wi20_book}
T.~Wick.
\newblock {\em Multiphysics Phase-Field Fracture: Modeling, Adaptive
  Discretizations, and Solvers}.
\newblock De Gruyter, Berlin, Boston, 2020.

\bibitem{winkler2001}
B.~J. Winkler.
\newblock {\em Traglastuntersuchungen von unbewehrten und bewehrten
  Betonstrukturen auf der Grundlage eines objektiven Werkstoffgesetzes für
  Beton}.
\newblock PhD thesis, Uni, 2001.

\bibitem{WuNgNgSuBoSi19}
J.-Y. Wu, V.~P. Nguyen, C.~Thanh~Nguyen, D.~Sutula, S.~Bordas, and S.~Sinaie.
\newblock Phase field modelling of fracture.
\newblock {\em Advances in Applied Mechanics}, 53:1--183, 09 2020.

\bibitem{WuRoLoMa21}
T.~Wu, B.~Rosic, L.~de~Lorenzis, and H.~Matthies.
\newblock Parameter identification for phase-field modeling of fracture: a
  {B}ayesian approach with sampling-free update.
\newblock {\em Computational Mechanics}, 67:435--453, 2021.

\bibitem{Zehnder2012}
A.~Zehnder.
\newblock {\em Fracture mechanics}.
\newblock Springer-Verlag, 2012.

\bibitem{Zowe1979}
J.~Zowe and S.~Kurcyusz.
\newblock Regularity and stability for the mathematical programming problem in
  banach spaces.
\newblock {\em Applied Mathematics and Optimization}, 5(1):49--62, Mar. 1979.

\end{thebibliography}

\end{document}